%% file: lie-theoretic-chains.tex
\documentclass[12pt,twoside,leqno,openany]{amsart}

\pdfoutput=1

\usepackage{amsbsy,amscd,amsfonts,amsmath,amssymb,amsthm,color,
fancybox,fancyhdr,footmisc,graphics,graphicx,ifthen,mathrsfs,
multicol,pdfpages,rotating,times,wasysym,pifont,blkarray}
\usepackage[dvipsnames,svgnames,x11names]{xcolor}

\usepackage[all]{xy}
\usepackage[utf8]{inputenc}
\usepackage[T1]{fontenc}
\usepackage{mathtools}
\newtagform{EngelLie}[\scriptstyle]{$}{$}
\makeatletter\newcommand{\leqnomode}{\tagsleft@true}
\newcommand{\reqnomode}{\tagsleft@false}\makeatother

\input macros.tex


\begin{document}

\setcounter{section}{0}

$\:$

\bigskip\bigskip\bigskip\bigskip\bigskip

\begin{center}

{\large\bf A Lie-Theoretic Construction\footnotemark[1] 
of Cartan-Moser Chains}
\label{Lie-theoretic-Moser-chains}

\footnotetext[1]{\,
This work was supported
in part by the Polish National Science Centre (NCN) 
via the grant number 2018/29/B/ST1/02583.}

\medskip

\bigskip\bigskip
\bigskip\bigskip

Joël~{\sc Merker}

\end{center}\bigskip

\begin{center}
\begin{minipage}[t]{12.5cm}
\parindent 0.53cm
\footnotesize
\noindent
{\sc Abstract}.
Let $M^3 \subset \mathbb{C}^2$ 
be a $\mathcal{C}^\omega$ Levi nondegenerate
hypersurface. In the literature, Cartan-Moser chains are detected from
rather advanced considerations: either from the construction of a
Cartan connection associated with the CR equivalence problem; or from
the construction of a formal or converging Poincar\'e-Moser normal
form.

This note provides an alternative direct elementary construction,
based on the inspection of the Lie prolongations of $5$ infinitesimal
holomorphic automorphisms to the space of second order jets of
CR-transversal curves. Within the $4$-dimensional jet fiber, the
orbits of these $5$ prolonged fields happen to have a simple cubic
$2$-dimensional degenerate exceptional orbit, the {\sl chain locus:}
\[
\Sigma_0
\,:=\,
\big\{
(x_1,y_1,x_2,y_2)
\in
\mathbb{R}^4
\colon\,\,
x_2
=
-2x_1^2y_1-2y_1^3,\,\,\,
y_2
=
2x_1y_1^2
+
2x_1^3
\big\}.
\]

Using plain translations, we may capture all points by working {\em
only at one point}, the origin, and computations become conceptually
enlightening and simple.

\end{minipage}
\end{center}

\Section{\bf Introduction}
\label{introduction-Lie-theoretic-Moser}
\HEAD{{\ref{introduction-Lie-theoretic-Moser}}.~{\sf Introduction}
}{
Joël {\sc Merker}}

The goal of this article is to present a simplified construction of
{\sl Cartan-Moser chains}, 
which are certain distinguished curves in Levi
nondegenerate Cauchy-Riemann (CR) manifolds of hypersurface type. 
We concentrate on
real-analytic embedded 
CR manifolds, because the interaction between the
{\em extrinsic} geometry of an ambient complex manifold $X$ and the
{\em intrinsic} geometry of a CR submanifold $M
\subset X$ is {\em richer} than in an abstract seetting.  Also, for
the sake of intuitive clarity and for elementariness, we restrict our
presentation to the $3$-dimensional case. The Lie-theoretical method
that we employ\,\,---\,\,which certainly has a wider
scope\,\,---\,\,drastically contracts all required computations by
working {\em only at one point}, as we shall rapidly see.

Thus, let $M^3 \subset \C^2$ be a $\mathcal{C}^\omega$ real
hypersurface. We are interested in results of a local nature, hence we
will allow to shrink neighborhoods of various points $p \in M$.  If $J
\colon T\C^2 \longrightarrow T\C^2$ is the standard complex structure,
with $J^2 = - \Id$, the complex tangent bundle $T^c M := TM \cap JTM$
is $J$-invariant of real rank $2$, hence at all point $p \in M$, the
$2$-planes $T_p^c M \subset T\C^2$ can be viewed as complex
affine sublines $\C \subset \C^2$.  Also, $T^{1,0}M := \big\{ X - i\,
JX\colon\, X \in T^cM \big\}$ and $T^{0,1}M := \big\{ X + i\,
JX\colon\, X \in T^cM \big\} = \overline{T^{1,0}}M$ are complex vector
subbundles of the complexified tangent bundle $\C \otimes_\R TM$.

We will always assume that $M^3 \subset
\C^2$ is {\sl Levi nondegenerate},
namely that $T^c M + [T^cM,\, T^cM] = TM$,
or equivalently~{\cite{Merker-Pocchiola-Sabzevari-2013-5-CR-II}}:
\[
\C\otimes_\R TM
\,=\,
T^{1,0}M
+
T^{0,1}M
+
[T^{1,0}M,\,T^{0,1}M].
\] 
For detailed foundations, the
reader may consult~{\cite{Merker-Pocchiola-Sabzevari-2013-5-CR-II}}.

These ``{\sl CR bundles}'' are invariant, in the sense that for any
(local) biholomorphism $h \colon \C^2 \longrightarrow {\C'}^2$ defined
in some neighborhood of $M$, with $M' := h(M)$ being a hypersurface of
${\C'}^2$, one has $h_\ast(T_p^cM) = T_{h(p)}^cM'$, and
$h_\ast(T_p^{1,0}M) = T_{h(p)}^{1,0}M'$ as well, where, by $h_\ast$,
we denote the differential of $h$ acting both on $TM$ and on $\C
\otimes_\R TM$, with the convention $\overline{h}_\ast = h_\ast$, {\em
  cf.}~{\cite{Merker-Pocchiola-Sabzevari-2013-5-CR-II}}.  Hence,
whenever $h$ is a (local) biholomorphism, $h\vert_M \colon M
\longrightarrow h(M)$ realizes a CR diffeomorphism.

So by definition, biholomorphic or CR equivalences stabilize some {\sl
horizontal} $2$-plane distribution $T^cM$, or the pair $T^{1,0} M
\oplus T^{0,1} M \subset \C TM$.  It seems that there is no reason
that there should exist some {\sl CR-transversal} structure which
would also be CR-invariant. For instance, 
does there exist a line field
$\{ \ell_p \}_{p\in M}$ with $\R \cong \ell_p \subset T_p M$ 
complementing $T_p^cM$ in
$T_pM = \ell_p \oplus T_p^c M$ 
which would be CR invariant\text{\bf ?} 
Yes of course in presence of some extra structure like 
{\em e.g.} a Riemannian metric on $M$\,\,---\,\,just 
take $\ell_p := [T_p^cM]^\bot$\,\,---, but {\em no} in general,
as is well known and 
as we will see.

\begin{center}
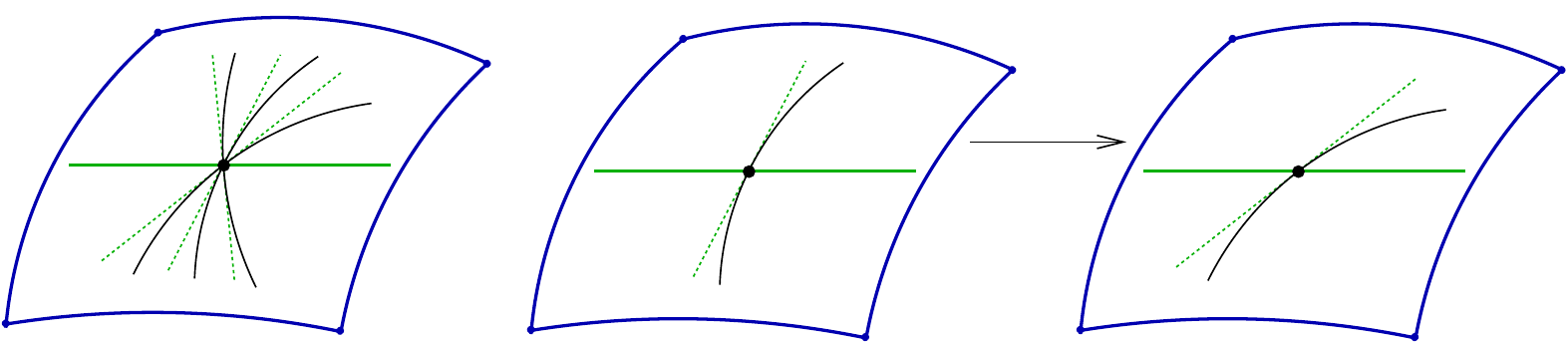

\nopagebreak
\begin{minipage}[t]{12.5cm}
\footnotesize
{\sc Figure~1:}
Left: representation of various chains at $p \in M$ directed by
various directions $\ell_p \subset T_p M$ with
$\R\ell_p + T_p^cM = T_pM$. Right: representation of
the transfer of a chain and its direction
through an ambient biholomorphism
$h \colon \C^2 \longrightarrow {\C'}^2$, making a 
CR-diffeomorphism $h\vert_M \colon M \longrightarrow
M' := h(M)$.
\end{minipage}
\end{center}

\'Elie Cartan~{\cite{Cartan-1932-groupe-hyperspherique,
Cartan-1932-I, Cartan-1932-II}} discovered that nevertheless,
there do exist certain invariant CR-transversal curves, called {\sl
chains}, namely unparametrized curves $\CC_{p, \ell_p}$ uniquely
determined at each $p \in M$ and for each line $\ell_p \ni p$
complementary to $T_p^cM$ such that the nonzero tangent vector
$\dot{\CC}_{p, \ell_p}$ is directed by $\ell_p$,
but their existence always 
remained a bit mysterious.
This unique determination is similar to that for a
scalar second order ODE $\ddot{y} = H(t, y, \dot{y})$
for which a starting
point $y(0)$ and a starting vector $\dot{y}(0)$ must be prescribed,
but here, since $M$ is $3$-dimensional, chains are defined by a {\em
system} of {\em two} scalar second order ODEs,
as we now explain.

One may equip $\C^2$ with affine coordinates $(z, w) = (x + i\,y,\,
u+i\,v)$, centered at some reference point $p_0 = 0 \in M$ so that the
projection $T_0M \longrightarrow \R_{x,y,u}^3$ gives a local chart on
$M^3$ near the origin and even so that $T_0 M = \C_z \times (\R_u +
i\,\{0\})$, whence $T_0^c M = \C_z \times \{0\}$.  Then $M$ can be
$\mathcal{C}^\omega$ graphed as:
\[
v
\,=\,
F(z,\overline{z},u)
\,=\,
\sum_{j+k+l\,\geqslant\,1}\,
F_{j,k,l}\,
z^j\overline{z}^ku^l
\eqno
{\scriptstyle{(\overline{F_{k,j,l}}\,=\,F_{j,k,l})}}.
\]

Since $T_0^c M = \{u=0\}$ within $T_0M = \R_{x,y,u}^3$, 
any CR-transversal curve may be parametrized
as $t \longmapsto (x(t), y(t), t)$ with
$u(t) \equiv t$. One may show that
there exist certain functions $A$ and $B$ such that
the equations of chains write as a system:
\[
\ddot{x}
\,=\,
A\big(t,x,y,\dot{x},\dot{y}\big),
\ \ \ \ \ \ \ \ \ \ \ \ \ \ \ 
\ddot{y}
\,=\,
B\big(t,x,y,\dot{x},\dot{y}\big),
\]
but 
the explicit expressions of $A$ and $B$ in terms of $F$ and its
derivatives are huge, never shown in the literature
[part of the mystery].  This is because
chains are considered at {\em every} point $p \in M$ near $p_0 = 0 \in
M$, which requires hard elimination computations in the commutative
differential ring with variables $\big\{ F_{z^j \overline{z}^k u^l}
\big\}_{j,k,l \in \N}$ generated by the derivatives of $F$.  As shown
in~{\cite{Aghasi-Merker-Sabzevari-2011, Merker-Sabzevari-2014}} the
explicit expression of Cartan's primary invariant $\Iaux_{\sf
Cartan}$, whose identical vanishing characterizes local
biholomorphic equivalence to the Heisenberg sphere $\{ v' =
z'\overline{z}' \}$, is even huger.

Fortunately, we will see that thanks to plain translations $(z,w)
\longmapsto (z-z_p,\, w-w_p)$, one may `decipher' chains {\em only at
the origin} for a family of hypersurfaces $\{ M^p \}_{p\in M}$
passing through $0 \in \C^2$ and parametrized by all points $p \in M$
in the original hypersurface.
Section~{\ref{point-normalizations-M3-C2}} presents this
start. 

\smallskip

In the literature, chains are detected 
from rather advanced considerations:

\smallskip\noindent$\square$\,
either from an almost complete construction of 
an $\{e\}$-structure or of a Cartan connection
associated with the CR equivalence 
problem~{\cite{Cartan-1932-I, Cartan-1932-II,
Nurowski-Sparling-2003, Merker-Sabzevari-2014}};

\smallskip\noindent$\square$\,
or from an almost complete construction of a formal
or converging Moser-like normal form~{\cite{Chern-Moser-1974,
Jacobowitz-1985, Jacobowitz-1990}} 
for $M^3 \subset \C^2$ at the origin $0 \in M$.

\smallskip

\begin{center}
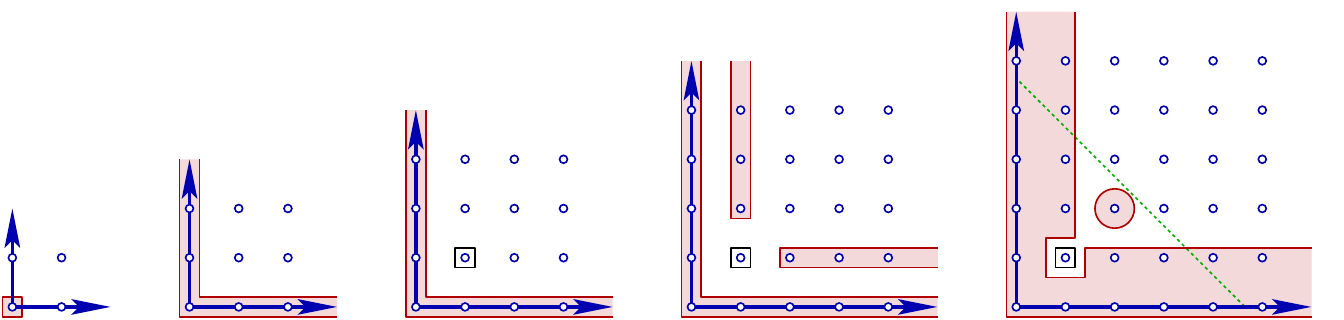

\nopagebreak
\begin{minipage}[t]{12.5cm}
\footnotesize
{\sc Figure~2:}
Successive annihilations (red dashed regions) of coefficient-functions 
$F_{j,k}(u)$ in the graphing function $v = \sum_{j,k}\,
z^j \overline{z}^k F_{j,k}(u)$ thanks to Moser's normalization
process, with $F_{1,1}(u) \equiv 1$, until first occurence of chains.
\end{minipage}
\end{center}

Let us comment only the second technique, which
proceeds in five steps.
At any reference point $p_0 = 0 \in M$, 
pick a curve $0 \in \gamma \subset M$
which is CR-transversal, namely $\dot{\gamma}(0) \not
\in T_0^c M$. Expand
$F$ in powers of $z, \overline{z}$ as:
\[
v
\,=\,
F(z,\overline{z},u)
\,=\,
\sum_{j,k}\,
z^j\overline{z}^k\,
F_{j,k}(u)
\eqno
{\scriptstyle{(F_{j,k}(u)\,:=\,
\sum_l\,F_{j,k,l}\,u^l)}}.
\] 

\medskip\noindent{\footnotesize\sf Step~1.}
Straighten $\gamma$ to be the $u$-axis, so that
$F(0,0,u) \equiv 0$, that is:
\[
v
\,=\,
\sum_{j+k\geqslant 1}\,
z^j\overline{z}^k\,
F_{j,k}(u).
\]

\medskip\noindent{\footnotesize\sf Step~2.}
Kill all harmonic terms $z^j\, F_{j,0}(u)$ and $\overline{z}^k\,
F_{0,k}(u)$, so that:
\[
v
\,=\,
\sum_{j\geqslant1\,\text{or}\,k\geqslant 1}\,
z^j\overline{z}^k\,
F_{j,k}(u).
\]

\medskip\noindent{\footnotesize\sf Step~3.}
Normalize $F_{1,1}(u) \longmapsto 1$, using the
assumption of Levi nondegeneracy, so that:
\[
v
\,=\,
z\overline{z}
+
\sum_{
j+k\,\geqslant\,3
\atop
j\geqslant2\,\text{or}\,k\geqslant2}\,
z^j\overline{z}^k\,
F_{j,k}(u).
\]

\medskip\noindent{\footnotesize\sf Step~4.}
Absorb all $z^1 \overline{z}^k\, F_{1,k}(u)$ and all
$z^j \overline{z}^1\, F_{j,1}(u)$ inside $z^1 \overline{z}^1$,
so that:
\[
v
\,=\,
z\overline{z}
+
\sum_{
j\,\geqslant\,2
\atop
k\,\geqslant\,2}\,
z^j\overline{z}^k\,
F_{j,k}(u).
\]

\medskip\noindent{\footnotesize\sf Step~5.}
Kill (in some way) $F_{2,2}(u)$, so that:
\[
v
\,=\,
z\overline{z}
+
z^3\overline{z}^2\,
F_{3,2}(u)
+
z^2\overline{z}^3\,
F_{2,3}(u)
+
\sum_{
j+k\,\geqslant\,6
\atop
j\,\geqslant\,2
\,{\rm and}\,
k\,\geqslant\,2}\,
z^j\overline{z}^k\,
F_{j,k}(u).
\]

\smallskip

Each one of these steps requires to perform an application of the
$\mathcal{C}^\omega$ implicit function theorem.  Next, what about
$F_{3,2}(u)$ and its conjugate $F_{2,3}(u) =
\overline{F_{3,2}(u)}${\bf ?}  One (known) paradox is that it is only
at an advanced stage of the progressive normalization process that one
can realize that the choice of a CR-transversal curve $\gamma$ should
{\em not} be made haphazardly.

Indeed, Proposition~6 in~{\cite[Chap.~4]{Jacobowitz-1990}}
states\,\,---\,\,not in the clearest thoughtful mathematical
way?\,\,---\,\,: {\em For each direction $\ell_p$ transverse to
$T_p^cM$ at $p \in M$, there exists a unique (unparametrized) real
analytic curve through $p$ and tangent to that direction such that
there exists some biholomorphism taking $M$ to:}
\[
v
\,=\,
\vert z\vert^2
+
\sum_{
j\geqslant 2
\atop
k\geqslant 2}\,
F_{j,k}(u)\,
z^j\overline{z}^k
\ \ \ \ \ \ \ \ \ \ \ \ \ \ \ 
\text{with}\ \ \
F_{3,2}(u)
\,\equiv\,
0,
\]
{\em and $\gamma$ to the $u$-axis}.

What are these curves? Why do they exist? Can one get them in advance?
Can one characterize them {\em geometrically}?  Without relying on the
existence of some normalizing biholomorphisms?

In fact, the proof of this Proposition~6 is the most technical and
difficult to follow in~{\cite{Chern-Moser-1974}}
or in~{\cite[Chap.~4]{Jacobowitz-1990}}. 
One first reason is that the argumentation appears almost at
the end of the normalization process, and a second reason is that 
it demands to perform biholomorphisms of the shape:
\[
z'
\,:=\,
\sum_{j=0}^\infty\,
z^j\,f_j(w),
\ \ \ \ \ \ \ \ \ \ \ \ \ \ \ 
w'
\,:=\,
\sum_{j=0}^\infty\,
z^j\,g_j(w),
\]
with $f_0(w) \neq 0$ required not to send the curve
$\{z = 0\} \cap M$ to 
the same curve $\{z' = 0\} \cap M'$\,\,---\,\,one 
really has to change the CR-transversal curve!\,\,---,
but this creates substantial computational obstacles. 

As an alternative, we will present a construction which is elementary,
simple, and requires almost no computation.  Furthermore, we will work
with power series in $3$ variables at one point, the origin, and only
up to order $5$ included.

Let therefore $0 \in M^3 \subset \C^2$ be $\mathcal{C}^\omega$ Levi
nondegenerate, graphed as $v = F(z, \overline{z}, u)$, with $0 \in M$.
We assign the weights $[x] := 1 =: [y]$ and $[u] := 2 =: [v]$. It is
well known that one can assume, with a {\em weighted} remainder, that
$M$ has equation:
\[
v
\,=\,
z\overline{z}
+
\sum_{3\leqslant\delta\leqslant 5}\,
\sum_{j+k+2l=\delta}\,
F_{j,k,l}\,
z^j\overline{z}^ku^l
+
{\rm O}(6).
\]

Anybody with a pen or a computer will reconstitute 
Proposition~{\ref{Prp-normalization-v-z-zbar-o-5}},
stating that there exists a change of holomorphic
coordinates in which $M$ becomes:
\[
v
\,=\,
z\overline{z}
+
{\rm O}(6).
\]

Next, the {\em key} fact is that the {\sl ambiguity}
of such a normalization up to (weighted) order $5$,
namely any biholomorphic equivalence:
\[
v
\,=\,
z\overline{z}
+
{\rm O}(6)
\ \ \ \ \
\xrightarrow[{\rule[0pt]{50pt}{0pt}}]{}
\ \ \ \ \
v'
\,=\,
z'\overline{z}'
+
{\rm O}(6),
\]
can be elementarily shown, by 
Proposition~{\ref{Prp-ambiguity-order-5}}, to coincide with 
the expansion, up to weighted order $5$, 
of the general isotropy group of the sphere 
$v = z\overline{z}$ $\longrightarrow$
$v' = z' \overline{z}'$ (without remainder),
which is known to be:
\[
z'
\,=\,
\frac{\lambda\,(z+\alpha\,w)}{
1-2i\overline{\alpha}\,z
-
(r+i\alpha\overline{\alpha})\,w},
\ \ \ \ \ \ \ \ \ \ \ \ \ \ \ \ 
w'
\,=\,
\frac{\lambda\overline{\lambda}\,w}{
1-2i\overline{\alpha}\,z
-
(r+i\alpha\overline{\alpha})\,w},
\]
with $\lambda \in \C\backslash \{0\}$,
$\alpha \in \C$, $r \in \R$.
(For $2$-nondegenerate constant Levi rank $1$
hypersurfaces $M^5 \subset \C^3$, this fact becomes
false, unfortunately~{\cite{Foo-Merker-Ta-2020}}.)

Then miraculously, the existence of Cartan-Moser chains
amounts to just understanding how the isotropy 
group of the model acts on CR-transversal objects!

This $5$-dimensional isotropy group has $5$ generators
${\sf D}$, ${\sf R}$, ${\sf I}_1$, ${\sf I_2}$,
${\sf J}$ which are $5$ linearly
independent holomorphic vector fields
${\sf X}$ with ${\sf X} + \overline{\sf X}$ tangent 
to $v = z \overline{z}$. Their expressions in
the intrinsic coordinates $(x,y,u) \in M^3$ read as
(Section~{\ref{intrinsic-automorphisms-sphere}}):
\[
\aligned
{\sf J}
+
\overline{\sf J}
&
\,=\,
(xu-x^2y-y^3)\,\partial_x
+
(x^3+xy^2+yu)\,\partial_y
+
\big(u^2-(x^2+y^2)^2\big)\,\partial_u,
\\
{\sf I}_2
+
\overline{\sf I}_2
&
\,=\,
(x^2-3y^2)\,\partial_x
+
(u+4xy)\,\partial_y
+
(2xu-2yx^2-2y^3)\,\partial_u,
\\
{\sf I}_1
+
\overline{\sf I}_1
&
\,=\,
(u-4xy)\,\partial_x
+
(3x^2-y^2)\,\partial_y
+
(-2x^3-2xy^2-2yu)\,\partial_u,
\\
{\sf R}
+
\overline{\sf R}
&
\,=\,
-\,y\partial_x
+
x\,\partial_y,
\\
{\sf D}
+
\overline{\sf D}
&
\,=\,
x\,\partial_x
+
y\,\partial_y
+
2u\,\partial_u.
\endaligned
\] 

Then according to the beautiful, highly conceptional, theory of
Lie~{\cite[Chap.~25]{Lie-Merker-2015}},
{\em see} also~{\cite{Olver-1995, Merker-2008,
Chen-Merker-2019}}, the action
of this group on first jets $(\dot{x}(t), \dot{y}(t))$
and on second jets $(\ddot{x}(t), \ddot{y}(t))$
of curves $t \longmapsto (x(t), y(t), t)$, equipped
with coordinates $(x_1, y_1)$ and
$(x_2, y_2)$, can be understood infinitesimally
by means of the {\sl prolongations} to the
second jet space $J_{1,2}^2$ of maps
$\R_u^1 \longrightarrow \R_{x,y}^2$,
thanks to straightforward universal formulas
(Sections~{\ref{prolongations-jet-1}} 
and~{\ref{prolongations-jet-2}}).
Since we work only at one point, namely above the origin,
it suffices to compute the coefficients,
in front of $\frac{\partial}{\partial x_1}$,
$\frac{\partial}{\partial y_1}$, 
$\frac{\partial}{\partial x_2}$, 
$\frac{\partial}{\partial y_2}$,
of these five prolonged vector fields
only for $x = y = u = 0$
(Section~{\ref{prolongations-jet-2}}):
\[
\begin{array}{ccccc}
& \partial_{x_1} & \partial_{y_1} & \partial_{x_2} & \partial_{y_2}
\\
{\sf D}^{(2)} & -x_1 & -y_1 & -3x_2 & -3y_2
\\
{\sf R}^{(2)} & -y_1 & x_1 & -y_2 & x_2
\\
{\sf I}_1^{(2)} & 1 & 0 & -4x_1y_1 & 6x_1^2+2y_1^2
\\
{\sf I}_2^{(2)} & 0 & 1 & -2x_1^2 - 6y_1^2 & 4x_1y_1
\\
{\sf J}^{(2)} & 0 & 0 & 0 & 0.
\end{array}
\]

From the first two columns that are 
everywhere of rank $2$,
it is clear that there does not exist any
{\em invariant} CR-transversal line $\ell_0 \ni 0$
with $\ell_0 \oplus T_0^cM = T_0M$.
Moreover, the action on such $\ell_0$ is {\em transitive}.

Next, 
by some kind of `{\sl algebraic miracle}' 
which can be verified by applying
a plain Gauss pivot to the above $4 \times 4$ submatrix:
\[
\left(\!
\begin{array}{cccc}
0 & 0 & -3x_2-6x_1^2y_1-6y_1^3 & 
-3y_2+6x_1y_1^2+6x_1^3
\\
0 & 0 & -y_2+2x_1y_1^2+2x_1^3 & 
x_2+2x_1^2y_1+2y_1^3
\\
1 & 0 & -2x_1^2-6y_1^2 & 4x_1y_1
\\
0 & 1 & -4x_1y_1 & 6x_1^2+2y_1^2
\end{array}
\!\right),
\]
there appears
to eyes
(Section~{\ref{prolongations-jet-2}}) 
a special surface $\Sigma_0^2 \subset \R_{x_1, y_1}^2
\times \R_{x_2, y_2}^2$, graphed as shown by the (redundant 
by pairs) entries $(1,3)$, $(1,4)$, $(2,3)$, $(2,4)$,
as:
\[
\Sigma_0^2
\,:=\,
\big\{
(x_1,y_1,x_2,y_2)
\in
\R^4
\colon\,\,
x_2
=
-2x_1^2y_1-2y_1^3,\,\,\,
y_2
=
2x_1y_1^2
+
2x_1^3
\big\},
\]
which is a $2$-dimensional orbit of the five prolonged vector fields
${\sf D}^2$, ${\sf R}^2$, ${\sf I}_1^{(2)}$, ${\sf I}_2^{(2)}$, ${\sf
J}^{(2)}$, while the complement $\R_{x_1, y_1,
x_2, y_2}^4 \big\backslash \Sigma_0^2$ is a
single orbit (Observation~{\ref{Obs-Sigma-0}}).

The existence of $\Sigma_0^2$ together with the normalizability
to $v' = z'\overline{z}' + {\rm O}(6)$ therefore explain
in an elementary manner the existence of 
Cartan-Moser chains above $0$.

\begin{center}
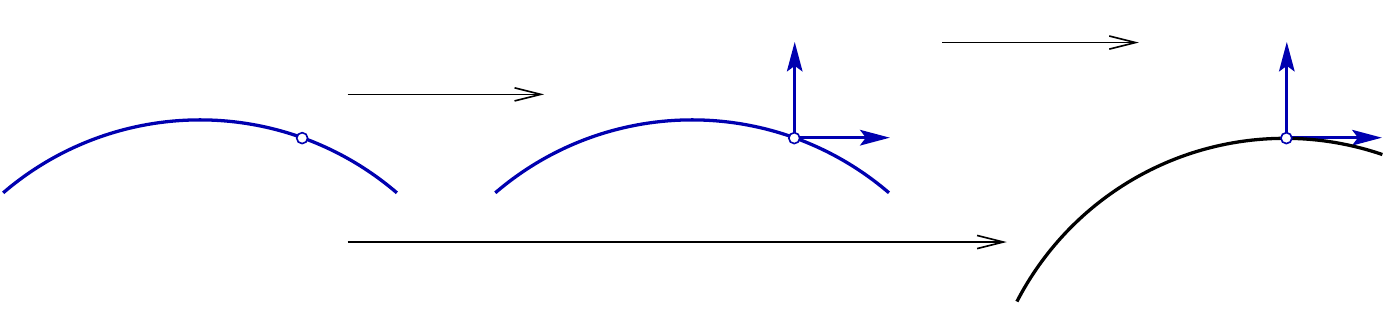

\nopagebreak
\begin{minipage}[t]{12.5cm}
\footnotesize
{\sc Figure~3:}
Centering (by translation) coordinates at an arbitrary point $p \in M$,
and sketching what any normalization map $\Phi_p$ does near $p = 0$.
\end{minipage}
\end{center}

Lastly, for any Levi nondegenerate hypersurface $M^3 \subset \C^3$,
we can define Cartan-Moser chains at any 
point $p \in M$ as follows.
Denote the translation map $\tau_p \colon
(M,p) \longrightarrow (M^p,0)$ by:
\[
\tau_p
\colon\ \ \ 
(z,w)
\,\,\longmapsto\,\,
\big(z-z_p,\,w-w_p\big)
\,=:\,
(z,w),
\]
denote {\em any} elementary normalization map as mentioned
above by:
\[
\Phi_p
\colon\ \ \
(M^p,0)
\,=\,
\Big\{
v
=
\smallsum{1\leqslant j+k+2l\leqslant 5}\,
F_{jkl}^p\,
z^j\overline{z}^ku^l
+
{\rm O}(6)
\Big\}
\,\,\longrightarrow\,\,
\big\{
v'
=
z'\overline{z}'
+
{\rm O}(6)
\big\}
\,\,=:\,\,
(N^p,0).
\]
Recall that the action of the $5$-dimensional isotropy group is
transitive on $1$-jets.

Given a $1$-jet $j_p^1$ at $p$, using {\em any} normalizing map
$\Phi_p \colon M^p \longrightarrow N^p$ which sends $(M^p, 0)$ to a
hypersurface $(N^p, 0)$ of equation $v' = z' \overline{z}' + {\rm
O}(6)$ {\em and also sends} $j_p^1$ to the flat $1$-jet $j_0^1 = (0,
0)$ at $0 \in N^p$, assign the $2$-jet $j_p^2$ of the Moser chain at
$p \in M$ associated with $j_p^1$ to be the inverse image of the {\em
flat} $2$-jet at $0 \in N^p$:
\[
j_p^2
\,:=\,
{\big(\Phi_p\circ\tau_p\big)^{(2)}}^{-1}
\big(0,0,0,0\big).
\]
It is not difficult to verify that this definition provides a map
$j_p^1 \longmapsto j_p^2(j_p^1)$ which is $\mathcal{C}^\omega$ on $M$.

Once chains are known, one can (re)start Step~1 above with the
CR-transversal curve $\gamma$ being a chain. Then Steps~2, 3, 4 go
without modification, while in Step~5, one realizes that $F_{3,2}(u)
\equiv 0$ automatically (Section~{\ref{link-F-3-2-origin}}), as a
consequence of the definition of chains
(Assertion~{\ref{Assertion-key-F-3-2-zero}}).

For self-contentness and for later
use in~{\cite{Foo-Merker-Ta-2020}},
although there is no originality,
we perform all these steps
in Section~{\ref{chain-straightening-harmonic-killing}},
{\ref{prenormalization}}, {\ref{Moser-normal-form-M3-C2}}
known as Propositions~1, 2, 3, 4, 5 
in~{\cite[Chap.~4]{Jacobowitz-1990}}.
We conclude by stating Moser's normal form theorem
in Section~{\ref{uniqueness-Moser-normal-form}} and by proving
some uniqueness property.

\medskip\noindent{\bf Acknowledgments.}
While the author was visiting Warsaw, 
Pawe{\l} Nurowski provided useful hints on how
certain {\sl distinguished curves} exist
in parabolic 
geometries~{\cite{Cap-Slovak-Zadnik-2004, Cap-Slovak-2009, 
Cap-Zadnik-2009}.

Grateful thanks are addressed to an anonymous referee
for a careful reading and for insightful suggestions.

\Section{\bf Point Normalizations of $\mathcal{C}^\omega$
Hypersurfaces $M^3 \subset \C^2$}
\label{point-normalizations-M3-C2}
\HEAD{{\ref{point-normalizations-M3-C2}}.~{\sf Point Normalizations 
of $\mathcal{C}^\omega$ Hypersurfaces $M^3 \subset \C^2$}
}{
Joël {\sc Merker}}

Consider a local real hypersurface $M^3 \subset \C^2$ of class (at
least) $\mathcal{C}^5$. In fact, we will mainly
work with $\mathcal{C}^\omega$ (real-analytic) objects,
and sometimes indicate what kind of lower regularity assumptions 
can be afforded.

In coordinates $(z,w) = (x+i\,y,\, u + i\,v)$,
assume $M$ is graphed as $v = F(z, \overline{z}, u)$, with $F \in
\mathcal{C}^5$. At all points $p = (z_p, w_p) \in M$ with $v_p =
F(z_p, \overline{z}_p, u_p)$, expand:
\[
v
\,=\,
F(z,\overline{z},u)
\,=\,
\sum_{j+k+l\leqslant 5}\,
{\textstyle{\frac{(z-z_p)^j}{j!}}}\,
{\textstyle{\frac{(\overline{z}-\overline{z}_p)^k}{k!}}}\,
{\textstyle{\frac{(u-u_p)^l}{l!}}}\,
F_{z^j\overline{z}^ku^l}(z_p,\overline{z}_p,u_p)
+
{\rm O}(6),
\]
subtract $v - v_p$, translate coordinates $z := z-z_p$, $w := w - w_p$,
and get a family
of hypersurfaces $M^p \subset \C^3$
passing through the origin: 
\[
v
\,=\,
F^p(z,\overline{z},u)
\,=\,
\sum_{1\leqslant j+k+l\leqslant 5}\,
z^j\overline{z}^ku^l\,
F_{j,k,l}^p
+
{\rm O}(6),
\]
namely with $F^p(0,0,0) = 0$, 
having coefficients $F_{j,k,l}^p := \frac{1}{j!}\, \frac{1}{k!}\,
\frac{1}{l!}\, F_{z^j\overline{z}^ku^l}(z_p,\overline{z}_p,u_p)$
smoothly parametrized by $p$. 
Thanks to this, working at {\em only one} point,
namely at the origin, we will treat {\em all} points $p \in M$.

Local biholomorphisms $h \colon M \longrightarrow M'$ between any two
CR manifolds respect by definition complex tangent bundles 
$h_\ast (T^c M) = T^c M'$.

\begin{Question}
\label{Question-CR-transversal}
{\sl Are there CR-transversal structures which are invariant
under biholomorphisms?}
\end{Question}

The goal of this note is to elaborate a simple, {\sl Lie-theoretic}
approach to this question which applies to any kind of CR structure,
does not require to fully solve any
equivalence problem, and does
not rest on the existence of Cartan-Tanaka connections.  To illustrate
the process on just one advanced example, we shall show how to recover
in a quite elementary way
the famous {\sl Moser chains} on Levi
nondegenerate hypersurfaces 
$M^3 \subset \C^2$. Forthcoming publications will
exhibit more about Lie's theoretical scope.

Since Question~{\ref{Question-CR-transversal}} is {\em invariant}, we
are allowed to perform {\sl normalizing} biholomorphisms in order to
`simplify' the equations $v = F^p(z, \overline{z}, u)$ of our
$p$-parametrized hypersurfaces $M^p$, before searching for
CR-transversal structures, if any.

After an elementary biholomorphism, it is well known that one can
assume:
\[
v
\,=\,
z\overline{z}
+
{\rm O}(3).
\]
This conducts to attribute weights $[z] := 1 =: [\overline{z}]$ and
$[w] := 2 =: [\overline{w}]$.  Up to order $5$, some monomials have
weight $> 5$, for instance $u^2 z^2$, and they will be
disregarded. Thus, with a now {\em weighted} remainder
${\rm O}(6)$:
\[
\aligned
v
\,=\,
F^p(z,\overline{z},u)
&
\,=\,
z\overline{z}
+
\sum_{3\leqslant \delta\leqslant 5}\,
\sum_{j+k+2l=\delta}\,
F_{j,k,l}^p\,
z^j\overline{z}^ku^l
+
{\rm O}(6).
\endaligned
\]

By performing biholomorphisms of the shape $z' = z + f_{\delta-1}
(z,w)$,
$w' = w + g_\delta(z,w)$, with appropriate polynomials $f_{\delta-1}$,
$g_\delta$ that are weighted homogeneous of degrees $\delta-1$,
$\delta$, it is not difficult to erase 
$F_\delta$ for $\delta = 3, 4, 5$.

\begin{Proposition}
\label{Prp-normalization-v-z-zbar-o-5}
Every $M^p$ can be normalized to $v = z\overline{z} + 0 + 0 + 0 + 
{\rm O}(6)$.\qed
\end{Proposition}

Of course, such a normalizing biholomorphism is not unique.

\begin{center}
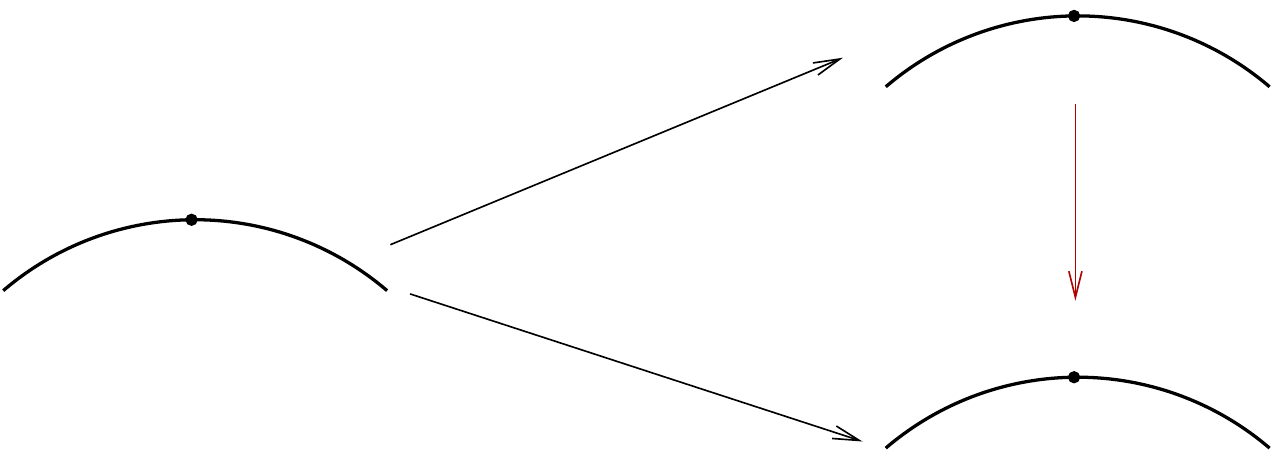

\nopagebreak
\begin{minipage}[t]{12.5cm}
\footnotesize
{\sc Figure~4:}
Representing two 6\textsuperscript{th}
order normalization maps at the origin
and calling `{\sl ambiguity}'
the `difference' (composition) between them.
\end{minipage}
\end{center}

The next statement\,\,---\,\,whose 
proof is also left as an exercise\footnote{\,
It turns out that all detailed proofs given later in
Sections~{\ref{chain-straightening-harmonic-killing}},
{\ref{prenormalization}},
{\ref{Moser-normal-form-M3-C2}}
do the job (solve the two exercises).}\,\,---\,\,determines 
the ambiguity transformation, 
which is obtained by expanding up to 
weight $5$ included the following two fractions in which
$\lambda \in \C^\ast$, $\alpha \in \C$, $r \in \R$ are free:
\leqnomode\usetagform{default}
\begin{align}
\label{group-G_0-5}
z'
\,=\,
\frac{\lambda\,(z+\alpha\,w)}{
1-2i\overline{\alpha}\,z
-
(r+i\alpha\overline{\alpha})\,w},
\ \ \ \ \ \ \ \ \ \ \ \ \ \ \ \ 
w'
\,=\,
\frac{\lambda\overline{\lambda}\,w}{
1-2i\overline{\alpha}\,z
-
(r+i\alpha\overline{\alpha})\,w}.
\end{align}

\begin{Proposition}
\label{Prp-ambiguity-order-5}
Every biholomorphism $z' = f(z,w) + {\rm O}(5)$, $w' = g(z,w) + {\rm
O}(6)$ with $f$, $g$ of weight $\leqslant 4, 5$ sending $v =
z\overline{z} + {\rm O}(6)$ to $v' = z' \overline{z}' + {\rm O}(6)'$ is
necessarily of the form:
\reqnomode\usetagform{EngelLie}
\begin{align}
z'
&
\,=\,
\lambda\,z
+
2i\lambda\overline{\alpha}\,
z^2
+
\big(
-4\lambda\overline{\alpha}^2
\big)\,
z^3
+
\big(
-8i\lambda\overline{\alpha}^3
\big)\,
z^4
\notag
\\
&
\ \ \ \ \ \ \ \ \ \ \ \
+
\lambda\alpha\,
w
+
\big(
3i\lambda\alpha\overline{\alpha}
+
\lambda r
\big)\,
zw
+
\big(
-8\lambda\alpha\overline{\alpha}^2
+
4i\overline{\alpha}\lambda r
\big)\,
z^2w
\notag
\\
&
\ \ \ \ \ \ \ \ \ \ \ \ \ \ \ \ \ \ \ \ \ \ \ \ \ \ \ \ \ \ \ \ \ \ 
\ \ \ \ \ \ \ \ \ \ \ \ \ \ \ \ \ \ \ \ \ \ \ \ \ \
+
\big(
\lambda\alpha r
+
i\lambda\alpha^2\overline{\alpha}
\big)\,
w^2
\notag
\\
w'
&
\,=\,
\lambda\overline{\lambda}\,
w
+
2i\lambda\overline{\lambda}\overline{\alpha}\,
zw
+
\big(
-4\lambda\overline{\lambda}\overline{\alpha}^2
\big)\,
z^2w
+
\big(
-8i\lambda\overline{\lambda}\overline{\alpha}^3
\big)\,
z^3w
\notag
\\
&
\ \ \ \ \ \ \ \ \ \ \ \ \ \ \ \ \ \ \ \
+
\big(
i\lambda\overline{\lambda}\alpha\overline{\alpha}
+
\lambda\overline{\lambda}r
\big)\,
w^2
+
\big(
4i\lambda\overline{\lambda}\overline{\alpha}r
-
4\lambda\overline{\lambda}\overline{\alpha}^2\alpha
\big)\,
zw^2.
\tag{\qed}
\end{align}
\end{Proposition}

But these formulas for this stability\big/ambiguity group 
are well known!

\Section{\bf Automorphisms of the Sphere 
$\{ {\rm Im}\, w = z\overline{z} \}$
Fixing the Origin}
\label{automorphisms-sphere-S3-C2}
\HEAD{{\ref{automorphisms-sphere-S3-C2}}.~{\sf Automorphisms 
of the Sphere $\{ {\rm Im}\, w = z\overline{z} \}$ Fixing the Origin}
}{
Joël {\sc Merker}}

Indeed, in $\C^2 \ni (z,w) = (x+i\,y,\, u + i\, v)$, 
consider the Heisenberg sphere:
\[
v
\,=\,
z\overline{z},
\]
which is biholomorphic, after a certain Cayley transform,
to the standard $3$-sphere $S^3 \subset \C^2$ 
minus one point sent to infinity.
It is known (details in {\cite[Sec.~3]{Aghasi-Merker-Sabzevari-2011}})
that the $5$-dimensional real Lie algebra $\mathfrak{g}^5$ of
holomorphic vector fields $X = a(z,w)\, \partial_z + b(z,w)\,
\partial_w$ with $a(0) = 0 = b(0)$ such that $X + \overline{X}$ is
tangent to $S_\ast^3$ consists of:
\[
\aligned
{\sf D}
&
\,:=\,
z\,\partial_z
+
2w\,\partial_w,
\\
{\sf R}
&
\,:=\,
iz\,\partial_z,
\\
{\sf I}_1
&
\,:=\,
(w+2iz^2)\,
\partial_z
+
2izw\,\partial_w,
\\
{\sf I}_2
&
\,:=\,
(iw+2z^2)\,\partial_z
+
2zw\,\partial_w,
\\
{\sf J}
&
\,:=\,
zw\,\partial_z
+
w^2\,\partial_w,
\endaligned
\]
with commutator table:
\begin{center}
\begin{tabular} [t] { l | l l l l l }
& ${\sf D}$ & ${\sf R}$ & ${\sf I}_1$ &
${\sf I}_2$ & ${\sf
J}$
\\
\hline
${\sf D}$ & $0$ & $0$ & ${\sf I}_1$ & ${\sf I}_2$ & $2\,{\sf J}$
\\
${\sf R}$ & $*$ & $0$ & $-{\sf I}_2$ & ${\sf I}_1$ & $0$
\\
${\sf I}_1$ & $*$ & $*$ & $0$ & $4\,{\sf J}$ & $0$
\\
${\sf I}_2$ & $*$ & $*$ & $*$ & $0$ & $0$
\\
${\sf J}$ & $*$ & $*$ & $*$ & $*$ & $0$ %
\end{tabular}
\end{center}
Integrating these fields, the finite equations of the
istropy Lie group $G^5 = {\sf Iso}(0)$ are:
\[
z'
\,=\,
\frac{\lambda\,(z+\alpha\,w)}{
1-2i\overline{\alpha}\,z
-
(r+i\alpha\overline{\alpha})\,w},
\ \ \ \ \ \ \ \ \ \ \ \ \ \ \ \ \ \ \ \ \ \ \ \ \ \
w'
\,=\,
\frac{\lambda\overline{\lambda}\,w}{
1-2i\overline{\alpha}\,z
-
(r+i\alpha\overline{\alpha})\,w},
\]
where $\lambda \in \C^\ast$, $\alpha \in \C$, $r \in \R$, as above.

So we know precisely the nonuniqueness (ambiguity)
in Proposition~{\ref{Prp-ambiguity-order-5}}.
Therefore,
we can pursue exploring our 
Question~{\ref{Question-CR-transversal}} 
by asking at first whether
some {\em tangential} 
(order $1$) CR-transversal invariant object exists.
 
\begin{Question}
\label{Question-v-CR-transversal}
{\sl 
Is there any vector 
$\vec{\bf v}_0 \in T_0M^p$
not complex-tangential $\vec{\bf v}_0 \not\in T_0^cM^p$
which would be invariant under biholomorphisms}\text{\bf ?}
\end{Question}

\begin{center}
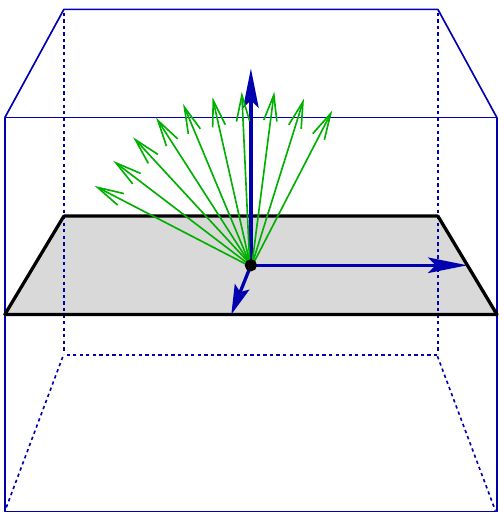

\nopagebreak
\begin{minipage}[t]{12.5cm}
\footnotesize
{\sc Figure~5:}
Representing horizontally the complex-tangential plane $T_0^c M^p$ of
$M^p$ at the origin within the $3$-dimensional $T_0 M^p$,
and drawing various vectors $\vec{\bf v} \in 
T_0 M^p \big\backslash T_0^c M^p$.
\end{minipage}
\end{center}

Predictably, the answer is {\em no}, because at order $1$, 
the above formulas read as linear transformations:
\[
z'
\,=\,
\lambda\,z
+
\lambda\,\alpha\,
w,
\ \ \ \ \ \ \ \ \ \ \ \ \ \ \ \ \ \ \ \ \ \ \ \ \ \
w'
\,=\,
\lambda\overline{\lambda}\,
w,
\]
and when $\alpha \in \C$ varies, the `{\sl slope}' of $\vec{\bf v}_0$
changes arbitrarily.
In fact, we must conceptualize carefully this intuition.

\Section{\bf Lie Jet Theory}
\label{Lie-jet-theory}
\HEAD{{\ref{Lie-jet-theory}}.~{\sf Lie Jet Theory}
}{
Joël {\sc Merker}}

The historical and philosophical monograph~{\cite{Merker-2010}}
explains how near 1870 Helmholtz involuntarily `{\sl invented}'
the so-called {\sl linearized
isotropy groups}, which were theoretically 
understood later by Sophus Lie after
finding a counterexample to Helmholtz's belief that any 
`macroscopic' (local) group action can be recovered 
`by integration' from
its `microscopic' (infinitesimal, linearized) behavior.

After Felix Klein's celebrated {\sl Erlanger program},
Lie indeed developped a fantastic theory of {\sl continuous
group actions}, having in mind applications to a new
`{\sl Galois theory}' of differential equations.
Lie erected a new theory of {\sl prolongations}
of group actions to jet spaces,
{\em see}~{\cite[Chap.~25]{Lie-Merker-2015}}. Lie also
conceptualized {\sl prolongations} of infinitesimal
transformations (vector fields) to jet spaces,
and this is exactly what we need here!

We must work with the three {\em intrinsic}, 
{\em real}, coordinates $(x,y,u)$ on 
$M$. A non CR-tangential vector $\vec{\bf v}_0 \in T_0M^p
\big\backslash T_0^c M^p$ 
can be represented as the derivative 
$\dot{\gamma}(0) = \vec{\bf v}_0$ 
of some 
parametrized real curve passing by the origin:
\[
t
\,\,\longmapsto\,\,
\big(x(t),y(t),u(t)\big)
\,=:\,
\gamma(t)
\eqno
{\scriptstyle{(\dot{\gamma}(0)\,\neq\,0)}}.
\]
Since $T_0^c M^p = \{ u = 0\}$, we have in fact
$\dot{u}(0) \neq 0$. 

So we are considering local curves $\R \longrightarrow \R^2$ 
graphed along the (vertical!) $u$-axis. We can then
represent by putting $u$ in he `horizontal' place as
$\big\{ (u, x(u), y(u)) \colon\, u \in \R \big\}$,
with {\em two}
graphing functions.

The associated jet space of order $2$\,\,---\,\,enough 
for our purposes\,\,---\,\,is equipped with further
independent coordinates corresponding
to $\dot{x}(u)$, $\dot{y}(u)$,
$\ddot{x}(u)$, $\ddot{y}(u)$:
\[
\big(
u,\,x,y,\,x_1,y_1,\,x_2,y_2
\big).
\]
We denote the first jet space by $J_{1,2}^1 \equiv \R^{1+2+2}$, and
this second jet space by $J_{1,2}^2 \equiv \R^{1+2+2+2}$.

Any diffeomorphism $(u,x,y) \longmapsto (u',x',y')$ lifts to jet
spaces of any order. The formulas rapidly become complicated
({\cite{Olver-1995, Merker-2008, Chen-Merker-2019}}).  Lie understood
this obstacle, and he linearized the formulas.

Indeed, 
by differentiating the prolongation to the second jet space of any
one-parameter diffeomorphism $\exp(t\vec{\bf v})(u,x,y)$ 
obtained as the flow
of a vector field $\vec{\bf v}$ on the base $\R^{1+2}$ , 
Lie introduced its {\sl
prolongations} $\vec{\bf v}^{(1)}$
to $J_{1,2}^1$ and $\vec{\bf v}^{(2)}$
to $J_{1,2}^2$.
A summarized presentation is available on pages~19--20
of~{\cite{Chen-Merker-2019}}.

Here, we just need to 
{\em apply} Lie's formulas. Start from a general vector field:
\[
\vec{\bf v}
\,:=\,
\xi(u,x,y)\,
\frac{\partial}{\partial u}
+
\varphi(u,x,y)\,
\frac{\partial}{\partial x}
+
\psi(u,x,y)\,
\frac{\partial}{\partial y},
\]
with smooth coefficients. Introduce the {\sl total differentiation
operator:}
\[
{\sf D}_u
\,:=\,
\frac{\partial}{\partial u}
+
x_1\,
\frac{\partial}{\partial x}
+
y_1\,
\frac{\partial}{\partial y}
+
x_2\,
\frac{\partial}{\partial x_1}
+
y_2\,
\frac{\partial}{\partial y_1}
+
x_3\,
\frac{\partial}{\partial x_2}
+
y_3\,
\frac{\partial}{\partial y_2}.
\]
Then the second prolongation of $\vec{\bf v}$:
\[
\aligned
\vec{\bf v}^{(2)}
&
\,=\,
\vec{\bf v}
+
\varphi_1\,
\frac{\partial}{\partial x_1}
+
\psi_1\,
\frac{\partial}{\partial y_1}
\\
&
\ \ \ \ \ \ \ \ \ 
+
\varphi_2\,
\frac{\partial}{\partial x_2}
+
\psi_2\,
\frac{\partial}{\partial y_2},
\endaligned
\]
has coefficients given uniquely 
by ({\cite{Olver-1995, Merker-2008, Chen-Merker-2019}}):
\[
\aligned
\varphi_1
&
\,:=\,
{\sf D}_u
\big(
\varphi
-
\xi\,x_1
\big)
+
\xi\,
x_2,
&
\ \ \ \ \ \ \ \ \ \ \ \ \ \ \ \ \ \ \ \ \ \ \ \ \ \
\psi_1
&
\,:=\,
{\sf D}_u
\big(
\psi
-
\xi\,y_1
\big)
+
\xi\,
y_2,
\\
\varphi_2
&
\,:=\,
{\sf D}_u
\big(
{\sf D}_u
\big(
\varphi
-
\xi\,x_1
\big)\big)
+
\xi\,
x_3,
&
\ \ \ \ \ \ \ \ \ \ \ \ \ \ \ \ \ \ \ \ \ \ \ \ \ \
\psi_2
&
\,:=\,
{\sf D}_u
\big(
{\sf D}_u
\big(
\psi
-
\xi\,y_1
\big)\big)
+
\xi\,
y_3.
\endaligned
\]

\Section{\bf Intrinsic Isotropy Automorphisms of the Sphere}
\label{intrinsic-automorphisms-sphere}
\HEAD{{\ref{intrinsic-automorphisms-sphere}}.~{\sf Intrinsic 
Isotropy Automorphisms of the Sphere}
}{
Joël {\sc Merker}}

Coming back to Question~{\ref{Question-v-CR-transversal}}, we must
apply Lie's prolongation formulas within the {\em first} jet space to
our $5$ vector fields ${\sf J}$, ${\sf I}_2$, ${\sf I}_1$, ${\sf R}$,
${\sf D}$.  But these vector fields ${\sf X} = a(z,w)\, \partial_z +
b(z,w)\, \partial_w$ were {\em extrinsic}, defined in $\C^2$, and
holomorphic!  Moreover, only their real parts $\frac{1}{2}\, \big(
{\sf X} + \overline{\sf X} \big)$ matter!

To apply Lie's theory, we must therefore write them up
in the {\em intrinsic} coordinates $(x,y,u) \in M^p$. 
We leave as an exercise to verify that 
the projection $\pi \colon (x,y,u,v) \longmapsto (x,y,u)$
is a chart on $S_\ast^3$ for which:
\[
\aligned
\pi_\ast\big(2\,\Re\,{\sf J}\big)
&
\,=\,
(xu-x^2y-y^3)\,\partial_x
+
(x^3+xy^2+yu)\,\partial_y
+
\big(u^2-(x^2+y^2)^2\big)\,\partial_u,
\\
\pi_\ast\big(2\,\Re\,{\sf I}_1\big)
&
\,=\,
(u-4xy)\,\partial_x
+
(3x^2-y^2)\,\partial_y
+
(-2x^3-2xy^2-2yu)\,\partial_u,
\\
\pi_\ast\big(2\,\Re\,{\sf I}_2\big)
&
\,=\,
(x^2-3y^2)\,\partial_x
+
(u+4xy)\,\partial_y
+
(2xu-2yx^2-2y^3)\,\partial_u,
\\
\pi_\ast\big(2\,\Re\,{\sf R}\big)
&
\,=\,
-\,y\partial_x
+
x\,\partial_y,
\\
\pi_\ast\big(2\,\Re\,{\sf D}\big)
&
\,=\,
x\,\partial_x
+
y\,\partial_y
+
2u\,\partial_u.
\endaligned
\] 
We will keep the same notation for these five
intrinsic vector fields.

\Section{\bf Prolongation to the Jet Space of Order $1$}
\label{prolongations-jet-1}
\HEAD{{\ref{prolongations-jet-1}}.~{\sf Prolongation
to the Jet Space of Order $1$}
}{
Joël {\sc Merker}}

As we said, it suffices to work above the origin $0 \in M^p$.  In
fact, the projectivization $\P(T_0 M^p) = \P^2$ of $T_0 M^p \cong
\R^3$ is a real projective plane. But excluding CR-tangential vectors,
we are considering only $\P^2 \backslash \P_\infty^1 = \R^2$, equipped
with affine coordinates $(x_1, y_1)$ as above.

This means that we are considering vectors $\vec{\bf v}_0 \in
T_0M^p \backslash T_0^c M^p$
of coordinates $(1, x_1^0, y_1^0)$, with unit coordinate $1$ along
the $u$-axis. Though we will not work in the projective space
$\P^2$, but only on its affine subset $\C^2 \subset \P^2$,
we mention that 
there are homogeneous coordinates 
$[U_1 \colon X_1 \colon Y_1]$ on $\P(T_0 M^p) = \P^2$
for which $\big[1 \colon \frac{X_1}{U_1} \colon \frac{Y_1}{U_1} \big]
=: (1, x_1, y_1)$.

\begin{center}
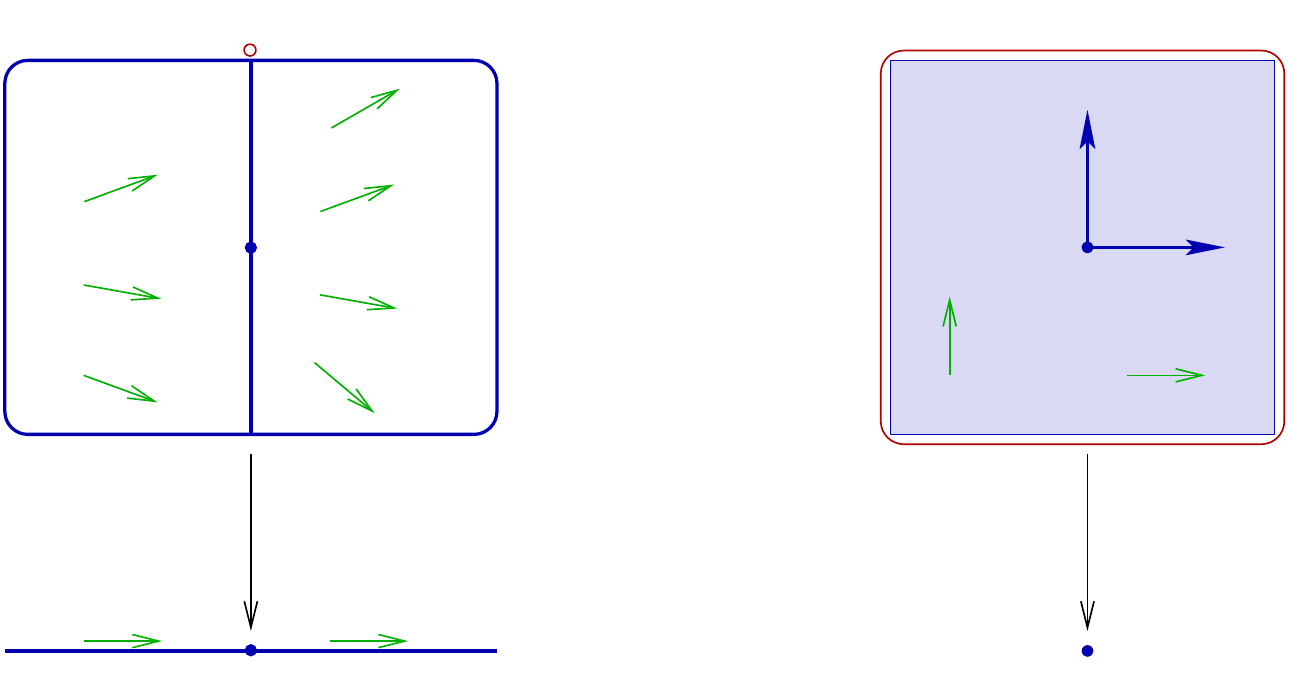

\nopagebreak
\begin{minipage}[t]{12.5cm}
\footnotesize
{\sc Figure~6:}
Left:
representing the first prolongation of a vector field $\vec{\bf v}$ on
$M^p$ to the first jet space $J_{1,2}^1$.  Right: Observing that, 
above the origin (only), the
first prolongations ${\sf I}_1^{(1)}$ and ${\sf I}_2^{(1)}$ of ${\sf
I}_1$ and ${\sf I}_2$ are straight (simple).
\end{minipage}
\end{center}

On the left, the figure represents this real $\P^2$ as a line, and on
the right, as a plane. The projective line
$\P_\infty^1$ at infinity is represented as a point,
and as a square perimeter.

By Lie's theory, any vector field $\vec{\bf v}$ on the base $M$ lifts
as a vector field $\vec{\bf v}^{(1)}$ on the first jet space
$J_{1,2}^1 = \R^{1+2+2}$.

Because our five intrinsic vector fields ${\sf J}$, ${\sf I}_1$, ${\sf
I}_2$, ${\sf R}$, ${\sf D}$ vanish at $u = x = y = 0$, their
prolongations will automatically be tangent to the fiber $\big\{ (0,
0, 0, x_1, y_1) \big\}$ above $(0,0,0)$ on the first jet space, a fiber
which identifies with $\R^2 = \P^2 \backslash \P_\infty^1$.

Lie's formulas yield the very simple values of these first
prolongations above the origin, namely for $x = y = u = 0$:
\[
\begin{array}{ccc}
& \partial_{x_1} & \partial_{y_1}
\\
{\sf D}^{(1)} & -x_1 & -y_1
\\
{\sf R}^{(1)} & -y_1 & x_1
\\
{\sf I}_1^{(1)} & 1 & 0
\\
{\sf I}_2^{(1)} & 0 & 1
\\
{\sf J}^{(1)} & 0 & 0
\end{array}
\]

Since the rank of the span of just 
${\sf I}_1^{(1)}$ and ${\sf I}_2^{(1)}$ is everywhere
equal to $2$, the orbit is the whole fiber $\R^2 = 
\{ (0, 0, 0, x_1, y_1)\}$, and this confirms what we already
guessed, namely 
that {\em there does not exist any 
biholomorphically invariant CR-transversal direction
$\ell_0 \subset T_0 M^p \big\backslash T_0^c M^p$}.

So what? All this for nothing? Let us keep hope by asking

\begin{Question}
{\sl 
Are there CR-transversal invariants of jet order $2$}{\bf ?}
\end{Question}

\Section{\bf Prolongation to the Jet Space of Order $2$}
\label{prolongations-jet-2}
\HEAD{{\ref{prolongations-jet-2}}.~{\sf Prolongation
to the Jet Space of Order $2$}
}{
Joël {\sc Merker}}

A non CR-tangential direction $\ell_0
\subset T_0 M^p \big\backslash T_0^c M^p$
can be represented as an order $1$ jet
$j_0^1 = (x_1^0, y_1^0)$.
A general jet of order two then writes as
$j_0^2 = \big( x_1^0, y_1^0, x_2^0, y_2^0 \big)$.

Since we just saw that the stability group
of the normalized equation $v = z \overline{z} + {\rm O}(6)$
for $M_p$, of dimension $5$, acts {\em transitively}
on first-order CR-transversal jets, it is clearly
impossible that a {\em unique} second order jet be invariant
under biholomorphisms. Anyway, it might be interesting
to see how the second order Lie prolongations
${\sf R}^{(2)}$, ${\sf D}^{(2)}$, ${\sf I}_1^{(2)}$, 
${\sf I}_2^{(2)}$, ${\sf J}^{(2)}$
act on second order jets.

Lie's formulas yield the very simple values of these first
prolongations above the origin, namely for $x = y = u = 0$:
\[
\begin{array}{ccccc}
& \partial_{x_1} & \partial_{y_1} & \partial_{x_2} & \partial_{y_2}
\\
{\sf D}^{(2)} & -x_1 & -y_1 & -3x_2 & -3y_2
\\
{\sf R}^{(2)} & -y_1 & x_1 & -y_2 & x_2
\\
{\sf I}_1^{(2)} & 1 & 0 & -4x_1y_1 & 6x_1^2+2y_1^2
\\
{\sf I}_2^{(2)} & 0 & 1 & -2x_1^2 - 6y_1^2 & 4x_1y_1
\\
{\sf J}^{(2)} & 0 & 0 & 0 & 0
\end{array}
\]

The key discovery, due to Cartan and then to Moser who expressed it
differently, now appears elementary.  But before writing the
statement, let us draw the key surface $\Sigma_0^2 \subset \R_{x_1,
y_1}^2 \times \R_{x_2, y_2}^2$ alluded to in the Introduction.

\begin{center}
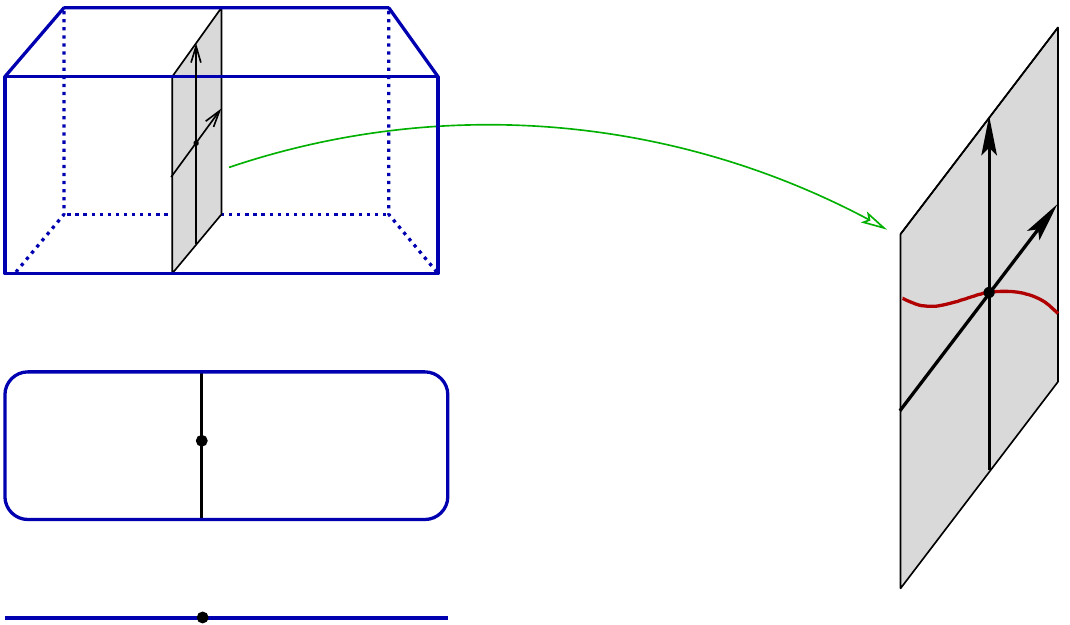

\nopagebreak
\begin{minipage}[t]{12.5cm}
\footnotesize
{\sc Figure~7:}
On the left, above $0 \in M^p$, we draw the first
jet fiber $J_{1,2}^1\big\vert_0 \cong
\R_{x_1, y_1}^2$ and the second jet fiber
$J_{1,2}^2\big\vert_0 \cong
\R_{x_1, y_1}^2 \times \R_{x_2, y_2}^2$. 
On the right, making a zoom, collapsing twice two dimensions into 
one dimension, we sketch
what the surface $\Sigma_0^2$ 
could be within $\R_{x_1, y_1}^2 \times \R_{x_2, y_2}^2$, 
representing it abusively as a $1$-curve in a $2$-plane.
\end{minipage}
\end{center}

\begin{Observation}
\label{Obs-Sigma-0}
On $\R^4 = \R_{x_1, y_1}^2 \times \R_{x_2, y_2}^2$, there
exists a unique invariant $2$-dimensional submanifold
$\Sigma_0^2 \subset \R^4$, algebraic, graphed as:
\[
x_2
\,=\,
-\,2\,
x_1^2y_1
-
2\,y_1^3,
\ \ \ \ \ \ \ \ \ \ \ \ \ \ \ \ \ \ \ \ \ \ \ \ \ \
y_2
\,=\,
2\,x_1y_1^2
+
2\,x_1^3.
\]
Moreover, the complement $\R^4 \backslash \Sigma_0^2$ is a unique
orbit under ${\sf D}^{(2)}$, ${\sf R}^{(2)}$,
${\sf I}_1^{(2)}$, ${\sf I}_2^{(2)}$, ${\sf J}^{(2)}$.
\end{Observation}

\proof
Any point of $\R^4$ can be represented as:
\[
x_2
\,=\,
-\,2\,
x_1^2y_1
-
2\,y_1^3
+
a_2,
\ \ \ \ \ \ \ \ \ \ \ \ \ \ \ \ \ \ \ \ \ \ \ \ \ \
y_2
\,=\,
2\,x_1y_1^2
+
2\,x_1^3
+
b_2,
\]
with some $(a_2, b_2) \in \R^2$.  A Gauss-pivot transforms the matrix
of the coefficients of the $4$ vector fields ${\sf D}^{(2)}$, ${\sf
R}^{(2)}$, ${\sf I}_1^{(2)}$, ${\sf I}_2^{(2)}$ into:
\[
\left(\!
\begin{array}{cccc}
0 & 0 & -3a_2 & -3b_2
\\
0 & 0 & -b_2 & a_2
\\
1 & 0 & -2x_1^2-6y_1^2 & 4x_1y_1
\\
0 & 1 & -4x_1y_1 & 6x_1^2+2y_1^2
\end{array}
\!\right).
\]
This matrix has maximal rank $4$ 
if and only if $(a_2, b_2) \neq (0, 0)$,
and constant rank $2$ for $(a_2, b_2) = (0, 0)$.
\endproof

In other words, to every (fixed) first order jet
$j_0^1 = (x_1, y_1)$ at the origin $0 \in M^p$ 
is associated a unique
second order jet at the origin:
\[
j_0^2
\,=\,
\Big(
x_1,y_1,\,\,
-2x_1^2y_1-2y_1^3,\,\,
2x_1y_1^2+2x_1^3
\Big),
\]
and since $\Sigma_0^2$ is invariant under the stability group
$G^5$ of $v = z\overline{z} + {\rm O}(6)$, 
this association is invariant under biholomorphic changes of
coordinates.

\Section{\bf Definition of Moser Chains}
\label{definition-Moser-chains}
\HEAD{{\ref{definition-Moser-chains}}.~{\sf Definition of Moser Chains}
}{
Joël {\sc Merker}}

Let us denote the translation map 
$\tau_p \colon (M,p) \longrightarrow (M^p,0)$
used in
Section~{\ref{point-normalizations-M3-C2}} by:
\[
\tau_p
\colon\ \ \ 
(z,w)
\,\,\longmapsto\,\,
\big(z-z_p,\,w-w_p\big)
\,=:\,
(z_0,w_0).
\]
Also, taking such coordinates $(z_0, w_0)$ around $(M^p, 0)$, 
let the punctual (at the origin) normalization
map offered by 
Proposition~{\ref{Prp-normalization-v-z-zbar-o-5}} be:
\[
\Phi_p
\colon\ \ \
(M^p,0)
\,=\,
\Big\{
v_0
=
\smallsum{1\leqslant j+k+2l\leqslant 5}\,
F_{0jkl}^p\,
z_0^j\overline{z}_0^ku_0^l
+
{\rm O}(6)
\Big\}
\,\,\longrightarrow\,\,
\big\{
v
=
z\overline{z}
+
{\rm O}(6)
\big\}
\,\,=:\,\,
(N^p,0),
\]
and abbreviate:
\[
\varphi 
\,:=\,
\Phi_p
\circ
\tau_p.
\]

\begin{center}
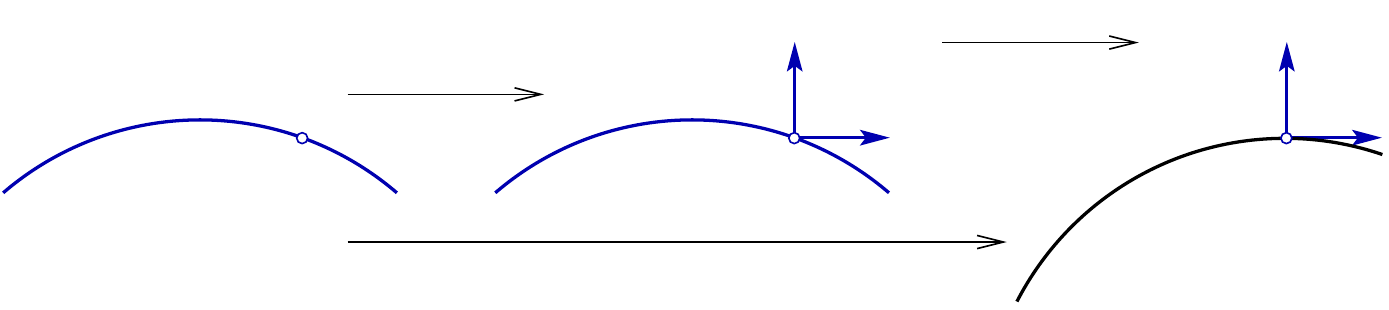

\nopagebreak
\begin{minipage}[t]{12.5cm}
\footnotesize
{\sc Figure~8:} Again, represent the translation map
$\tau_p$ and a normalizing map $\Phi_p$.
\end{minipage}
\end{center}

As in Observation~{\ref{Obs-Sigma-0}}, in the $2$-jet fiber
above $0 \in N^p$, introduce the surface:
\[
\Sigma_0
\,:=\,
\big\{
(x_1,y_1,x_2,y_2)
\in 
J_{N^p,0}^2
\colon\,
x_2
=
-2x_1^2y_1-2y_1^3,\,\,
y_2
=
2x_1y_1^2+2x_1^3
\big\}.
\]
Using the second prolongation $\varphi^{(2)}$, define the 
$2$-dimensional submanifold of $J_{M,p}^2$:
\[
\Sigma_p
\,:=\,
{\varphi^{(2)}}^{-1}
(\Sigma_0).
\]
Since $\varphi^{(1)}$ is a diffeomorphism $J_{M,p}^1
\overset{\sim}{\longrightarrow} J_{N^p,0}^1$, and 
the same about $\varphi^{(2)} \colon J_{M,p}^2
\overset{\sim}{\longrightarrow} J_{N^p,0}^2$,
this $\Sigma_p$
is also a graph, say of the form:
\[
x_2^p
\,=\,
A(x_1^p,y_1^p),
\ \ \ \ \ \ \ \ \ \ \ \ \ \ \ \ \ \ \ \ \ \ \ \ \ \
y_2^p
\,=\,
B(x_1^p,y_1^p),
\] 
with
$(x_1^p, y_1^p, x_2^p, y_2^p) \in 
J_{M,p}^2$, and 
with two functions $A$, $B$ which depend on $p$
{\em and also {\em a priori} on the normalizing map $\varphi$}.
\[
\xymatrix{
\Sigma_p
\ar@{^{(}->}[r]
&
J_{M,p}^2 
\ar[rr]^{\varphi^{(2)}}
\ar[d]
& &
J_{N^p,0}^2
\ar[d]
&
\ar@{^{(}->}[l]_f
\Sigma_0
\ar@/_3pc/[llll]_{{\varphi^{(2)}}^{-1}}
\\
&
J_{M,p}^1 
\ar[rr]^{\varphi^{(1)}}
\ar[d]
& &
J_{N^p,0}^1
\ar[d]
&
\\
&
(M,p)
\ar[rr]^{\varphi}
& &
(N^p,0)
&
}
\]

\begin{Assertion}
This graphed surface $\Sigma_p \subset J_{M,p}^2 \cong \R^4$ is {\em
independent} of the map $\varphi = \Phi_p \circ \tau_p$ normalizing
$v = F(z, \overline{z}, u)$ 
near $p$ to $v = z\overline{z} + {\rm O}(6)$ near $0$.
\end{Assertion}

\proof
Suppose another such normalizing map is given:
\[
\xymatrix{
& &
(N^p,0)
\ar[ddr]^{\psi\,:=\,\varphi_\prime\circ\varphi^{-1}}
\\
(M,p)
\ar[rru]^{\varphi}
\ar[rrrd]^{\varphi_\prime}
&
&
&
\\
& & &
(N_\prime^p,0),
}
\]
with $(N_\prime^p, 0)$ also of equation $v' = z' \overline{z}' + {\rm
O}(6)$. Define the special surface $\Sigma_0^\prime \subset
J_{N_\prime^p,0}^2$ by the {\em same} two graphed cubic equations
$x_2' = -2{x_1'}^2 y_1' - 2 {y_1'}^3$,\, 
$y_2' = 2\, x_1' {y_1'}^2 + 2\,
{x_1'}^3$, and then define similarly:
\[
\Sigma_p^\prime
\,:=\,
{\varphi_\prime^{(2)}}^{-1}
\big(\Sigma_0^\prime\big).
\]
Is it really true that $\Sigma_p^\prime = \Sigma_p${\bf ?}

Thanks to Proposition~{\ref{Prp-ambiguity-order-5}},
the {\sl relation map} $\psi := \varphi_\prime \circ \varphi^{-1}$
is a composition of flows of 
the five vector fields ${\sf D}$, ${\sf R}$,
${\sf I}_1$, ${\sf I}_2$, ${\sf J}$. But 
because the second
prolongations ${\sf D}^{(2)}$, ${\sf R}^{(2)}$,
${\sf I}_1^{(2)}$, ${\sf I}_2^{(2)}$, ${\sf J}^{(2)}$
of these fields are tangent to $\Sigma_0$
thanks to Observation~{\ref{Obs-Sigma-0}},
the map $\psi^{(2)}$ stabilizes the special surface:
\[
{\psi^{(2)}}^{-1}
\big(\Sigma_0^\prime\big)
\,=\,
\Sigma_0.
\]
Then as asserted:
\begin{align*}
\Sigma_p^\prime
&
\,=\,
{\varphi_\prime^{(2)}}^{-1}
\big(
\Sigma_0^\prime
\big)
\\
&
\,=\,
{\varphi_\prime^{(2)}}^{-1}
\big(
\psi^{(2)}(\Sigma_0)
\big)
\\
&
\,=\,
{\varphi_\prime^{(2)}}^{-1}
\Big(
\big(\varphi_\prime\circ\varphi^{-1}\big)^{(2)}(\Sigma_0)
\Big)
\\
&
\,=\,
\zero{
{\varphi_\prime^{(2)}}^{-1}
\circ
\varphi_\prime^{(2)}}
\circ
\big(
\varphi^{-1}
\big)^{(2)}
\big(\Sigma_0\big)
\\
&
\,=\,
{\varphi^{(2)}}^{-1}
\big(\Sigma_0\big)
\\
&
\,=\,
\Sigma_p.
\qedhere
\end{align*}
\endproof

\begin{Proposition}
There exist two $\mathcal{C}^\omega$ functions $A$ and $B$
such that $2$-jets are invariantly associated to 
CR-transversal $1$-jets as:
\reqnomode\usetagform{EngelLie}
\begin{align}
x_2
&
\,=\,
A\big(u,x,y,x_1,y_1\big),
\notag
\\
y_2
&
\,=\,
B\big(u,x,y,x_1,y_1\big).
\tag{\qed}
\end{align}
\end{Proposition}

These functions $A$ and $B$ can be made explicit in terms of $\big\{
F_{j,k,l}^p \big\}_{1 \leqslant j + k + l \leqslant 5}$, but
expressions are huge.  To these two jet equations is naturally
associated a system of two second order ordinary differential
equations:
\[
\aligned
\ddot{x}
&
\,=\,
A\big(u,x,y,\dot{x},\dot{y}\big),
\notag
\\
\ddot{y}
&
\,=\,
B\big(u,x,y,\dot{x},\dot{y}\big).
\endaligned
\]

\begin{Definition}
\label{Def-Moser-chain}
At a point $(u_p, x_p, y_p) \in M$, a {\sl Moser chain} directed
by some $1$-jet $(1, x_1^p, y_1^p)$ is
the unique solution $u \longmapsto (x(u), y(u))$ to the above 
$\mathcal{C}^\omega$ ODE system 
satisfying the initial conditions:
\[
\big(
x(u_p),\,y(u_p) 
\big) 
\,=\,
(x_p, y_p)
\ \ \ \ \ \ \ \ \ \ \ \ \ \ \ \ \ \ \ \
\text{and}
\ \ \ \ \ \ \ \ \ \ \ \ \ \ \ \ \ \ \ \
\big(\dot{x}(u_p),\,\dot{y}(u_p)\big)
\,=\,
(x_1^p,y_1^p).
\]
\end{Definition}

Equivalently, Moser chains 
$\big\{(u,x(u),y(u))\big\}$
are projections onto
the base space $M \ni (u, x, y)$ of integral curves of the
vector field on $J_{1,2}^1$:
\[
\frac{\partial}{\partial u}
+
x_1\,\frac{\partial}{\partial x}
+
y_1\,\frac{\partial}{\partial y}
+
A\big(u,x,y,x_1,y_1\big)\,
\frac{\partial}{\partial x_1}
+
B\big(u,x,y,x_1,y_1\big)\,
\frac{\partial}{\partial y_1}.
\]

Another equivalent, alternative, definition of $2$-jets of 
Moser chains uniquely associated with $1$-jets will be
useful later. 
Recall that first prolongations $\psi^{(1)}$ 
of maps like $\psi = \varphi_\prime \circ \varphi$
described in Proposition~{\ref{Prp-ambiguity-order-5}}
are {\em transitive} on $1$-jets, 
according to Section~{\ref{prolongations-jet-1}}.

So we can restrict considerations to normalizing
maps $\varphi = \tau_p \circ \Phi_p$ which 
send any $1$-jet $j_p^1$ at $p \in M$ 
to the flat $1$-jet $j_0^1 = (0, 0)$ at
$0 \in N^p$. 

\begin{Definition}
\label{Def-2-jet-chain-pullback-flat}
Given a hypersurface $M^3 \subset \C^2$, a point $p \in M$,
a $1$-jet $j_p^1$ at $p$, 
given 
the translation map $\tau_p \colon (M,p) 
\longrightarrow (M^p,0)$, 
and using {\em any} normalizing map $\Phi_p \colon M^p 
\longrightarrow N^p$ which sends $(M^p, 0)$ to
a hypersurface
$(N^p, 0)$ of equation
$v = z \overline{z} + {\rm O}(6)$
{\em and also sends} $j_p^1$ to the flat $1$-jet
$j_0^1 = (0, 0)$ at $0 \in N^p$, 
assign the $2$-jet $j_p^2$
of the Moser chain at $p \in M$ associated with $j_p^1$ 
to be
the inverse image of the {\em flat} $2$-jet at $0 \in N^p$:
\[
j_p^2
\,:=\,
{\big(\Phi_p\circ\tau_p\big)^{(2)}}^{-1}
\big(0,0,0,0\big).
\]
\end{Definition}

Thanks to the preceding reasonings, the result $j_p^2$ is
independent of the normalizing map 
$\Phi_p \circ \tau_p$ 
satisfying ${(\Phi_p \circ \tau_p)}^{(1)} (j_p^1) = (0,0)$,
the flat $1$-jet at $0 \in N^p$. 

\Section{\bf Link of Chains with $F_{3,2,0}^p$ at the Origin}
\label{link-F-3-2-origin}
\HEAD{{\ref{link-F-3-2-origin}}.~{\sf Link with $F_{3,2,0}^p$ 
at the Origin}
}{
Joël {\sc Merker}}

Once a point $p \in M$ and a CR-transversal $1$-jet $j_p^1$ at $p$ are
chosen, by known existence theorems, 
there is a unique local $\mathcal{C}^\omega$ curve $\gamma
\colon I \longrightarrow M$ passing through $p$ directed by $j_p^1$
which is a Moser chain.

Because such a chain is invariant under biholomorphisms, if one wants
to {\em normalize} the equation of a hypersurface $M^3 \subset \C^2$,
the very first natural normalization to perform is to {\em straighten}
(to normalize) such a chain. This can be done for {\em
any} CR-transversal curve, not necessarily a Moser chain.

\begin{Lemma}
Given any $\mathcal{C}^\omega$ curve $\gamma \colon (-1,1) 
\longrightarrow
M$ with $\gamma(0) = p \in M$ and $\dot{\gamma}(0) \not\in T_p^c M$,
there exist holomorphic coordinates $(z,w)$ centered at $p$ with $w = u
+ iv$ in which $M^p$ is graphed as $v = F^p(z, \overline{z}, u)$ such
that:
\[
\gamma(t)
\,=\,
\big(0,\,t+i\,0\big)
\eqno
{\scriptstyle{(t\,\in\,I)}}.
\]
\end{Lemma}

The (easy) proof will be written later in
Section~{\ref{chain-straightening-harmonic-killing}}.
So we may assume that $\{ (0,u) \}$ is a chain,
contained in $M^p$, whence $0 \equiv F^p (0, 0, u)$. 

In our preliminary Proposition~{\ref{Prp-normalization-v-z-zbar-o-5}},
the existence of Moser chains was unknown. Only successive Taylor
coefficients annihilations were performed.  Consequently, it is
necessary to restart the proof of
Proposition~{\ref{Prp-normalization-v-z-zbar-o-5}} with the
supplementary constraint to {\em keep invariant the straightened Moser
chain $\{(0,u)\}$}.

First of all, to annihilate 
all monomials except $z \overline{z}$ up to
weight $4$ is again possible by transformations $(z,w) \longmapsto
(z', w')$ sending (stabilizing) the $u$-axis to the
$u'$-axis\,\,---\,\,exercise\footnote{\, Again, it turns out that all
detailed proofs given later in
Sections~{\ref{chain-straightening-harmonic-killing}},
{\ref{prenormalization}}, {\ref{Moser-normal-form-M3-C2}} show how
to do it.}.

Furthermore, in weight $5$, all the monomials:
\[
z^5,\ \
z^4\overline{z},\ \
z\overline{z}^4,\ \
\overline{z}^5,\ \ \ \ \
z^3u,\ \
z^2\overline{z}u,\ \
z\overline{z}^2u,\ \
\overline{z}^3u,\ \ \ \ \
zu^2,\ \
\overline{z}u^2,
\]
can similarly be killed without modifying the unparametrized
straightened Moser chain $\{ z = v = 0\}$.  {\em Only the two
monomials $z^3\overline{z}^2$ and $z^2\overline{z}^3$ remain as
causing troubles}. In the notations of
Section~{\ref{point-normalizations-M3-C2}}, let us therefore formulate
a

\begin{Lemma}
\label{Lm-partly-normalize-order-5-keep-chain}
Every hypersurface $0 \in M^p \subset \C^3$ of equation:
\[
v_0
\,=\,
F_0^p(z_0,\overline{z}_0,u_0)
\ \ \ \ \
\text{with}
\ \ \ \ \
0
\,\equiv\,
F_0^p(0,0,u_0),
\]
having a Moser chain straightened to be $\{ (0, u_0) \}$,
can be normalized without deforming the chain being
$\{(0,u)\}$, into a hypersurface of equation:
\leqnomode\usetagform{default}
\begin{align}
\label{eq-with-F-3-2}
N^p
\colon\ \ \ \ \
v
\,=\,
F^p(z,\overline{z},u)
\,=\,
z\overline{z}
+
F_{3,2,0}^p\,
z^3\overline{z}^2
+
\overline{F}_{3,2,0}^p\,
z^2\overline{z}^3
+
{\rm O}(6).
\end{align}
\end{Lemma}

Now, remember that Proposition~{\ref{Prp-normalization-v-z-zbar-o-5}}
asserted that the remaining coefficient $F_{3,2,0}^p$, 
can be also killed.
However, there is a supplementary constraint, now.

\begin{Question}
{\sl Can one annihilate $F_{3,2,0}^p$ 
without unstraightening the chain?}
\end{Question}

It turns out that the answer is `{\sl no-becomes-yes}'! Indeed, for
some subtle reason which 
lies in the definition of
chains, it will soon turn out that {\em this coefficient
$F_{3,2,0}^p$ needs not be annihilated}, because it will be shown to
be already zero for free! Let us explain this key fact
which will be very useful later in
Assertion~{\ref{Assertion-crucial-F-3-2-zero}}.

\begin{Assertion}
\label{Assertion-key-F-3-2-zero}
If $F^p(z, \overline{z}, u)$ is as 
in~({\ref{eq-with-F-3-2}}) 
with $0 \equiv F^p(0, 0, u)$ and with $\{(0,u)\}$ being a chain,
then $F_{3,2,0}^p = 0$.
\end{Assertion}

\proof
Denote $h_0 \colon M^p \longrightarrow N^p$ one incomplete
normalizing map given by 
Lemma~{\ref{Lm-partly-normalize-order-5-keep-chain}}.
Since $h_0$ sends $\{(0,u_0)\}$ to $\{(0,u)\}$,
it sends the flat $1$-jet $j_{M^p, 0}^1 = (0,0)$
to the flat $1$-jet $j_{N^p, 0}^1 = (0,0)$.
We will apply Definition~{\ref{Def-2-jet-chain-pullback-flat}}
to $(N^p, 0)$ with $j_{N^p, 0}^1 = (0,0)$.

We know by Proposition~{\ref{Prp-normalization-v-z-zbar-o-5}},
that it is possible to continue to
perform normalizations by means of a further map: 
\[
\xymatrix{
M^p
\ar[rr]^{h_0}
& &
N^p
\ar[rr]^{h}
& &
N_\prime^p,
}
\]
in order that $N_\prime^p$ has equation $v' = z' \overline{z}' + {\rm
O}(6)$. In fact, the map
$h = (z + f_4,\, w+g_5) = (z', w')$ with:
\leqnomode\usetagform{default}
\begin{equation}
\label{imposed-form-F-3-2}
\begin{aligned}
z'
&
\,:=\,
z
-
i\,F_{3,2,0}^p
z^2w
-
{\textstyle{\frac{1}{4}}}\,
F_{2,3,0}^p\,
w^2
+
{\rm O}(5),
\\
w'
&
\,:=\,
w
-
{\textstyle{\frac{i}{2}}}\,
F_{3,2,0}^p\,
zw^2
+
{\rm O}(6),
\end{aligned}
\end{equation}
works. Because $h = (z,w) + {\rm O}_{z,w}(2)$,
this maps sends the flat $1$-jet at $0 \in
N^p$ to the flat $1$-jet at $0 \in N_\prime^p$. Then according
to
Definition~{\ref{Def-2-jet-chain-pullback-flat}} of a Moser chain, the
$2$-jet of the Moser chain at $0 \in N^p$ along $\{ (0,u)
\}$\,\,---\,\,which is flat!\,\,---\,\,{\em must be} 
the inverse image,
through ${h^{(2)}}^{-1}$, of the flat $2$-jet at $0 \in
N_\prime^p$.  Equivalently, $h$ must send the flat $2$-jet
at $0 \in N^p$ to the flat $2$-jet at $0 \in N_\prime^p$.

Let us write a flat $2$-jet at $0 \in N^p$ as a parametrized
curve $\R_u \longrightarrow \R_{x,y}^2$:
\[
x
\,=\,
{\rm O}_3(u),
\ \ \ \ \ 
y
\,=\,
{\rm O}_3(u).
\]
Then at $0 \in N_\prime^p$, 
do we also have $x' = {\rm O}_3(u')$ and $y' = {\rm
O}_3(u')$ through the
map~({\ref{imposed-form-F-3-2}})? 
We claim: {\em No if $F_{3,2,0}^p \neq 0$}!

Indeed, it comes $z = x + i\, y = {\rm O}_3(u)$, 
hence $w = u + i\, z\overline{z} + {\rm O}(5) = u + {\rm O}_3(u)$,
and also $u' = u + {\rm O}_3(u)$
or inversely $u' + {\rm O}_3(u') = u$, 
whence:
\[
x'+i\,y'
\,=\,
\,-
\tfrac{1}{4}\,
F_{3,2,0}^p\,
{u'}^2
+
{\rm O}_3(u).
\qedhere
\]
\endproof

This Lie-theoretic construction of Moser chains can be applied to any
CR manifold, and the paper could certainly stop at this point.

Ideed, we would like to mention that the normalizations 
applied in the
remainder of this paper, {\em i.e.} in the
next Sections~{\ref{chain-straightening-harmonic-killing}},
{\ref{prenormalization}},
{\ref{Moser-normal-form-M3-C2}},
{\ref{uniqueness-Moser-normal-form}},
are known to be done in the general case of hypersurfaces
$M^{2n+1} \subset \C^{n+1}$ in any
CR dimension $n \geqslant 1$ by Chern-Moser
in their celebrated work~{\cite{Chern-Moser-1974}}. 
More particularly, in part (d), page
246 of~{\cite{Chern-Moser-1974}}, Chern-Moser 
briefly concentrate on the
specific case of real hypersurfaces in $\C^2$.

Although Chern-Moser did not mention precisely all the intermediate
normalizations which are applicable in $\C^2$, Jacobowitz in Chapter~4
of his monograph~{\cite{Jacobowitz-1990}} endeavoured to detect and to
explain in $\C^2$ those appropriate normalizations.

But since, to the best of our knowledge, 
there is no considerable work in
the literature specifying such normalizations, 
we hope that the rest of the paper may raise interest 
of readers who want to learn Chern-Moser's normalizations 
in the specific case of $\C^2$.
Proofs are neither straightforward, nor elementary, 
because they require an
intensive, repeated use of the {\sl implicit function theorem}.

Thus, although the next results can not be regarded as new, 
for self-contentness reasons, 
and in order to prepare forthcoming works 
on new kinds of CR structures ({\em cf.}
{\em e.g.}~{\cite{Foo-Merker-Ta-2020}}), 
let us start to reconstitute the
Chern-Moser normalization theory in $\C^2$, 
setting up fully detailed arguments readable by
non-experts.

\Section{\bf Chain Straightening and Harmonic Killing}
\label{chain-straightening-harmonic-killing}
\HEAD{{\ref{chain-straightening-harmonic-killing}}.~{\sf Chain 
Straightening and Harmonic Killing}
}{
Joël {\sc Merker}}

The main feature being that Moser chains are 
{\em biholomorphically invariant}, it is natural to 
take them as a starting point for the process
of normalization. 

Let $M^3 \subset \C^2$ be a Levi nondegenerate hypersurface
passing by the origin $0 \in M$. Since $T_0^c M \cong \C$,
an appropriate $\C$-linear transformation makes $T_0^c M = 
\C_z \times \{ 0\}$ in coordinates $(z, w) \in \C^2$.

Our goal is to transform $M$ into certain {\sl normal forms},
by performing biholomorphisms fixing the origin:
\[
\aligned
\C^2
\,\supset\,
M^3
\ \ \ \ \ \ \ \ \ \ 
&
\xrightarrow[{\rule[0pt]{50pt}{0pt}}]{{\sf normalize}}
\ \ \ \ \ \ \ \ \ \ 
{M'}^3
\,\subset\,
{\C'}^2,
\\
\big(z,w\big)
\ \ \ \ \ \ \ \ \ \ 
&
\xrightarrow[{\rule[0pt]{50pt}{0pt}}]{}
\ \ \ \ \ \ \ \ \ \
\big(
f(z,w),\,
g(z,w)
\big)
\,\,=:\,
\big(z',w'\big).
\endaligned
\]

All objects will be real analytic ($\mathcal{C}^\omega$).
Thus with $w = u + i\, v$ and $w' = u' + i\, v'$, 
both hypersurfaces $M$ and $M'$ are $\mathcal{C}^\omega$-graphed as:
\[
v
\,=\,
F\big(
z,\overline{z},u
\big)
\ \ \ \ \ \ \ \ \ \ \ \ \ \ \ \ \ \ \ \
\text{and}
\ \ \ \ \ \ \ \ \ \ \ \ \ \ \ \ \ \ \ \
v'
\,=\,
F'\big(
z',\overline{z}',u'
\big).
\]
We also assume $T_0^c M' = \{ w' = 0\}$.

Expand $F$ as:
\[
F\big(
z,\overline{z},u
\big)
\,=\,
\sum_{j+k+l\geqslant 1}\,
F_{j,k,l}\,
z^j\overline{z}^ku^l,
\]
with $F_{j,k,l} \in \C$. Define:
\[
\overline{F}
\big(
z,\overline{z},u
\big)
\,:=\,
\sum_{j+k+l\geqslant 1}\,
\overline{F}_{j,k,l}\,
z^j\overline{z}^ku^l.
\]
From $\overline{v} = v$, it comes $\overline{ F(z, \overline{z},
u) } = F (z, \overline{z}, u)$, whence:
\leqnomode\usetagform{default}
\begin{align}
\label{reality-F}
\overline{F}
\big(
\overline{z},z,u
\big)
\,\equiv\,
F\big(
z,\overline{z},u
\big).
\end{align}
Applying $\frac{1}{j!} \partial_z^j \frac{1}{k!} 
\partial_{\overline{z}}^k \frac{1}{l!}\, \partial_u^l$ at
$(z, \overline{z}, u) = (0, 0, 0)$ we get:
\[
\overline{F}_{k,j,l}
\,=\,
F_{j,k,l}.
\]

The hypothesis that
the biholomorphism $(z, w) \longmapsto \big( f(z,w), g(z,w) \big)
=: (z', w')$ fixing the origin sends $M$ to $M'$ expresses
as a {\sl fundamental identity:}
\leqnomode\usetagform{default}
\begin{small}
\begin{align}
\label{fundamental-identity}
0
&
\,\equiv\,
-\,
{\textstyle{\frac{1}{2i}}}\,
g\big(z,\,u+i\,F(z,\overline{z},u)\big)
+
{\textstyle{\frac{1}{2i}}}\,
\overline{g}
\big(
\overline{z},\,
u-i\,F(z,\overline{z},u)
\big)
\,+
\notag
\\
&
\ \ \ \ \
+
F'
\Big(
f\big(z,\,u+i\,F(z,\overline{z},u)\big),\,\,
\overline{f}\big(\overline{z},\,u-i\,F(z,\overline{z},u)\big),\,\,
{\textstyle{\frac{1}{2}}}\,
g\big(z,\,u+i\,F(z,\overline{z},u)\big)
+
{\textstyle{\frac{1}{2}}}\,
\overline{g}\big(\overline{z},\,u-i\,F(z,\overline{z},u)\big)
\Big),
\end{align}
\end{small}
which holds in $\C\{z, \overline{z}, u\}$.

According to the preceding sections, for any CR-transversal
$1$-jet $j_0^1$ at $0 \in M$, there exists a Moser chain
directed by $j_0^1$ at $0$. We let $\gamma \colon I \longrightarrow 
M$ with $\gamma(0) = 0$ and $0 \in I \subset \R$ an interval,
be such a chain. 
In fact, the next statement is true for any local CR-transversal
curve.

\begin{Lemma}
Let $\gamma \colon I \longrightarrow M$ be a local 
$\mathcal{C}^\omega$ curve with $\gamma(0) = 0 \in M$ and 
$\dot{\gamma}(0) \not\in T_0^cM = \{ w= 0\}$. 
Then there exists a biholomorphism $(z, w) \longmapsto (z', w')$
stabilizing $T_0^c M' = \{ w' = 0\}$ which sends $\gamma$
to the curve $\gamma_\prime(t) = (0, t)$ straightened along the
$v'$-axis.
\end{Lemma}

Notice that a third direction 
$\dot{\gamma}_\prime(0) \in T_0 M' \big\backslash
T_0^c M'$ implies $T_0 M' = \{ u' = 0\}$.

\proof
Write:
\[
\gamma(t)
\,=\,
\big(
\varphi(t),\,
\psi(t)
\big)
\,\in\,
\C\times \C.
\]
By assumption, $\dot{\psi}(0) \neq 0$. Thus the map:
\[
z
\,:=\,
z'
+
\varphi(w'),
\ \ \ \ \ \ \ \ \ \ \ \ \ \ \ \ \ \ \ \ \ \ \ \ \ \
w
\,:=\,
\psi(w'),
\]
establishes a biholomorphism (inverse).

Similarly, the target curve writes $\gamma_\prime(t) = \big(
\varphi_\prime(t), \psi_\prime(t) \big)$. Thus for all $t \in I$:
\[
\varphi(t)
\,\equiv\,
\varphi_\prime(t)
+
\varphi\big(\psi_\prime(t)\big)
\ \ \ \ \ \ \ \ \ \ \ \ \ \ \ \ \ \ \ \
\text{and}
\ \ \ \ \ \ \ \ \ \ \ \ \ \ \ \ \ \ \ \
\psi(t)
\,\equiv\,
\psi\big(\psi_\prime(t)\big).
\]
The second equation and the invertibility of $\psi$ forces
$t \equiv \psi_\prime(t)$. Replacing this in the first equation
yields $0 \equiv \varphi_\prime(t)$.
\endproof

Consequently, the graphing function of the transformed hypersurface 
writes, after erasing the primes, as:
\[
M
\colon\ \ \ \ \
v
\,=\,
F\big(z,\overline{z},u\big),
\]
with $F = {\rm O}(2)$ and $F(0,0,u) \equiv 0$.

\begin{Lemma}
There exists a biholomorphism of the form:
\[
z'
\,:=\,
z,
\ \ \ \ \ \ \ \ \ \ \ \ \ \ \ \ \ \ \ \ \ \ \ \ \ \
w'
\,:=\,
w
+
g(z,w),
\]
with $g = {\rm O}(2)$ and $g(0, w) \equiv 0$, which transforms
$\{ v = F\}$ into $\{v' = F'\}$ satisfying:
\[
0
\,\equiv\,
F'\big(z',0,u'\big)
\,\equiv\,
F'\big(0,\overline{z}',u'\big).
\]
\end{Lemma}

The second vanishing follows from the first, by~({\ref{reality-F}}).
Notice that $F'(0, 0, u') \equiv 0$ is preserved.

\proof
If such a biholomorphism exists, the fundamental identity writes
for it:
\leqnomode\usetagform{default}
\begin{align}
\label{F-F-prime-g-bar}
0
&
\,\equiv\,
-\,
F\big(z,\overline{z},u\big)
-
{\textstyle{\frac{1}{2i}}}\,
g\big(z,\,u+i\,F(z,\overline{z},u)\big)
+
{\textstyle{\frac{1}{2i}}}\,
\overline{g}
\big(
\overline{z},\,
u-i\,F(z,\overline{z},u)
\big)
\,+
\notag
\\
&
\ \ \ \ \
+
F'
\Big(
z,\,
\overline{z},\,\,
u
+
{\textstyle{\frac{1}{2}}}\,
g\big(z,\,u+i\,F(z,\overline{z},u)\big)
+
{\textstyle{\frac{1}{2}}}\,
\overline{g}\big(\overline{z},\,u-i\,F(z,\overline{z},u)\big)
\Big).
\end{align}
We want $F'(z', 0, u') \equiv 0$. 
If this goal would be reached, putting $\overline{z} := 0$,
we would deduce:
\leqnomode\usetagform{default}
\begin{align}
\label{F-zbar-zero}
0
\,\equiv\,
-\,F(z,0,u)
-
{\textstyle{\frac{1}{2i}}}\,
g\big(z,\,u+i\,F(z,0,u)\big)
+
{\textstyle{\frac{1}{2i}}}\,
\overline{g}\big(0,\,u-i\,F(z,0,u)\big)
+
0.
\end{align}
By luck, such an equation can be used to defined $g(z,w)$ uniquely,
even with the supplementary condition that 
the last term be identically zero.

Indeed, by $F = {\rm O}(2)$, 
the implicit function theorem enables to invert:
\[
u
+
i\,F(z,0,u)
\,=:\,
\omega
\ \ \ \ \ \ \ \ \ \ \ \ \ \ \ \ \ \ \ \
\Longleftrightarrow
\ \ \ \ \ \ \ \ \ \ \ \ \ \ \ \ \ \ \ \
u
\,=\,
\TT(z,\omega)
\,=\,
\omega
+
{\rm O}(2).
\]
Define therefore $g(z,\omega)$, after erasing the third
term $\frac{1}{2i}\, \overline{g}$ above, by:
\[
0
\,\equiv\,
-\,
F\big(z,0,\TT(z,\omega)\big)
-
{\textstyle{\frac{1}{2i}}}\,
g(z,\omega)
+
0,
\]
and notice then that {\em because $F(0,0,u) \equiv 0$ by assumption},
we fulfill by setting $z := 0$, :
\[
0
\,\equiv\,
g
\big(
0,\,
\omega
\big).
\]
Thus, ({\ref{F-zbar-zero}}) really holds with
$\frac{1}{2i}\, \overline{g} = 0$, and then
coming back to  
({\ref{F-F-prime-g-bar}})$\big\vert_{\overline{z} = 0}$,
we get as desired:
\[
0
\,\equiv\,
0
+
F'
\Big(
z,\,0,\,
u
+
{\textstyle{\frac{1}{2}}}\,
g\big(z,\,u+i\,F(z,0,u)\big)
\Big).
\qedhere
\]
\endproof

\Section{\bf Prenormalization}
\label{prenormalization}
\HEAD{{\ref{prenormalization}}.~{\sf Prenormalization}
}{
Joël {\sc Merker}}

Now, erase the primes, and assume $0 \equiv F(z, 0, u)$. Write:
\[
v
\,=\,
F\big(z,\overline{z},u\big)
\,\,=\,\,
z\overline{z}\,F_{1,1}(u)
+
\sum_{j+k\geqslant 3
\atop
j\geqslant1,\,k\geqslant 1}\,
z^j\overline{z}^kF_{j,k}(u)
\,\,=\,\,
z\overline{z}\,F_{1,1}(u)
+
z^2\overline{z}\,
\big(\cdots\big)
+
\overline{z}^2z\,
\big(\cdots\big).
\] 
Since $M$ is Levi nondegenerate at $0$, after a $\C$-linear 
transformation, we make:
\[
F_{1,1}(0)
\,=\,
1.
\]
This equality $F_{1,1} (0) = 1$ is known as Poincar\'e's realization
of nondegenerate hypersurfaces in $\C^2$. It is quite crucial in the
Chern-Moser normal form construction.

\begin{Lemma}
There exists a biholomorphism of the form:
\[
z'
\,:=\,
z\,\varphi(w),
\ \ \ \ \ \ \ \ \ \ \ \ \ \ \ \ \ \ \ \ \ \ \ \ \ \
w'
\,:=\,
w,
\]
which transforms $M = \{ v = F \}$ into $M'$ with:
\[
v'
\,=\,
F'
\,=\,
z'\overline{z}'
+
{z'}^2\overline{z}'\,
\big(\cdots\big)
+
{\overline{z}'}^2z'\,
\big(\cdots\big).
\]
\end{Lemma}

So we may normalize $F_{1,1}'(u') \equiv 1$.
Notice that since $z' (\cdots) = z (\cdots)$, the 
preceding normalization
is preserved, namely $F'(z', 0, u') \equiv 0$.

\proof
Expanding:
\[
\varphi\big(u+i\,F(z,\overline{z},u)\big)
\,=\,
\varphi\big(u+i\,z\overline{z}\,(\cdots)\big)
\,=\,
\varphi(u)
+
z\overline{z}\,\big(\cdots\big),
\]
the fundamental identity writes:
\[
\footnotesize
\!\!\!\!\!\!\!\!\!\!\!\!\!\!\!
\aligned
0
&
\,\equiv\,
-\,
F\big(z,\overline{z},u\big)
+
F'
\Big(
z\,\varphi\big(u+i\,F(z,\overline{z},u)\big),\,\,
\overline{z}\,
\overline{\varphi}\big(u-i\,F(z,\overline{z},u)\big),\,\,
u
\Big)
\\
&
\,\equiv\,
-\,
z\overline{z}\,
F_{1,1}(u)
+
z^2\overline{z}\,
\big(\cdots\big)
+
z\overline{z}^2\,
\big(\cdots\big)
+
z\,
\big(
\varphi(u)
+
z\overline{z}\,(\cdots)
\big)\,
\overline{z}\,
\big(
\overline{\varphi(u)}
+
z\overline{z}\,(\cdots)
\big)\,
F_{1,1}'(u)
+
z^2\overline{z}\,
\big(\cdots\big)
+
\overline{z}^2z\,
\big(\cdots\big)
\\
&
\,\equiv\,
z\overline{z}\,
\Big[
-
F_{1,1}(u)
+
\varphi(u)\,\overline{\varphi}(u)\,
F_{1,1}'(u)
\Big]
+
z^2\overline{z}\,\big(\cdots\big)
+
\overline{z}^2z\,\big(\cdots\big).
\endaligned
\]
To have $F_{1,1}'(u) \equiv 1$, it suffices to take:
\[
\varphi(u)
\,:=\,
\sqrt{F_{1,1}(u)}
\eqno
{\scriptstyle{(\text{\rm remind}\,\,F_{1,1}(0)\,=\,1)}}, 
\]
which is real on the $u$-axis, and then to define
$\varphi(w) := \varphi(u) \big\vert_{u:=w}$, replacing
$u$ by $w$ in the (converging)
power series of $\varphi$.
\endproof

Thus, erasing the primes, still with $0 \equiv F(z, 0, u)$, we have:
\[
v
\,=\,
F\big(z,\overline{z},u\big)
\,=\,
z\overline{z}
+
\sum_{j+k\geqslant 3
\atop
j\geqslant1,\,k\geqslant 1}\,
z^j\overline{z}^k\,
F_{j,k}(u).
\]

\begin{Lemma}
There exists a biholomorphism of the form:
\[
z'
\,:=\,
z
+
\Lambda(z,w)
\,=\,
z
+
z^2\,\big(\cdots\big),
\ \ \ \ \ \ \ \ \ \ \ \ \ \ \ \ \ \ \ \ \ \ \ \ \ \
w'
\,:=\,
w,
\]
which transforms $M = \{ v = F\}$ into $M'$:
\[
v'
\,=\,
F'
\,=\,
z'\overline{z}'
+
\sum_{j\geqslant 2,\,k\geqslant 2}\,
{z'}^j
{\overline{z}'}^k\,
F_{j,k}'(u')
\,\,=\,\,
z'\overline{z}'
+
{z'}^2{\overline{z}'}^2\,
\big(\cdots\big).
\]
\end{Lemma}

Any such biholomorphism with $z' = z + z^2(\cdots)$ preserves the
already achieved normalizations. 

\proof
Single out all monomials with $k = 1$:
\[
\aligned
v
&
\,=\,
z\overline{z}
+
\sum_{j\geqslant 2}\,
z^j\overline{z}^1\,
F_{j,1}(u)
+
\sum_{j+k\geqslant 3
\atop
j\geqslant 1,\,k\geqslant 2}\,
z^j\overline{z}^k\,
F_{j,k}(u)
\\
&
\,=\,
\overline{z}\,
\bigg(
z
+
\underbrace{
\sum_{j\geqslant 2}\,
z^j\,
F_{j,1}(u)}_{=:\,\,\Lambda(z,u)}
\bigg)
+
\overline{z}^2\,
\big(\cdots\big).
\endaligned
\]
Expand:
\[
z'
\,=\,
z
+
\Lambda(z,w)
\,=\,
z
+
\Lambda\big(z,\,u+i\,F(z,\overline{z},u)\big)
\,=\,
z
+
\Lambda\big(z,\,u+iz\overline{z}\,(\cdots)\big)
\,=\,
z
+
\Lambda(z,u)
+
z\overline{z}\,\big(\cdots\big),
\]
and get:
\[
v
\,=\,
\overline{z}\,
\Big(
z'
-
z\overline{z}\,\big(\cdots\big)
\Big)
+
\overline{z}^2\,\big(\cdots\big)
\,=\,
\overline{z}\,z'
+
\overline{z}^2\,\big(\cdots\big).
\]

Next, write the inverse as:
\[
z'
+
{z'}^2\,\big(\cdots\big)
\,=\,
z'
+
\Lambda'(z',w')
\,=\,
z,
\]
so that $\overline{z}^2 (\cdots) = {\overline{z}'}^2 (\cdots)$,
and continue:
\[
\aligned
v'
\,=\,
v
\,=\,
\overline{z}\,z'
+
\overline{z}^2\,\big(\cdots\big)
&
\,=\,
\Big(
\overline{z}'
+
\overline{\Lambda}'\big(\overline{z}',\overline{w}'\big)
\Big)\,
z'
+
\overline{z}^2\,\big(\cdots\big)
\\
&
\,=\,
\Big(
\overline{z}'
+
{\overline{z}'}^2\,\big(\cdots\big)
\Big)\,z'
+
{\overline{z}'}^2\,\big(\cdots\big)
\\
&
\,=\,
z'\overline{z}'
+
{\overline{z}'}^2\,\big(\cdots\big).
\endaligned
\]
The remainder after $z' \overline{z}'$ being real,
it must be also a multiple of ${z'}^2$.
\endproof

\Section{\bf Complete Moser Normal Form for Hypersurfaces 
$M^3 \subset \C^2$}
\label{Moser-normal-form-M3-C2}
\HEAD{{\ref{Moser-normal-form-M3-C2}}.~{\sf Complete Moser Normal 
Form for Hypersurfaces $M^3 \subset \C^2$}
}{
Joël {\sc Merker}}

Thus:
\leqnomode\usetagform{default}
\begin{align}
\label{v-F-2-2-u}
v
\,=\,
z\overline{z}
+
z^2\overline{z}^2\,
F_{2,2}(u)
+
\sum_{j+k\geqslant 5
\atop
j\geqslant 2,\,k\geqslant 2}\,
z^j\overline{z}^k\,
F_{j,k}(u).
\end{align}

\begin{Lemma}
There exists a biholomorphism of the form:
\[
z'
\,:=\,
z\,\lambda(w),
\ \ \ \ \ \ \ \ \ \ \ \ \ \ \ \ \ \ \ \ \ \ \ \ \ \
w'
\,:=\,
w,
\]
with $\lambda(u)\, \overline{\lambda(u)} \equiv 1$ and
$\lambda(0) = 1$, 
such that the new $M'$ has vanishing $F_{2,2}'(u') \equiv 0$:
\[
v'
\,=\,
z'\overline{z}'
+
0
+
\sum_{j+k\geqslant 5
\atop
j\geqslant 2,\,k\geqslant 2}\,
{z'}^j{\overline{z}'}^k\,
F_{j,k}'(u').
\]
\end{Lemma}

The condition $\big\vert \lambda(u) \big\vert^2 \equiv 1$
for $w = u \in \R$ guarantees that all the previously achieved
normalizations are preserved.

\proof
Expand:
\[
\aligned
\lambda\big(u+i\,F(z,\overline{z},u)\big)
\,=\,
\lambda\big(u+i\,z\overline{z}+z^2\overline{z}^2\,(\cdots)\big)
&
\,=\,
\lambda(u)\,
+
\lambda_u(u)\,
\big[
i\,z\overline{z}
+
z^2\overline{z}^2\,(\cdots)
\big]
+
z^2\overline{z}^2\,
\big(\cdots\big)
\\
&
\,=\,
\lambda(u)\,
\Big(
1
+
\frac{\lambda_u(u)}{\lambda(u)}\,
i\,z\overline{z}
+
z^2\overline{z}^2\,(\cdots)
\Big).
\endaligned
\]
Since we assume $\big\vert \lambda(u) \big\vert^2 \equiv 1$,
{\em i.e.} $\lambda(u) = e^{i\,\varphi(u)}$ with $\varphi(u)$
real, the quotient $\frac{\lambda_u(u)}{\lambda(u)}$ is purely
imaginary, hence:
\[
\big\vert
\lambda(u+i\,F)
\big\vert^2
\,=\,
1
+
2i\,z\overline{z}\,
\frac{\lambda_u(u)}{\lambda(u)}
+
z^2\overline{z}^2\,
\big(\cdots\big).
\]
Also, it is clear that ${z'}^j {\overline{z}'}^k\, (\cdots) 
= z^j \overline{z}^k\, (\cdots)$.

Thanks to these preliminaries:
\[
\aligned
v'
&
\,=\,
z'\overline{z}'
+
{z'}^2{\overline{z}'}^2\,
F_{2,2}'(u')
+
{z'}^3{\overline{z}'}^2\,
\big(\cdots\big)
+
{z'}^2{\overline{z}'}^3\,
\big(\cdots\big)
\\
&
\,=\,
\big\vert
\lambda(u+i\,F)
\big\vert^2\,
z\overline{z}
+
\big\vert
\lambda(u+i\,F)
\big\vert^4\,
z^2\overline{z}^2\,
F_{2,2}'(u)
+
z^3\overline{z}^2\,
\big(\cdots\big)
+
z^2\overline{z}^3\,
\big(\cdots\big)
\\
&
\,=\,
z\overline{z}
+
z^2\overline{z}^2\,
2i\,
\frac{\lambda_u(u)}{\lambda(u)}
+
z^3\overline{z}^3\,
\big(\cdots\big)
+
z^2\overline{z}^2\,
\big(1+z\overline{z}\,(\cdots)\big)\,
F_{2,2}'(u)
+
z^3\overline{z}^2\,
\big(\cdots\big)
+
z^2\overline{z}^3\,
\big(\cdots\big)
\\
&
\,=\,
z\overline{z}
+
z^2\overline{z}^2\,
\Big[
2i\,
\frac{\lambda_u(u)}{\lambda(u)}
+
F_{2,2}'(u)
\Big]
+
z^3\overline{z}^2\,
\big(\cdots\big)
+
z^2\overline{z}^3\,
\big(\cdots\big),
\endaligned
\]
and since $v' = v$ with $v$ given by~({\ref{v-F-2-2-u}}), 
an identification yields:
\[
2i\,
\frac{\lambda_u(u)}{\lambda(u)}
+
F_{2,2}'(u)
\,\equiv\,
F_{2,2}(u).
\]
In order to annihilate $F_{2,2}'(u') := 0$, it suffices
therefore to set:
\[
\lambda(u)
\,:=\,
\exp\,
\bigg(
\frac{1}{2i}
\int_0^u\,
F_{2,2}(t)\,dt
\bigg).
\qedhere
\]
\endproof

Now we come to a crucial moment offering a key simplification
which was prepared in advance
by Assertion~{\ref{Assertion-key-F-3-2-zero}}.

\begin{Assertion}
\label{Assertion-crucial-F-3-2-zero}
After having normalized:
\[
0
\,\equiv\,
F_{j,0}(u)
\,\equiv\,
F_{0,k}(u)
\,\equiv\,
\underset{(j\neq 1)}{F_{j,1}(u)}
\,\equiv\,
\underset{(1\neq k)}{F_{1,k}(u)},
\ \ \ \ \ \ \ \ \ \ \ \ \ \ \ \ \ \ \ \
1
\,\equiv\,
F_{1,1}(u),
\ \ \ \ \ \ \ \ \ \ \ \ \ \ \ \ \ \ \ \
0
\,\equiv\,
F_{2,2}(u),
\]
the fact that the $u$-axis, contained in $M$, is a Moser chain,
offers without any further work:
\[
0
\,\equiv\,
F_{3,2}(u)
\,\equiv\,
F_{2,3}(u).
\]
\end{Assertion}

\proof
At each point $p = (0, u_p) \in M$ with any (small) $u_p \in \R$
in the straightened Moser chain, the equation of $M$ normalized
up to this point and truncated after
weighted order $6$ writes exactly:
\[
v
\,=\,
z\overline{z}
+
z^3\overline{z}^2\,
F_{3,2}(u_p)
+
z^2\overline{z}^3\,
F_{2,3}(u_p)
+
{\rm O}(6),
\] 
under the form considered in 
Assertion~{\ref{Assertion-key-F-3-2-zero}},
which then yields $F_{3,2}(u_p) = 0 = F_{2,3}(u_p)$,
this for any $u_p$.
\endproof

Thus:
\[
v
\,=\,
z\overline{z}
+
z^4\overline{z}^2\,
F_{4,2}(u)
+
z^3\overline{z}^3\,
F_{3,3}(u)
+
z^2\overline{z}^4\,
F_{2,4}(u)
+
\sum_{j+k\geqslant 7
\atop
j\geqslant 2,\,k\geqslant 2}\,
z^j\overline{z}^k\,
F_{j,k}(u).
\]

\begin{Lemma}
There exists a biholomorphism of the form:
\[
z'
\,:=\,
z\,\sqrt{\psi_w(w)},
\ \ \ \ \ \ \ \ \ \ \ \ \ \ \ \ \ \ \ \ \ \ \ \ \ \
w'
\,:=\,
\psi(w),
\]
with $\psi(\R) \subset \R$, 
with $\psi(0) = 0$, with $\psi_w(0) \in \R_{>0}$, such that the
new $M'$ has vanishing $F_{3,3}'(u') \equiv 0$:
\[
v'
\,=\,
z'\overline{z}'
+
{z'}^4{\overline{z}'}^2\,
F_{4,2}'(u')
+
0
+
{z'}^2{\overline{z}'}^4\,
F_{2,4}'(u')
+
\sum_{j+k\geqslant 7
\atop
j\geqslant 2,\,k\geqslant 2}\,
{z'}^j{\overline{z}'}^k\,
F_{j,k}'(u').
\]
\end{Lemma}

We will see in the proof why such a biholomorphism preserves
all previously achieved normalizations.
The function $\psi = \psi(u)$ will be solution of the
ODE:
\[
\psi_{uuu}(u)
\,=\,
\frac{3}{2}\,
\frac{\psi_{uu}^2(u)}{\psi_u(u)}
-
3\,F_{3,3}(u)\,
\psi_u(u).
\]

\proof
More generally, we perform a biholomorphism of the form:
\[
z'
\,:=\,
z\,\varphi(w),
\ \ \ \ \ \ \ \ \ \ \ \ \ \ \ \ \ \ \ \ \ \ \ \ \ \
w'
\,:=\,
\psi(w),
\]
assuming that $\varphi(u) \in \R$, $\varphi(0) \neq 0$, 
and $\psi(u) \in \R$, $\psi_w(0) \in \R_{\neq 0}$. We let
$v' = F'(z', \overline{z}', u')$ be the transformed hypersurface
equation. Many computations are needed.

Firstly:
\[
\!\!\!\!\!\!\!\!\!\!\!\!\!\!\!\!\!\!\!\!\!\!\!\!\!
\scriptsize
\aligned
v'
\,=\,
\Im\,\psi\big(u+i\,F\big)
&
\,=\,
\Im\,
\Big\{
\psi(u)
+
\psi_u(u)\,i\,F
+
\psi_{uu}(u)\,
\frac{(i\,F)^2}{2!}
+
\psi_{uuu}(u)\,
\frac{(i\,F)^3}{3!}
+
F^4\,\big(\cdots\big)
\Big\}
\\
&
\,=\,
\psi_u(u)\,F
-
\frac{1}{6}\,
\psi_{uuu}(u)\,F^3
+
z^4\overline{z}^4\,
\big(\cdots\big)
\\
&
\,=\,
\psi_u(u)\,
\Big[
z\overline{z}
+
z^4\overline{z}^2\,
F_{4,2}(u)
+
z^3\overline{z}^3\,
F_{3,3}(u)
+
z^2\overline{z}^4\,
F_{2,4}(u)
+
{\rm O}_{z,\overline{z}}(7)
\Big]
-
\frac{1}{6}\,
\psi_{uuu}(u)\,
\big[
z^3\overline{z}^3
+
{\rm O}_{z,\overline{z}}(10)
\big]
+
z^4\overline{z}^4\,
\big(\cdots\big),
\endaligned
\]
so that no terms of order $3$, $4$, $5$ in $(z, \overline{z})$
are present:
\leqnomode\usetagform{default}
\begin{align}
\label{v-prime-last-lemma}
\!\!\!\!\!\!\!\!\!\!
v'
\,=\,
z\overline{z}\,
\psi_u(u)
+
z^4\overline{z}^2\,
\psi_u(u)\,F_{4,2}(u)
+
z^3\overline{z}^3\,
\big[
\psi_u(u)\,
F_{3,3}(u)
-
{\textstyle{\frac{1}{6}}}\,
\psi_{uuu}(u)
\big]
+
z^2\overline{z}^4\,
\psi_u(u)\,F_{2,4}(u)
+
{\rm O}_{z,\overline{z}}(7).
\end{align}

Secondly, one can convince oneself that the normalization
$v' = z'\overline{z}' + {z'}^2 {\overline{z}'}^2 \,
\big( \cdots \big)$ is preserved, so that the equation of the
transformed hypersurface is:
\[
v'
\,=\,
z'\overline{z}'
+
\sum_{j\geqslant 2,\,k\geqslant 2}\,
{z'}^j{\overline{z}'}^k\,
F_{j,k}'(u').
\]

Thirdly, using $\varphi(u) \in \R$ and $F = z\overline{z} + 
{\rm O}_{z,\overline{z}}(6)$:
\[
\footnotesize
\aligned
z'\overline{z}'
&
\,=\,
z\overline{z}
\Big(
\varphi(u)
+
\varphi_u(u)\,
i\,F
+
\varphi_{uu}(u)\,
\frac{(i\,F)^2}{2!}
+
F^3\,\big(\cdots\big)
\Big)\,
\Big(
\varphi(u)
+
\varphi_u(u)\,\big(-i\,F\big)
+
\varphi_{uu}(u)\,
\frac{(-i\,F)^2}{2!}
+
F^3\,\big(\cdots\big)
\Big)
\\
&
\,=\,
z\overline{z}\,
\varphi(u)^2
+
z\overline{z}\,
\Big[
\varphi_u(u)^2
-
\varphi(u)\,\varphi_{uu}(u)
\Big]\,
F^2
+
z\overline{z}\,
F^3\,\big(\cdots\big)
\\
&
\,=\,
z\overline{z}\,
\varphi(u)^2
+
z^3\overline{z}^3\,
\big(
\varphi_u(u)^2
-
\varphi(u)\,\varphi_{uu}(u)
\big)
+
{\rm O}_{z,\overline{z}}(8).
\endaligned
\]

Fourthly, for every $j \geqslant 2$ and every $k \geqslant 2$:
\[
\aligned
{z'}^j{\overline{z}'}^k
&
\,=\,
z^j\overline{z}^k\,
\Big(
\varphi(u)
+
\varphi_u(u)\,i\,F
+
z^2\overline{z}^2\,(\cdots)
\Big)^j\,
\Big(
\varphi(u)
+
\varphi_u(u)\,(-i\,F)
+
z^2\overline{z}^2\,(\cdots)
\Big)^k
\\
&
\,=\,
z^j\overline{z}^k\,
\Big(
\varphi(u)^j
+
j\,\varphi(u)^{j-1}\varphi_u(u)\,i\,z\overline{z}
+
z^2\overline{z}^2\,(\cdots)
\Big)\,
\Big(
\varphi(u)^k
-
k\,\varphi(u)^{k-1}\varphi_u(u)\,
i\,z\overline{z}
+
z^2\overline{z}^2\,(\cdots)
\Big)
\\
&
\,=\,
z^j\overline{z}^k\,
\Big(
\varphi(u)^{j+k}
+
i\,(j-k)\,\varphi(u)^{j+k-1}\,
\varphi_u(u)\,
z\overline{z}
+
z^2\overline{z}^2\,
(\cdots)
\Big).
\endaligned
\]

Fifthly:
\[
\aligned
F_{j,k}'(u')
\,=\,
F_{j,k}'
\Big(
\Re\,\psi\big(u+i\,F\big)
\Big)
&
\,=\,
F_{j,k}'
\Big(
\Re\,
\big[
\psi(u)
+
\psi_u(u)\,i\,F
+
F^2\,(\cdots)
\big]
\Big)
\\
&
\,=\,
F_{j,k}'
\Big(
\psi(u)
+
0
+
z^2\overline{z}^2\,
\big(\cdots\big)
\Big)
\\
&
\,=\,
F_{j,k}'\,
\big(
\psi(u)
\big)
+
z^2\overline{z}^2\,
\big(\cdots\big).
\endaligned
\]

Thanks to all this:
\[
\!\!\!\!\!\!\!\!\!\!\!\!\!\!\!
\!\!\!\!\!\!\!\!\!\!\!\!\!\!\!
\scriptsize
\aligned
F\big(z',\overline{z}',u'\big)
&
\,=\,
z'\overline{z}'
+
{z'}^2{\overline{z}'}^2\,
F_{2,2}'(u')
+
{z'}^3{\overline{z}'}^2\,
F_{3,2}'(u')
+
{z'}^2{\overline{z}'}^3\,
F_{2,3}'(u')
\,+
\\
&
\ \ \ \ \ 
\ \ \ \ \ \ \ \ \ \ \ \ \ \ \ \ \ \ \ \ \ \ \ \ \ \
+
{z'}^4{\overline{z}'}^2\,
F_{4,2}'(u')
+
{z'}^3{\overline{z}'}^3\,
F_{3,3}'(u')
+
{z'}^2{\overline{z}'}^4\,
F_{2,4}'(u')
+
{\rm O}_{z',\overline{z}'}(7)
\\
&
\,=\,
z\overline{z}\,\varphi(u)^2
+
z^3\overline{z}^3\,
\big(
\varphi_u(u)^2
-
\varphi(u)\,\varphi_{uu}(u)
\big)
+
{\rm O}_{z,\overline{z}}(8)
\,+
\\
&
\ \ \ \ \ 
+
z^2\overline{z}^2\,
\Big(
\varphi(u)^4
+
0
+
z^2\overline{z}^2\,(\cdots)
\Big)\,
\Big(
F_{2,2}'\big(\psi(u)\big)
+
z^2\overline{z}^2\,(\cdots)
\Big)
\\
&
\ \ \ \ \ 
+
z^3\overline{z}^2\,
\Big(
\varphi(u)^5
+
i\,\varphi(u)^4\,\varphi_u(u)\,z\overline{z}
+
z^2\overline{z}^2\,(\cdots)
\Big)\,
F_{3,2}'\big(\psi(u)\big)
+
z^2\overline{z}^3\,
\Big(
\varphi(u)^5
-
i\,
\varphi(u)^4\,\varphi_u(u)\,z\overline{z}
+
z^2\overline{z}^2\,(\cdots)
\Big)\,
F_{2,3}'\big(\psi(u)\big)
\\
&
\ \ \ \ \
+
z^4\overline{z}^2\,
\varphi(u)^6\,
F_{4,2}'\big(\psi(u)\big)
+
z^3\overline{z}^3\,
\varphi(u)^6\,
F_{3,3}'\big(\psi(u)\big)
+
z^2\overline{z}^4\,
\varphi(u)^6\,
F_{2,4}'\big(\psi(u)\big)
+
{\rm O}_{z,\overline{z}}(7)
\\
&
\,=\,
z\overline{z}\,
\varphi(u)^2
+
z^2\overline{z}^2\,\varphi(u)^4\,
F_{2,2}'\big(\psi(u)\big)
\,+
\\
&
\ \ \ \ \ \ \ \ \ \ \ \ \ \ \ \ \ \ \ \ \ \,
+
z^3\overline{z}^2\,\varphi(u)^5\,F_{3,2}'\big(\psi(u)\big)
+
z^2\overline{z}^3\,\varphi(u)^5\,F_{2,3}'\big(\psi(u)\big)
\,+
\\
&
\ \ \ \ \ \ \ \ \ \ \ \ \ \ \ \ \ \ \ \ \ \,
+
z^4\overline{z}^2\,\varphi(u)^6\,F_{4,2}'\big(\psi(u)\big)
+
z^3\overline{z}^3\,
\Big[
\varphi_u(u)^2
-
\varphi(u)\,\varphi_{uu}(u)
+
\varphi(u)^6\,F_{3,3}'\big(\psi(u)\big)
\Big]
+
z^2\overline{z}^4\,\varphi(u)^6\,F_{2,4}'\big(\psi(u)\big)
+
{\rm O}_{z,\overline{z}}(7).
\endaligned
\]
By identifying powers $z^j \overline{z}^k$ 
with~({\ref{v-prime-last-lemma}}), we get:
\leqnomode\usetagform{default}
\begin{align}
\psi_u(u)
&
\,\equiv\,
\varphi(u)^2,
\tag{1,1}
\\
0
&
\,\equiv\,
\varphi(u)^4\,F_{2,2}'\big(\psi(u)\big),
\tag{2,2}
\\
0
&
\,\equiv\,
\varphi(u)^5\,F_{3,2}'\big(\psi(u)\big),
\tag{3,2}
\\
0
&
\,\equiv\,
\varphi(u)^5\,F_{2,3}'\big(\psi(u)\big),
\tag{2,3}
\\
\psi_u(u)\,F_{4,2}(u)
&
\,\equiv\,
\varphi(u)^6\,F_{4,2}'\big(\psi(u)\big),
\tag{4,2}
\\
\psi_u(u)\,F_{3,3}(u)
-
{\textstyle{\frac{1}{6}}}\,
\psi_{uuu}(u)
&
\,\equiv\,
\varphi_u(u)^2
-
\varphi(u)\,\varphi_{uu}(u)
+
\varphi(u)^6\,F_{3,3}'\big(\psi(u)\big),
\tag{3,3}
\\
\psi_u(u)\,F_{2,4}(u)
&
\,\equiv\,
\varphi(u)^6\,F_{2,4}'\big(\psi(u)\big).
\tag{2,4}
\end{align}

Visibly, to annihilate $F_{3,3}'(u')$, it suffices to
fulfill:
\[
\aligned
\psi_u(u)
&
\,\equiv\,
\varphi(u)^2,
\\
\psi_u(u)\,F_{3,2}(u)
-
{\textstyle{\frac{1}{6}}}\,
\psi_{uuu}(u)
&
\,\equiv\,
\varphi_u(u)^2
-
\varphi(u)\,\varphi_{uu}(u)
+
0.
\endaligned
\]
Assuming $\psi_u(0) = 1$, choosing $\varphi(u) := \sqrt{ \psi_u(u)}$,
and replacing, it suffices in conclusion that $\psi$ satisfies the
solvable ODE:
\[
\psi_{uuu}(u)
\,=\,
\frac{3}{2}\,
\frac{\psi_{uu}(u)^2}{\psi_u(u)}
-
3\,F_{3,3}(u)\,
\psi_u(u).
\qedhere
\]
\endproof

In summary, we have fully reproved with expository details what 
is actually the equation~(3.18) of Chern-Moser's
celebrated work~{\cite{Chern-Moser-1974}}.

\begin{Proposition}
\label{Prp-final-normal-form}
Given a Levi nondegenerate $\mathcal{C}^\omega$ hypersurface $M^3
\subset \C^2$, for every $p \in M$ and every CR-transversal $1$-jet
$j_p^1$ at $p$, if $\gamma_p \ni p$ denotes the unique piece of Moser
chain directed by $j_p^1$ at $p$, then there exist local holomorphic
coordinates $(z, w = u + i\, v)$ centered at $p$ in which $\gamma_p$
is the $u$-axis and such that $M$ is graphed as:
\[
v
\,=\,
z\overline{z}
+
z^4\overline{z}^2\,F_{4,2}(u)
+
z^2\overline{z}^4\,F_{2,4}(u)
+
\sum_{j+k\geqslant 7
\atop j\geqslant 2,\,k\geqslant 2}\,
z^j\overline{z}^k
F_{j,k}(u).
\eqno\qed
\]
\end{Proposition}

\Section{\bf Uniqueness of Moser Normal Form}
\label{uniqueness-Moser-normal-form}
\HEAD{{\ref{uniqueness-Moser-normal-form}}.~{\sf Uniqueness
of Moser Normal Form}
}{
Joël {\sc Merker}}

Starting with a $\mathcal{C}^\omega$ Levi nondegenerate hypersurface
$M^3 \subset \C^2$, at any point $p \in M$, it is elementary to find
holomorphic coordinates $(z, w)$ vanishing at $p$ in which $M$ has
equation $v = F = z\overline{z} + {\rm O}_{z, \overline{z}, u}(3)$.
Such an equation can hence freely 
be taken as the starting point towards a complete
normalization of $F$.

In the preceding sections, we have in fact established the {\em
existence} of a {\sl normal form} for $M$. We can now present 
the known
{\sl uniqueness} statement.

\begin{Theorem}
{\cite{Chern-Moser-1974, Jacobowitz-1990}}
Given a $\mathcal{C}^\omega$ 
Levi nondegenerate hypersurface $M^3 \subset \C^2$
with $0 \in M$ of the form:
\[
v
\,=\,
z\overline{z}
+
{\rm O}_{z,\overline{z},u}(3),
\]
there exists a biholomorphism $(z,w) \longmapsto (z', w')$
fixing $0$ which maps $(M, 0)$ into $(M', 0)$ of 
normalized equation:
\[
v'
\,=\,
z'\overline{z}'
+
F_{4,2}'(u')\,
{z'}^4{\overline{z}'}^2
+
F_{2,4}'(u')\,
{z'}^2{\overline{z}'}^4
+
{z'}^2{\overline{z}'}^2\,
{\rm O}_{z',\overline{z}'}(3).
\]

Furthermore, the map exists and is {\em unique} 
if it is assumed to be of the form:
\[
\aligned
z'
\,:=\,
&\
z
+
f(z,w),
\ \ \ \ \ \ \ \ \ \ \ \ \ \ \ \ \ \ \ \
w'
\,:=\,
w
+
g(z,w),
\\
f_z(0)
\,=\,
&\,
f_w(0)
\,=\,
0,
\ \ \ \ \ \ \ \ \ \ \ \ \ \ \ \ \ \ \ \
g_z(0)
\,=\,
g_w(0)
\,=\,
\Re\,g_{ww}(0)
\,=\,
0.
\endaligned
\]
\end{Theorem}

\proof
By choosing a chain at $0 \in M$ whose first jet is flat,
directed along the $u$-axis, one can verify (exercise) that
all the constructions done in the preceding sections do 
indeed give a biholomorphism of this specific form.
So our job is to establish uniqueness.

\smallskip

Suppose that two such normalizations $h_\iota \colon (z,w)
\longmapsto (z+f_\iota,\, w+g_\iota)$, $\iota = 1, 2$, 
are given:
\[
\xymatrix{
& & &
M_1'
\ar[dd]^{h_2\circ h_1^{-1}}
\\
\ar[urrr]^{h_1}
M
\ar[drrr]_{h_2}
\\
& & &
M_2',}
\]
with $0 = f_{\iota,z}(0) = f_{\iota,w}(0)$
and $0 = g_{\iota,z}(0) = g_{\iota,w}(0) = {\rm Re}\, 
g_{\iota,ww}(0)$.
On ${\C'}^2 \supset M_1'$, let us take 
for simplicity coordinates with the same name
$(z,w)$, and coordinates $(z', w')$ on the ${\C'}^2 \supset M_2'$.

\begin{Assertion}
Then $h_2 \circ h_1^{-1} =: (z + f,\, w+g)$ also satisfies
$0 = f_z(0) = f_w(0)$ and $0 = g_z(0) = g_w(0) = 
{\rm Re}\, g_{ww}(0)$.
\end{Assertion}

\proof
Since both $h_1$ and $h_2$ are the identity plus 
${\rm O}_{z,w}(2)$ terms,
the same holds for $h_2 \circ h_1^{-1}$. 
It remains only to show $\Re\, g_{ww}(0) = 0$.

The following lemma then applies to the map $h_2 \circ h_1^{-1}$,
since $M_1'$ and $M_2'$ are in normal form.

\begin{Lemma}
\label{Lm-g-ww-0-real}
If $(z,w) \longmapsto (z+f,\, w+g)$ with $f, g = {\rm O}_{z,w}(2)$,
maps $v = z\overline{z} + {\rm O}_{z, \overline{z}, u}(3)$
to $v' = z'\overline{z}' + {\rm O}_{z', \overline{z}', u'}(3)$,
then $g_{zz}(0) = g_{zw}(0) = 0$ and $g_{ww}(0) \in \R$, so that:
\[
g(z,w)
\,=\,
w
+
\tfrac{1}{2}\,
g_{ww}(0)\,
w^2
+
{\rm O}_{z,w}(3).
\]
\end{Lemma}

\proof
Writing $w' = w + g = w + \alpha\, z^2 + \beta\, zw + (a+ib)\, w^2 + 
{\rm O}_{z,w}(3)$, we have:
\[
\aligned
v'
&
\,=\,
v
+
\Im\,(\alpha\,z^2)
+
\Im\,\big(
\beta\,z\,(u+iv)
\big)
+
2a\,uv
+
b\,u^2
-
b\,v^2
+
{\rm O}_{z,w}(3)
\\
&
\,=\,
z\overline{z}
+
\Im\,(\alpha\,z^2)
+
\Im\,(\beta\,zu)
+
b\,u^2
+
{\rm O}_{z,\overline{z},u}(3),
\endaligned
\]
hence using the inversion $z = z' + {\rm O}_{z',w'}(2)$,
$w = w' + {\rm O}_{z', w'}(2)$, we get
$\alpha = \beta = b = 0$ from:
\[
v'
\,=\,
z'\overline{z}'
+
\Im\,(\alpha\,{z'}^2)
+
\Im\,(\beta\,z'u')
+
b\,{u'}^2
+
{\rm O}_{z',\overline{z}',u'}(3).
\qedhere
\]
\endproof

Thus, the assumption 
$\Re\, g_{\iota, ww}(0) = 0$, $\iota = 1, 2$,
implies that
the $h_\iota$, are both of the form
$\big(z + {\rm O}_{z,w}(2),\,\, w + {\rm O}_{z,w}(3) \big)$.
Such a form is stable under composition and inversion,
hence $h_2 \circ h_1^{-1}$ is also of this form,
and in particular, one has $\Re\, g_{ww}(0) = 0$.
\endproof

Our uniqueness goal is to obtain $h_1 = h_2$. Equivalently, 
$h_2 \circ h_1^{-1} = \Id$. This will be offered by 
the next independent key uniqueness statement.
\endproof

\begin{Theorem}
\label{Thm-uniqueness-normal-form}
If two $\mathcal{C}^\omega$ Levi nondegenerate hypersurfaces
$0 \in M^3 \subset \C^2$ and $0 \in {M'}^3 \subset {\C'}^2$
are both in normal form:
\[
\aligned
v
&
\,=\,
F
\,=\,
z\overline{z}
+
z^4\overline{z}^2\,F_{4,2}(u)
+
z^2\overline{z}^4\,F_{2,4}(u)
+
\sum_{j+k\geqslant 7
\atop j\geqslant 2,\,k\geqslant 2}\,
z^j\overline{z}^k
F_{j,k}(u),
\\
v'
&
\,=\,
F'
\,=\,
z'\overline{z}'
+
{z'}^4{\overline{z}'}^2\,F_{4,2}'(u')
+
{z'}^2{\overline{z}'}^4\,F_{2,4}'(u')
+
\sum_{j+k\geqslant 7
\atop j\geqslant 2,\,k\geqslant 2}\,
{z'}^j{\overline{z}'}^k
F_{j,k}'(u'),
\endaligned
\]
and if there exists a biholomorphism $(M, 0) \longrightarrow 
(M', 0)$ of the form:
\[
\aligned
z'
\,:=\,
&\
z
+
f(z,w),
\ \ \ \ \ \ \ \ \ \ \ \ \ \ \ \ \ \ \ \
w'
\,:=\,
w
+
g(z,w),
\\
f_z(0)
\,=\,
&\,
f_w(0)
\,=\,
0,
\ \ \ \ \ \ \ \ \ \ \ \ \ \ \ \ \ \ \ \
g_z(0)
\,=\,
g_w(0)
\,=\,
\Re\,g_{ww}(0)
\,=\,
0,
\endaligned
\]
then $(f,g) \equiv (0,0)$, and the biholomorphism is the identity.
\end{Theorem}

\proof
Equivalently, the graphing function $F = \sum_{j,k}\, 
F_{j,k}(u)\, z^j \overline{z}^k$ of $M$ satisfies the 
general {\sl prenormalization conditions:}
\[
0
\,\equiv\,
F_{j,0}(u)
\,\equiv\,
F_{0,k}(u),
\ \ \ \ \ \ \ \ \ \ \ \ \ \ \ \ \ \ \ \
0
\,\equiv\,
F_{j,1}(u)
\,\equiv\,
F_{1,k}(u)
\eqno
{\scriptstyle{(j,\,k\,\in\,\N)}},
\]
except of course $1 \equiv F_{1,1}(u)$, together with the
{\sl sporadic normalization conditions:}
\[
0
\,\equiv\,
F_{2,2}(u)
\,\equiv\,
F_{3,2}(u)
\,\equiv\,
F_{2,3}(u)
\,\equiv\,
F_{3,3}(u),
\]
and the same holds about $F'$.

Accordingly, let us introduce:
\[
S
\,:=\,
\big\{
(j,0),\,\,
(0,k),\,\,
(j,1),\,\,
(1,k)
\big\}
\cup
\big\{
(2,2),\,\,
(3,2),\,\,
(2,3),\,\,
(3,3)
\big\}.
\]
For a general real converging power series vanishing 
at $(z, \overline{z}, u) = (0, 0, 0)$:
\[
G
\,=\,
\sum_{j,k,l}\,
G_{j,k,l}\,
z^j\overline{z}^ku^l
\eqno
{\scriptstyle{(\overline{G_{k,j,l}}\,=\,G_{j,k,l})}},
\]
{\em i.e.} with $G_{0,0,0} = 0$, introduce the projection:
\[
\Pi_S(G)
\,:=\,
\sum_{(j,k)\in S}\,
\sum_{l=0}^\infty\,
G_{j,k,l}\,
z^j\overline{z}^ku^l,
\]
so that:
\[
\Pi_S(F)
\,=\,
z\overline{z}
\ \ \ \ \ \ \ \ \ \ \ \ \ \ \ \ \ \ \ \
\text{and}
\ \ \ \ \ \ \ \ \ \ \ \ \ \ \ \ \ \ \ \
\Pi_S(F')
\,=\,
z'\overline{z}'.
\]

Also, reminding that granted our current assumption $\Re\, g_{ww}(0) =
0$, we already understood in Lemma~{\ref{Lm-g-ww-0-real}} that we have
in fact $g = w + {\rm O}_{z,w}(3)$.  Next, taking integers $\nu
\geqslant 3$, reminding weights $[z] = 1$, 
$[w] = 2$, 
let us decompose in weighted homogeneous components:
\[
\aligned
f(z,w)
\,=\,
\sum_{j+l\geqslant 2}\,
f_{j,l}\,
z^jw^l
\,=\,
&\,
\sum_{\nu\geqslant 3}\,
f_{\nu-1},
\ \ \ \ \ \ \ \ \ \
&
\ \ \ \ \ \ \ \ \ \
g(z,w)
\,=\,
\sum_{j+l\geqslant 3}\,
g_{j,l}\,z^jw^l
\,=\,
&\,
\sum_{\nu\geqslant 3}\,
g_\nu,
\\
f_{\nu-1}
\,:=\,
&\,
\sum_{j+2l=\nu-1}\,
f_{j,l}\,z^jw^l
\ \ \ \ \ \ \ \ \ \
&
\ \ \ \ \ \ \ \ \ \
g_\nu
\,:=\,
&\,
\sum_{j+2l=\nu}\,
g_{j,l}\,z^jw^l.
\endaligned
\]
Still for any $\nu \geqslant 3$, introduce the projections:
\[
\pi_{\nu-1}(f)
\,:=\,
f_{\nu-1},
\ \ \ \ \ \ \ \ \ \ \ \ \ \ \ \ \ \ \ \
\pi_\nu(g)
\,:=\,
g_\nu,
\ \ \ \ \ \ \ \ \ \ \ \ \ \ \ \ \ \ \ \
\pi_\nu(G)
\,:=\,
G_\nu
\,:=\,
\sum_{j+k+2l=\nu}\,
G_{j,k,l}\,z^j\overline{z}^ku^l,
\]
so that:
\[
\Pi_S\big(\pi_\nu(F)\big)
\,=\,
0
\,=\,
\Pi_S\big(\pi_\nu(F')\big)
\eqno
{\scriptstyle{(\nu\,\geqslant\,3)}}.
\]
Also, introduce:
\[
\pi^\nu
\,:=\,
\pi_2+\cdots+\pi_\nu.
\]
For later use, observe that
for any holomorphic function $e_\mu = e_\mu(z,w)$ which is
weigthed $\mu$-homogeneous, it holds (exercise):
\leqnomode\usetagform{default}
\begin{align}
\label{pi-mu-truncate-z-zbar}
\pi^\mu\,
\Big(
e_\mu
\big(
z,\,
u+i\,[z\overline{z}+{\rm O}_{z,\overline{z},u}(3)]
\big)
\Big)
\,\,=\,\,
e_\mu\big(z,\,u+iz\overline{z}\big).
\end{align}

Next, since $f = f_2 + f_3 + \cdots$ and $g = g_3 + g_4 + \cdots$,
the fundamental identity writes:
\[
0
\,\equiv\,
-\,\Im\,
\big(
w+g_3+g_4+\cdots
\big)
+
F'
\Big(
z+f_2+f_3+\cdots,\,\,
\overline{z}+\overline{f}_2+\overline{f}_3+\cdots,\,\,
\Re
\big(
w+g_3+g_4+\cdots
\big)
\Big),
\]
identically in $\C\{z, \overline{z}, u\}$
after replacing $(z,w) = \big(z, u + i\, F(z, \overline{z}, u) \big)$.

To prove $(f, g) = (0, 0)$, we may proceed progressively:

\smallskip\noindent$\bullet$\,
$(f_2, g_3) = (0,0)$;

\smallskip\noindent$\bullet$\,
$(f_3, g_4) = (0,0)$;

\smallskip\noindent$\bullet$\,
$(f_{\mu-1},\, g_\mu) = (0, 0)$ for $\mu = 3, \dots, \nu-1$
and some $\nu \geqslant 5$ implies $(f_{\nu-1},\, g_\nu) = (0,0)$.

\begin{Assertion}
One has $(f_2, g_3) = (0, 0)$.
\end{Assertion}

\proof
Applying $\pi^3$ to the fundamental identity gives,
using~({\ref{pi-mu-truncate-z-zbar}}):
\[
\aligned
0
&
\,\equiv\,
\pi^3
\bigg(
-\,\Im\,\big(w+g_3\big)
+
F'\Big(
z+f_2,\,\overline{z}+\overline{f}_2,\,
\Re\,\big(w+g_2\big)
\Big)
\bigg)
\\
&
\,\equiv\,
\pi^3
\bigg(
-\,v
-
\Im\,g_3
+
\big(z+f_2\big)\,
\big(\overline{z}+\overline{f}_2\big)
+
F_3'
\big(
z,\overline{z},u
\big)
\bigg)
\\
&
\,\equiv\,
\pi^3
\bigg(
-\,
\zerozero{z\overline{z}}
-
\zero{
F_3(z,\overline{z},u)}
-
\Im\,g_3
+
\zerozero{z\overline{z}}
+
z\overline{f}_2
+
\overline{z}f_2
+
f_2\overline{f}_2
+
\zero{
F_3'\big(z,\overline{z},u\big)}
\bigg),
\endaligned
\]
and since $M$ and $M'$ are normalized by assumption,
with $\pi^3(f_2 \overline{f}_2) \equiv 0$,
it remains only:
\[
0
\,\equiv\,
\Re\,
\Big\{
i\,g_3\big(z,\,u+iz\overline{z}\big)
+
2\,\overline{z}\,
f_2\big(z,\,u+iz\overline{z}\big)
\Big\}.
\]

Replacing $f_2 = f_{2,0}\, z^2 + f_{0,1}\, w$
with $f_{0,1} = 0$ by assumption and replacing $g_3 = 
g_{3,0}\, z^3 + g_{1,1}\, zw$, this is:
\[
0
\,\equiv\,
\tfrac{i}{2}\,
g_{3,0}\,z^3
-
\tfrac{i}{2}\,
\overline{g}_{3,0}\,
\overline{z}^3
+
\big(
f_{2,0}
-
\tfrac{1}{2}\,
g_{1,1}
\big)\,
z^2\overline{z}
+
\big(
\overline{f}_{2,0}
-
\tfrac{1}{2}\,
\overline{g}_{1,1}
\big)\,
z\overline{z}^2
+
\tfrac{i}{2}\,
g_{1,1}\,
zu
-
\tfrac{i}{2}\,
\overline{g}_{1,1}\,
\overline{z}u,
\]
and starting from the end, this forces $0 = g_{1,1} = f_{2,0} = 
g_{3,0}$, so as asserted $0 = f_2 = g_3$.
\endproof

\begin{Assertion}
One has $(f_3, g_4) = (0, 0)$.
\end{Assertion}

\proof
Applying now $\pi^4$ to the fundamental identity, 
taking into account that $F$ and $F'$ are normalized, we compute:
\[
\footnotesize
\aligned
0
&
\,\equiv\,
\pi^4
\bigg(
-\,\Im\,
\big(
w+0+g_4
\big)
+
\big(z+0+f_3\big)\,
\big(\overline{z}+0+\overline{f}_3\big)
+
\sum_{3\leqslant\mu\leqslant 4}\,
F_\mu'
\Big(
z+0+f_3,\,\,
\overline{z}+0+\overline{f}_3,\,\,
\Re\,\big(w+0+g_4\big)
\Big)
\bigg)
\\
&
\,\equiv\,
\pi^4
\bigg(
-\,z\overline{z}
-
\zero{F_3(z,\overline{z},u)}
-
\zero{F_4(z,\overline{z},u)}
+
\Re\,\big(i\,g_4\big)
+
z\overline{z}
+
z\overline{f}_3
+
\overline{z}f_3
+
\zero{f_3\overline{f}_3}
+
\zero{F_3'\big(z,\overline{z},u\big)}
+
\zero{F_4'\big(z,\overline{z},u\big)}
\bigg)
\\
&
\,\equiv\,
\pi^4
\bigg(
\Re\,\Big\{
i\,g_4
\big(
z,\,u+i[z\overline{z}+{\rm O}_{z,\overline{z},u}(3)]
\big)
+
2\,\overline{z}\,f_3
\big(
z,\,u+i[z\overline{z}+{\rm O}_{z,\overline{z},u}(3)]
\big)
\Big\}
\bigg)
\\
&
\,\equiv\,
\Re\,
\Big\{
i\,
g_4\big(z,\,u+i\,z\overline{z}\big)
+
2\,\overline{z}\,
f_3\big(z,\,u+i\,z\overline{z}\big)
\Big\}.
\endaligned
\]

Replacing $f_3 = f_{3,0}\, z^3 + f_{1,1}\, zw$ 
and $g_4 = 
g_{4,0}\, z^4 + g_{2,1}\, z^2w + g_{0,2}\, w^2$
with $\Re\, g_{0,2} = 0$ by assumption
(or even $g_{0,2} = 0$, but only null real part will suffice), 
this is:
\[
\footnotesize
\aligned
0
\,\equiv\,
\tfrac{i}{2}\,
g_{4,0}\,z^4
-
\tfrac{i}{2}\,
\overline{g}_{4,0}\,\overline{z}^4
&
+
\big(
f_{3,0}
-
\tfrac{1}{2}\,g_{2,1}
\big)\,
z^3\overline{z}
+
\big(
\overline{f}_{3,0}
-
\tfrac{1}{2}\,
\overline{g}_{2,1}
\big)\,
z\overline{z}^3
+
\big(
i\,f_{1,1}
-
i\,\overline{f}_{1,1}
-
\tfrac{i}{2}\,
g_{0,2}
+
\tfrac{i}{2}\,
\overline{g}_{0,2}
\big)\,
z^2\overline{z}^2
\\
&
+
\tfrac{i}{2}\,
g_{2,1}\,z^2u
-
\tfrac{i}{2}\,
\overline{g}_{2,1}\,
\overline{z}^2u
+
\big(
f_{1,1}
+
\overline{f}_{1,1}
-
g_{0,2}
-
\overline{g}_{0,2}
\big)\,
z\overline{z}u
+
\big(
\tfrac{i}{2}\,
g_{0,2}
-
\tfrac{i}{2}\,
\overline{g}_{0,2}
\big)\,
u^2,
\endaligned
\]
and starting from the end, since $g_{0,2}$ is purely imaginary, this
forces $0 = g_{0,2}$, then $f_{1,1} + \overline{f}_{1,1} = 0$, then $0
= g_{2,1}$, then $0 = f_{1,1}$, then $0 = f_{3,0}$, and lastly $0 =
g_{4,0}$, so as asserted $0 = f_3 = g_4$.
\endproof

Now, we discuss the induction vanishing process.  Assuming therefore
that $(f_{\mu-1},\, g_\mu) = (0, 0)$ for $\mu = 3, \dots, \nu-1$ and
some $\nu \geqslant 5$, we want to have $(f_{\nu-1},\, g_\nu) =
(0,0)$.

At first, it is not difficult to verify (left to the reader) that,
then:
\[
F_\mu'\big(z,\overline{z},u\big)
\,\equiv\,
F_\mu\big(z,\overline{z},u\big)
\eqno
{\scriptstyle{(\mu\,=\,3,\dots,\,\nu-1)}}.
\]
Using this, the fundamental identity then reads:
\[
\aligned
0
&
\,\equiv\,
\pi^\nu
\bigg(
-\,\Im\,\big(w+g_\nu\big)
+
\big(z+f_{\nu-1}\big)\,
\big(\overline{z}+\overline{f}_{\nu-1}\big)
+
\sum_{3\leqslant\mu\leqslant\nu}\,
F_\mu'
\Big(
z+f_{\nu-1},\,
\overline{z}+\overline{f}_{\nu-1},\,
u+\Re\,g_\nu
\Big)
\bigg)
\\
&
\,\equiv\,
\pi^\nu
\bigg(
-\,\zero{z\overline{z}}
-
\sum_{3\leqslant\mu\leqslant\nu-1}\,
\zerozero{F_\mu\big(z,\overline{z},u\big)}
-
F_\nu\big(z,\overline{z},u\big)
-
\Im\,g_\nu
+
\zero{z\overline{z}}
+
z\overline{f}_{\nu-1}
+
\overline{z}f_{\nu-1}
+
\zero{
f_{\nu-1}\overline{f}_{\nu-1}}
\\
&
\ \ \ \ \ \ \ \ \ \ \ \ \ \ \ \ \ \ \ \ \ \ \ \ \ \ \ \ \ \ \ \ \ \ \
\ \ \ \ \ \ \ \ \ \ \ \ \ \ \ \ \ \ \ \ \ \ \ \ \ \ \ \ \ \ \ \ \ \ \
\ \ \ \ \ \ \ \ \ \ \ \ \ \ \ \ \ \ \ \ \ \ \ \ \
+
\sum_{3\leqslant\mu\leqslant\nu-1}\,
\zerozero{F_\mu'\big(z,\overline{z},u\big)}
+
F_\nu'\big(z,\overline{z},u\big)
\bigg)
\\
&
\,\equiv\,
\pi^\nu
\bigg(
\Re\,
\Big\{
i\,g_\nu
\big(
z,\,u+i\,[z\overline{z}+{\rm O}_{z,\overline{z},u}(3)]
\big)
+
2\,\overline{z}\,f_{\nu-1}
\big(
z,\,u+i\,[z\overline{z}+{\rm O}_{z,\overline{z},u}(3)]
\big)
\Big\}
-
F_\nu\big(z,\overline{z},u\big)
+
F_\nu'\big(z,\overline{z},u\big)
\bigg)
\\
&
\,\equiv\,
\Re\,\Big\{
i\,g_\nu\big(z,\,u+iz\overline{z}\big)
+
2\,\overline{z}\,f_{\nu-1}\big(z,\,u+iz\overline{z}\big)
\Big\}
-
F_\nu\big(z,\overline{z},u\big)
+
F_\nu'\big(z,\overline{z},u\big).
\endaligned
\]

Now, we project further 
this equation by applying
to it $\Pi_S (\centersmallbullet)$.
Since $F$ and $F'$ are in normal form, we obtain, still
for any $\nu \geqslant 5$:
\[
0
\,\equiv\,
\Pi_S
\bigg(
\Re\,\Big\{
i\,g_\nu\big(z,\,u+iz\overline{z}\big)
+
2\,\overline{z}\,f_{\nu-1}\big(z,\,u+iz\overline{z}\big)
\Big\}
\bigg)
-
0
+
0.
\]
This is a linear system of equations in the coefficients 
$g_{j', l'}$
of $g_\nu$
and $f_{j', l'}$ of $f_{\nu-1}$. 
Instead of solving this linear system for any fixed
$\nu \geqslant 5$ (the cases $\nu = 3, 4$ have been done
above), we will solve in one stroke {\em all such systems
for any $\nu \geqslant 3$}, and this will simplify 
our job, especially by lightening a bit the combinatorics.

In any case, by taking the coefficients of all the monomials $z^j
\overline{z}^k u^l$ with $(j,k) \in S$ and $j + k + 2l = \nu$, we know
that there exist linear forms $L_{j,k,l}$ such that the above
system writes:
\[
0
\,=\,
L_{j,k,l}
\Big(
\big\{f_{j',l'}\big\}_{j'+2l'=\nu-1},\,\,
\big\{g_{j',l'}\big\}_{j'+2l'=\nu}
\Big),
\]
a system that we may abbreviate as:
\[
({\sf E}_\nu)
\colon\ \ \ \ \ \ \ \ \ \ 
0
\,=\,
L_{j,k,l}
\big(
f_{\smallbullet,\smallbullet},\,
g_{\smallbullet,\smallbullet}
\big)
\eqno
{\scriptstyle{((j,k)\,\in\,S,\,\,j+k+2l\,=\,\nu)}}.
\]
From now on, $\nu \geqslant 3$, so we incorporate $\nu = 3, 4$
in the discussion. 

On the other hand, by considering
the complete $f = f_2 + f_3 + \cdots$ and the
complete $g = g_3 + g_4 + \cdots$, we can introduce the
analog `complete' linear system:
\[
0
\,\equiv\,
\Pi_S
\bigg(
\Re\,
\Big\{
i\,g\big(z,\,u+iz\overline{z}\big)
+
2\,\overline{z}\,f\big(z,\,u+iz\overline{z}\big)
\Big\}
\bigg),
\]
which, similarly, after extracting the coefficients of
all monomials $z^j \overline{z}^k u^l$ with 
$(j,k) \in S$ and any $l \in \N$, 
can be abbreviated as:
\[
({\sf E})
\colon\ \ \ \ \ \ \ \ \ \ 
0
\,=\,
L_{j,k,l}
\big(
f_{\smallbullet,\smallbullet},\,
g_{\smallbullet,\smallbullet}
\big)
\eqno
{\scriptstyle{((j,k)\,\in\,S,\,\,l\,\in\,\N)}}.
\]
The key and elementary observation 
is that, because $u + iz \overline{z}$ is $2$-homogeneous,
the full system $({\sf E})$
{\em splits in the linear subsystems 
$({\sf E}_\nu)$
having separate unknowns $\big(f_{\nu-1},\, g_{\nu} \big)$}:
\[
({\sf E})
\,=\,
({\sf E}_3)
\cup
({\sf E}_4)
\cup\cdots\cup
({\sf E}_\nu)
\cup
\cdots.
\]
Therefore:
\[
\Big(
({\sf E})
\,\,\,\Longrightarrow\,\,\,
(f,g)=(0,0)
\Big)
\ \ \ \ \ \ \ \ \ \ \ \ \ \ \ \ \ \ \ \
\Longleftrightarrow
\ \ \ \ \ \ \ \ \ \ \ \ \ \ \ \ \ \ \ \
\Big(
({\sf E}_\nu)
\,\,\,\Longrightarrow\,\,\,
\big(f_{\nu-1},g_\nu\big)=(0,0)
\ \ \ \ \
\text{for all}\,\,
\nu\geqslant 3
\Big).
\]

Thus, we are left with establishing the following main technical
statement, which will close the proof of
Theorem~{\ref{Thm-uniqueness-normal-form}}.
\endproof

\begin{Theorem}
Let $f(z,w)$ and $g(z,w)$ be holomorphic 
of weights $\geqslant 2$ and $\geqslant 3$, namely
$f = f_2 + f_3 + \cdots$ and $g = g_3 + g_4 + \cdots$, 
and with:
\[
0
\,=\,
f_w(0),
\ \ \ \ \ \ \ \ \ \ \ \ \ \ \ \ \ \ \ \ \ \ \ \ \ \
0
\,=\,
\Re\,g_{ww}(0).
\]
If for all $(j,k) \in S$ and all $l \in \N$:
\[
0
\,=\,
\big[z^j\overline{z}^ku^l\big]\,
\bigg(
\Re\,\Big\{
i\,g\big(z,\,u+iz\overline{z}\big)
+
2\,\overline{z}\,f\big(z,\,u+iz\overline{z}\big)
\Big\}
\bigg),
\]
then $(f,g) \equiv (0,0)$.
\end{Theorem}

\proof
The key simplification is to gather all powers $u^l$
in the linear system so as to deal with finitely many
functions of the CR-transversal variable $u$.

Indeed, given a holomorphic function $e = e(w)$, we may expand:
\[
e\big(
u+i\,z\overline{z}
\big)
\,=\,
e(u)
+
e_w(u)\,
i\,z\overline{z}
+
e_{ww}(u)\,
\tfrac{1}{2!}\,
\big(i\,z\overline{z}\big)^2
+
e_{www}(u)\,
\tfrac{1}{3!}\,
\big(i\,z\overline{z}\big)^3
+\cdots,
\]
and we will write $e'(u)$, $e''(u)$, $e'''(u)$, {\em etc.},
instead of $e_w(u)$, $e_{ww}(u)$, $e_{www}(u)$, {\em etc.}
Thus:
\[
\aligned
f\big(z,\,u+i\,z\overline{z}\big)
&
\,=\,
\sum_{k\geqslant0}\,
z^k\,f_k\big(u+i\,z\overline{z}\big)
\\
&
\,=\,
\sum_{k\geqslant 0}\,
z^k\,
\Big[
f_k(u)
+
f_k'(u)\,
i\,z\overline{z}
+
f_k''(u)\,
\tfrac{1}{2!}\,
\big(i\,z\overline{z}\big)^2
+
f_k'''(u)\,
\tfrac{1}{3!}\,
\big(i\,z\overline{z}\big)^3
+\cdots
\Big],
\endaligned
\]
and similarly:
\[
\aligned
g\big(z,\,u+i\,z\overline{z}\big)
&
\,=\,
\sum_{k\geqslant0}\,
z^k\,g_k\big(u+i\,z\overline{z}\big)
\\
&
\,=\,
\sum_{k\geqslant 0}\,
z^k\,
\Big[
g_k(u)
+
g_k'(u)\,
i\,z\overline{z}
+
g_k''(u)\,
\tfrac{1}{2!}\,
\big(i\,z\overline{z}\big)^2
+
g_k'''(u)\,
\tfrac{1}{3!}\,
\big(i\,z\overline{z}\big)^3
+\cdots
\Big],
\endaligned
\]

Hence our zero equation is:
\[
\aligned
0
&
\,\equiv\,
2\,\Re\,
\Big\{
2\,\overline{z}\,
f\big(z,\,u+i\,z\overline{z}\big)
+
i\,g\big(z,\,u+i\,z\overline{z}\big)
\Big\}
\\
&
\,\equiv\,
2\,\overline{z}\,f
+
2\,z\,\overline{f}
+
i\,g
-
i\,\overline{g}
\\
&
\,\equiv\,
\sum_{k\geqslant 0}\,
\Big(
2\,f_k\,z^k\overline{z}
+
2\,i\,f_k'\,z^{k+1}\overline{z}^2
-
f_k''\,z^{k+2}\overline{z}^3
-
\tfrac{i}{3}\,f_k'''\,z^{k+3}\overline{z}^4
+\cdots
\Big)
\\
&\ \ \ 
+
\sum_{k\geqslant 0}\,
\Big(
2\,\overline{f}_k\,z\overline{z}^k
-
2\,i\,\overline{f}_k'\,z^2\overline{z}^{k+1}
-
\overline{f}_k''\,z^3\overline{z}^{k+2}
+
\tfrac{i}{3}\,\overline{f}_k'''\,z^4\overline{z}^{k+3}
+\cdots
\Big)
\\
&\ \ \ 
+
\sum_{k\geqslant 0}\,
\Big(
i\,g_k\,z^k
-
g_k'\,z^{k+1}\overline{z}
-
\tfrac{i}{2}\,g_k''\,z^{k+2}\overline{z}^2
+
\tfrac{1}{6}\,g_k'''\,z^{k+3}\overline{z}^3
+\cdots
\Big)
\\
&\ \ \ 
+
\sum_{k\geqslant 0}\,
\Big(
-i\,\overline{g}_k\,\overline{z}^k
-
\overline{g}_k'\,z\overline{z}^{k+1}
+
\tfrac{i}{2}\,\overline{g}_k''\,z^2\overline{z}^{k+2}
+
\tfrac{1}{6}\,\overline{g}_k'''\,z^3\overline{z}^{k+3}
+\cdots
\Big),
\endaligned
\]
where the common argument of all $f_k$, $f_k'$, $f_k''$,
$g_k$, $g_k'$, $g_k''$, $g_k'''$ is $u \in \R$.

We are thus capturing the coefficients 
$[z^j \overline{z}^k] (\centersmallbullet)$ of this
identity, not anymore all $[z^j \overline{z}^k u^l] 
(\centersmallbullet)$. This means that we are extracting
identities satisfied by functions of $u$.

Let us therefore list the coefficients of $z^j \overline{z}^k$,
indicating plainly $(j,k)$. Note that we can restrict the
considerations to only $j \geqslant k$, since the above
zero equation is real.

\begin{center}
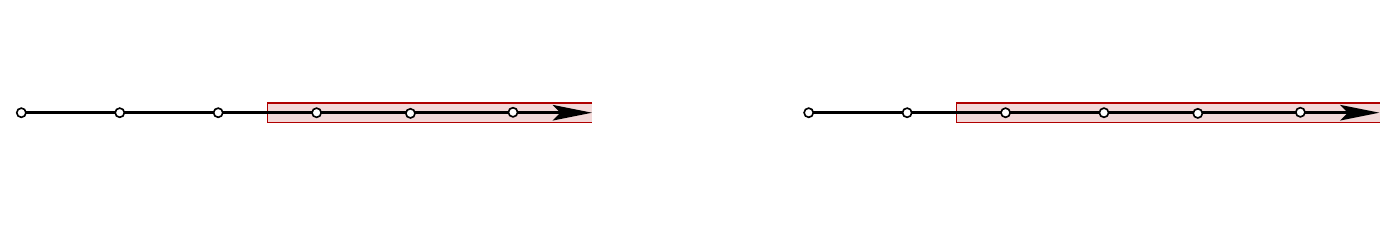

\nopagebreak
\begin{minipage}[t]{12.5cm}
\footnotesize
{\sc Figure~9:}
Two infinite red-shaded families of coefficients 
$f_3(u) \equiv f_4(u) \equiv \cdots \equiv 0$ and 
$g_2(u) \equiv g_3(u) \equiv \cdots \equiv 0$
easily shown to vanish identically.
\end{minipage}
\end{center}

Firstly, we extract the coefficients of $z^k$ with
$k \geqslant 2$ and of $z^k \overline{z}$ for $k \geqslant 3$:
\leqnomode\usetagform{default}
\begin{align}
0
&
\,=\,
i\,g_k,
\tag{$k^{\geqslant 2},0$}
\\
0
&
\,=\,
2\,f_k
-
g_{k-1}'.
\tag{$k^{\geqslant 3},1$}
\end{align}
So $g_k (u) \equiv 0$ for all $k \geqslant 2$ and
$f_k(u) \equiv 0$ for all $k \geqslant 3$, and therefore:
\[
f
\,=\,
f_0(w)
+
z\,f_1(w)
+
z^2\,f_2(w),
\ \ \ \ \ \ \ \ \ \ \ \ \ \ \ \ \ \ \ \ \ \ \ \ \ \
g
\,=\,
g_0(w)
+
z\,g_1(w).
\]

Next, we extract the remaining coefficients of $z^j \overline{z}^k$,
and we get $7$ equations:
\leqnomode\usetagform{default}
\begin{align}
0
&
\,=\,
i\,g_0
-
i\,\overline{g}_0,
\tag{$0,0$}
\\
0
&
\,=\,
2\,\overline{f}_0
+
i\,g_1,
\tag{$1,0$}
\\
0
&
\,=\,
2\,f_1
+
2\,\overline{f}_1
-
g_0'
-
\overline{g}_0',
\tag{$1,1$}
\\
0
&
\,=\,
2\,f_2
-
2i\,\overline{f}_0'
-
g_1',
\tag{$2,1$}
\\
0
&
\,=\,
2i\,f_1'
-
2i\,\overline{f}_1'
-
\tfrac{i}{2}\,g_0''
+
\tfrac{i}{2}\,\overline{g}_0'',
\tag{$2,2$}
\\
0
&
\,=\,
2i\,f_2'
-
\overline{f}_0''
-
\tfrac{i}{2}\,g_1'',
\tag{$3,2$}
\\
0
&
\,=\,
-\,f_1''
-
\overline{f}_1''
+
\tfrac{1}{6}\,g_0'''
-
\tfrac{1}{6}\,\overline{g}_0'''.
\tag{$3,3$}
\end{align}
Now, since:
\[
\aligned
0
&
\,=\,
f(0)
\ \ \ \ \
&
\Longleftrightarrow
\ \ \ \ \
f_0(0)
\,=\,
0,
&
\ \ \ \ \ \ \ \ \ \ \ \ \ \ \ \ \ \ \ \ \ \ \ \ \ \
0
&
\,=\,\ 
g(0)
\ \ \ \ \
&
\Longleftrightarrow
\ \ \ \ \
g_0(0)
&
\,=\,
0,
\\
0
&
\,=\,
f_z(0)
\ \ \ \ \
&
\Longleftrightarrow
\ \ \ \ \
f_1(0)
\,=\,
0,
&
\ \ \ \ \ \ \ \ \ \ \ \ \ \ \ \ \ \ \ \ \ \ \ \ \ \
0
&
\,=\,
g_z(0)
\ \ \ \ \
&
\Longleftrightarrow
\ \ \ \ \
g_1(0)
&
\,=\,
0,
\\
0
&
\,=\,
f_w(0)
\ \ \ \ \
&
\Longleftrightarrow
\ \ \ \ \
f_0'(0)
\,=\,
0,
&
\ \ \ \ \ \ \ \ \ \ \ \ \ \ \ \ \ \ \ \ \ \ \ \ \ \
0
&
\,=\,
g_w(0)
\ \ \ \ \
&
\Longleftrightarrow
\ \ \ \ \
g_0'(0)
&
\,=\,
0,
\endaligned
\]
and since:
\[
0
\,=\,
\Re\,g_{ww}(0)
\ \ \ \ \
\Longleftrightarrow
\ \ \ \ \
\Re\,g_0''(0)
\,=\,
0,
\]
the assumptions of the theorem are equivalent
to the ones formulated in the next statement, 
which will finish everything. 
\endproof

\begin{Assertion}
If five functions $f_0$, $f_1$, $f_2$, $g_0$, $g_1$
of the real variable $u \in \R$ with:
\[
\aligned
0
&
\,=\,
f_0(0)
\,=\,
f_0'(0),
&
\ \ \ \ \ \ \ \ \ \ \ \ \ \ \ \ \ \ \ \ \ \ \ \ \ \
0
&
\,=\,
g_0(0)
\,=\,
g_0'(0)
\,=\,
\Re\,g_0''(0),
\\
0
&
\,=\,
f_1(0),
&
\ \ \ \ \ \ \ \ \ \ \ \ \ \ \ \ \ \ \ \ \ \ \ \ \ \
0
&
\,=\,
g_1(0),
\endaligned
\]
satisfy the above $7$ linear ordinary differential equations,
then they all vanish identically:
\[
0
\,\equiv\,
f_0(u)
\,\equiv\,
f_1(u)
\,\equiv\,
f_2(u),
\ \ \ \ \ \ \ \ \ \ \ \ \ \ \ \ \ \ \ \ \ \ \ \ \ \
0
\,\equiv\,
g_0(u)
\,\equiv\,
g_1(u).
\]
\end{Assertion}

\proof
From $(0,0)$, solve $\overline{g}_0 := g_0$. From
$(1,0)$, solve $g_1 := 2i\, \overline{f}_0$. Then
the five remaining equations become:
\leqnomode\usetagform{default}
\begin{align}
0
&
\,=\,
2\,f_1(u)
+
2\,\overline{f}_1(u)
-
2\,g_0'(u),
\tag{$1,1$}
\\
0
&
\,=\,
2\,f_2(u)
-
4i\,\overline{f}_0'(u),
\tag{$2,1$}
\\
0
&
\,=\,
2i\,f_1'(u)
-
2i\,\overline{f}_1'(u),
\tag{$2,2$}
\\
0
&
\,=\,
2i\,f_2'(u),
\tag{$3,2$}
\\
0
&
\,=\,
-\,f_1''(u)
-
\overline{f}_1''(u).
\tag{$3,3$}
\end{align}

From $(3,2)$, we see $f_2 = \alpha \in \C$ is constant.
From $(2,1)$ at $u = 0$, since $f_0'(0) = 0$ by assumption,
we get $\alpha = 0$. So $f_2(u) \equiv 0$ in $(2,1)$ gives
$f_0(u) \equiv 0$ too. Thus $g_1(u) \equiv 0$ as well.

From $(3,3)$ and $\frac{d}{du} (2,2)$, it comes $f_1''(u) \equiv 0$,
and since $f_1(0) = 0$ by assumption, $f_1(u) = c\, u$
with $c \in \R$ by $(2,2)$. From $(1,1)$, it comes $g_0'(u) = 2\,
c\, u$,
and since $\Re\, g_0''(0) = 0$, we get $c = 0$.
Thus $f_1(u) \equiv 0$. 

From $g_0'(u) \equiv 0$ and $g_0(0) = 0$ by assumption,
we get $g_0(u) \equiv 0$. This concludes.
\endproof



\vfill\end{document}

%% file: macros.tex
\newtheorem{Theorem}[equation]{Theorem}

\newtheorem{Proposition}[equation]{Proposition}

\newtheorem{Lemma}[equation]{Lemma}

\newtheorem{Assertion}[equation]{Assertion}

\newtheorem{Observation}[equation]{Observation}

\theoremstyle{definition}

\newtheorem{Definition}[equation]{Definition}

\newtheorem{Question}[equation]{Question}

\newcommand{\C}{\mathbb{C}}

\newcommand{\N}{\mathbb{N}}

\renewcommand{\P}{\mathbb{P}}

\newcommand{\R}{\mathbb{R}}

\newcommand{\CC}{\text{\sc c}}

\newcommand{\TT}{\text{\sc t}}

\newcommand{\Iaux}{{\text{\usefont{T1}{qcs}{m}{sl}I}}}

\definecolor{blue}{cmyk}{1.,1.,0.,0.63}
\definecolor{red}{cmyk}{0.,1.,1.,0.63}
\definecolor{green}{cmyk}{1.,0.,1.,0.63}
\definecolor{black}{cmyk}{1.,1.,1.,1.}

\newcommand{\blue}{\textcolor{blue}}
\newcommand{\green}{\textcolor{green}}
\newcommand{\red}{\textcolor{red}}

\makeatletter
\renewcommand{\@fnsymbol}[1]
{\ensuremath{\ifcase#1\or $*$\or $**$\or $***$\or $****$\or $*****$
\else\@ctrerr\fi}}
\makeatother

\newcommand{\HEAD}[2]{%
\pagestyle{fancy}
\fancyhead[RO]{\tiny\sf\thepage}
\fancyhead[CO]{{\tiny\sf #1}}
\fancyhead[LE]{\tiny\sf\thepage}
\fancyhead[CE]{{\tiny\sf #2}}
\fancyfoot{}}

\numberwithin{equation}{section}

\newcommand{\Section}[1]{
\renewcommand{\thesection}{\bf\arabic{section}}
\section{#1}
\renewcommand{\thesection}{\arabic{section}}}

\newcommand{\style}[1]{\text{\footnotesize{\sf #1}}}

\renewcommand{\exp}{\style{exp}}

\newcommand{\Id}{\style{Id}}

\renewcommand{\Im}{\style{Im}}

\renewcommand{\lim}{\style{lim}}

\renewcommand{\Re}{\style{Re}}

\newcommand{\centersmallbullet}{{}_{{}^{{}^{
\scriptscriptstyle{\bullet\!}}}}}

\newcommand{\Hall}{\Hall}

\newcommand{\smallbullet}{{\scriptscriptstyle{\bullet}}}

\newcommand{\smallsum}[1]{
\underset{#1}{\raisebox{1pt}{$\sum$\,}}
}

\newcommand{\vf}{\vfill\end{document}}

\newcommand{\zero}[1]{\underline{#1}_{\red{\circ}}}

\newcommand{\zerozero}[1]{\underline{#1}_{\red{\circ\circ}}}


%% file: various-chains.pdf_t
\begin{picture}(0,0)%
\includegraphics{various-chains.pdf}%
\end{picture}%
\setlength{\unitlength}{4144sp}%
\begingroup\makeatletter\ifx\SetFigFont\undefined%
\gdef\SetFigFont#1#2#3#4#5{%
  \reset@font\fontsize{#1}{#2pt}%
  \fontfamily{#3}\fontseries{#4}\fontshape{#5}%
  \selectfont}%
\fi\endgroup%
\begin{picture}(7188,1566)(427,-2311)
\put(615,-892){\makebox(0,0)[lb]{\smash{{\SetFigFont{9}{10.8}{\familydefault}{\mddefault}{\updefault}{\color[rgb]{0,0,0}\blue{$\C^2$}}%
}}}}
\put(785,-1427){\makebox(0,0)[lb]{\smash{{\SetFigFont{8}{9.6}{\familydefault}{\mddefault}{\updefault}{\color[rgb]{0,0,0}\green{$T_p^cM$}}%
}}}}
\put(1485,-1577){\makebox(0,0)[lb]{\smash{{\SetFigFont{8}{9.6}{\familydefault}{\mddefault}{\updefault}{\color[rgb]{0,0,0}$p$}%
}}}}
\put(499,-2166){\makebox(0,0)[lb]{\smash{{\SetFigFont{9}{10.8}{\familydefault}{\mddefault}{\updefault}{\color[rgb]{0,0,0}\blue{$M$}}%
}}}}
\put(3024,-921){\makebox(0,0)[lb]{\smash{{\SetFigFont{9}{10.8}{\familydefault}{\mddefault}{\updefault}{\color[rgb]{0,0,0}\blue{$\C^2$}}%
}}}}
\put(3194,-1456){\makebox(0,0)[lb]{\smash{{\SetFigFont{8}{9.6}{\familydefault}{\mddefault}{\updefault}{\color[rgb]{0,0,0}\green{$T_p^cM$}}%
}}}}
\put(3894,-1606){\makebox(0,0)[lb]{\smash{{\SetFigFont{8}{9.6}{\familydefault}{\mddefault}{\updefault}{\color[rgb]{0,0,0}$p$}%
}}}}
\put(2908,-2195){\makebox(0,0)[lb]{\smash{{\SetFigFont{9}{10.8}{\familydefault}{\mddefault}{\updefault}{\color[rgb]{0,0,0}\blue{$M$}}%
}}}}
\put(4096,-986){\makebox(0,0)[lb]{\smash{{\SetFigFont{9}{10.8}{\familydefault}{\mddefault}{\updefault}{\color[rgb]{0,0,0}\green{$\ell_p$}}%
}}}}
\put(5428,-2195){\makebox(0,0)[lb]{\smash{{\SetFigFont{9}{10.8}{\familydefault}{\mddefault}{\updefault}{\color[rgb]{0,0,0}\blue{$M'$}}%
}}}}
\put(5544,-921){\makebox(0,0)[lb]{\smash{{\SetFigFont{9}{10.8}{\familydefault}{\mddefault}{\updefault}{\color[rgb]{0,0,0}\blue{${\C'}^2$}}%
}}}}
\put(1696,-970){\makebox(0,0)[lb]{\smash{{\SetFigFont{9}{10.8}{\familydefault}{\mddefault}{\updefault}{\color[rgb]{0,0,0}\green{$\ell_p$}}%
}}}}
\put(2012,-1074){\makebox(0,0)[lb]{\smash{{\SetFigFont{9}{10.8}{\familydefault}{\mddefault}{\updefault}{\color[rgb]{0,0,0}\green{$\ell_p$}}%
}}}}
\put(1373,-966){\makebox(0,0)[lb]{\smash{{\SetFigFont{9}{10.8}{\familydefault}{\mddefault}{\updefault}{\color[rgb]{0,0,0}\green{$\ell_p$}}%
}}}}
\put(5157,-1353){\makebox(0,0)[lb]{\smash{{\SetFigFont{8}{9.6}{\familydefault}{\mddefault}{\updefault}{\color[rgb]{0,0,0}$h$}%
}}}}
\put(6390,-1641){\makebox(0,0)[lb]{\smash{{\SetFigFont{8}{9.6}{\familydefault}{\mddefault}{\updefault}{\color[rgb]{0,0,0}$h(p)$}%
}}}}
\put(5731,-1443){\makebox(0,0)[lb]{\smash{{\SetFigFont{8}{9.6}{\familydefault}{\mddefault}{\updefault}{\color[rgb]{0,0,0}\green{$T_{h(p)}^cM'$}}%
}}}}
\put(6943,-1100){\makebox(0,0)[lb]{\smash{{\SetFigFont{9}{10.8}{\familydefault}{\mddefault}{\updefault}{\color[rgb]{0,0,0}\green{$h_\ast(\ell_p)$}}%
}}}}
\put(6025,-1999){\makebox(0,0)[lb]{\smash{{\SetFigFont{9}{10.8}{\familydefault}{\mddefault}{\updefault}{\color[rgb]{0,0,0}$\CC_{h(p),h_\ast(\ell_p)}^{h(M)}$}%
}}}}
\put(3753,-2014){\makebox(0,0)[lb]{\smash{{\SetFigFont{9}{10.8}{\familydefault}{\mddefault}{\updefault}{\color[rgb]{0,0,0}$\CC_{p,\ell_p}^M$}%
}}}}
\end{picture}%

%% file: kill-Moser.pdf_t
\begin{picture}(0,0)%
\includegraphics{kill-Moser.pdf}%
\end{picture}%
\setlength{\unitlength}{4144sp}%
\begingroup\makeatletter\ifx\SetFigFont\undefined%
\gdef\SetFigFont#1#2#3#4#5{%
  \reset@font\fontsize{#1}{#2pt}%
  \fontfamily{#3}\fontseries{#4}\fontshape{#5}%
  \selectfont}%
\fi\endgroup%
\begin{picture}(6136,1448)(799,-1423)
\put(4023,-353){\makebox(0,0)[lb]{\smash{{\SetFigFont{8}{9.6}{\familydefault}{\mddefault}{\updefault}{\color[rgb]{0,0,0}\blue{$k$}}%
}}}}
\put(4668,-348){\makebox(0,0)[lb]{\smash{{\SetFigFont{9}{10.8}{\familydefault}{\mddefault}{\updefault}{\color[rgb]{0,0,0}Step~4}%
}}}}
\put(3175,-547){\makebox(0,0)[lb]{\smash{{\SetFigFont{9}{10.8}{\familydefault}{\mddefault}{\updefault}{\color[rgb]{0,0,0}Step~3}%
}}}}
\put(2764,-561){\makebox(0,0)[lb]{\smash{{\SetFigFont{8}{9.6}{\familydefault}{\mddefault}{\updefault}{\color[rgb]{0,0,0}\blue{$k$}}%
}}}}
\put(2071,-781){\makebox(0,0)[lb]{\smash{{\SetFigFont{9}{10.8}{\familydefault}{\mddefault}{\updefault}{\color[rgb]{0,0,0}Step~2}%
}}}}
\put(5828,-635){\makebox(0,0)[lb]{\smash{{\SetFigFont{8}{9.6}{\familydefault}{\mddefault}{\updefault}{\color[rgb]{0,0,0}\blue{$F_{2,3}$}}%
}}}}
\put(6045,-866){\makebox(0,0)[lb]{\smash{{\SetFigFont{8}{9.6}{\familydefault}{\mddefault}{\updefault}{\color[rgb]{0,0,0}\blue{$F_{3,2}$}}%
}}}}
\put(5489,-122){\makebox(0,0)[lb]{\smash{{\SetFigFont{8}{9.6}{\familydefault}{\mddefault}{\updefault}{\color[rgb]{0,0,0}\blue{$k$}}%
}}}}
\put(6367,-124){\makebox(0,0)[lb]{\smash{{\SetFigFont{9}{10.8}{\familydefault}{\mddefault}{\updefault}{\color[rgb]{0,0,0}Step~5}%
}}}}
\put(6678,-1311){\makebox(0,0)[lb]{\smash{{\SetFigFont{8}{9.6}{\familydefault}{\mddefault}{\updefault}{\color[rgb]{0,0,0}\blue{$j$}}%
}}}}
\put(4968,-1284){\makebox(0,0)[lb]{\smash{{\SetFigFont{8}{9.6}{\familydefault}{\mddefault}{\updefault}{\color[rgb]{0,0,0}\blue{$j$}}%
}}}}
\put(3495,-1278){\makebox(0,0)[lb]{\smash{{\SetFigFont{8}{9.6}{\familydefault}{\mddefault}{\updefault}{\color[rgb]{0,0,0}\blue{$j$}}%
}}}}
\put(1735,-788){\makebox(0,0)[lb]{\smash{{\SetFigFont{8}{9.6}{\familydefault}{\mddefault}{\updefault}{\color[rgb]{0,0,0}\blue{$k$}}%
}}}}
\put(2236,-1286){\makebox(0,0)[lb]{\smash{{\SetFigFont{8}{9.6}{\familydefault}{\mddefault}{\updefault}{\color[rgb]{0,0,0}\blue{$j$}}%
}}}}
\put(1079,-981){\makebox(0,0)[lb]{\smash{{\SetFigFont{9}{10.8}{\familydefault}{\mddefault}{\updefault}{\color[rgb]{0,0,0}Step~1}%
}}}}
\put(882,-969){\makebox(0,0)[lb]{\smash{{\SetFigFont{7}{8.4}{\familydefault}{\mddefault}{\updefault}{\color[rgb]{0,0,0}\blue{$k$}}%
}}}}
\put(1205,-1314){\makebox(0,0)[lb]{\smash{{\SetFigFont{8}{9.6}{\familydefault}{\mddefault}{\updefault}{\color[rgb]{0,0,0}\blue{$j$}}%
}}}}
\end{picture}%

%% file: re-translation-normalization.pdf_t
\begin{picture}(0,0)%
\includegraphics{re-translation-normalization.pdf}%
\end{picture}%
\setlength{\unitlength}{4144sp}%
\begingroup\makeatletter\ifx\SetFigFont\undefined%
\gdef\SetFigFont#1#2#3#4#5{%
  \reset@font\fontsize{#1}{#2pt}%
  \fontfamily{#3}\fontseries{#4}\fontshape{#5}%
  \selectfont}%
\fi\endgroup%
\begin{picture}(6341,1405)(661,-1943)
\put(4994,-688){\makebox(0,0)[lb]{\smash{{\SetFigFont{8}{9.6}{\familydefault}{\mddefault}{\updefault}{\color[rgb]{0,0,0}Normalization}%
}}}}
\put(5311,-846){\makebox(0,0)[lb]{\smash{{\SetFigFont{9}{10.8}{\familydefault}{\mddefault}{\updefault}{\color[rgb]{0,0,0}$\Phi_p$}%
}}}}
\put(2004,-1271){\makebox(0,0)[lb]{\smash{{\SetFigFont{9}{10.8}{\familydefault}{\mddefault}{\updefault}{\color[rgb]{0,0,.69}\blue{$p$}}%
}}}}
\put(826,-1411){\makebox(0,0)[lb]{\smash{{\SetFigFont{9}{10.8}{\familydefault}{\mddefault}{\updefault}{\color[rgb]{0,0,.69}\blue{$M$}}%
}}}}
\put(2335,-922){\makebox(0,0)[lb]{\smash{{\SetFigFont{9}{10.8}{\familydefault}{\mddefault}{\updefault}{\color[rgb]{0,0,0}Translation}%
}}}}
\put(2572,-1045){\makebox(0,0)[lb]{\smash{{\SetFigFont{8}{9.6}{\familydefault}{\mddefault}{\updefault}{\color[rgb]{0,0,0}$\tau_p$}%
}}}}
\put(4258,-1307){\makebox(0,0)[lb]{\smash{{\SetFigFont{9}{10.8}{\familydefault}{\mddefault}{\updefault}{\color[rgb]{0,0,.69}\blue{$0$}}%
}}}}
\put(3069,-1414){\makebox(0,0)[lb]{\smash{{\SetFigFont{9}{10.8}{\familydefault}{\mddefault}{\updefault}{\color[rgb]{0,0,.69}\blue{$M^p$}}%
}}}}
\put(4617,-1108){\makebox(0,0)[lb]{\smash{{\SetFigFont{9}{10.8}{\familydefault}{\mddefault}{\updefault}{\color[rgb]{0,0,.69}\blue{$z,\overline{z},u$}}%
}}}}
\put(6508,-1307){\makebox(0,0)[lb]{\smash{{\SetFigFont{9}{10.8}{\familydefault}{\mddefault}{\updefault}{\color[rgb]{0,0,.69}\blue{$0$}}%
}}}}
\put(4322,-769){\makebox(0,0)[lb]{\smash{{\SetFigFont{9}{10.8}{\familydefault}{\mddefault}{\updefault}{\color[rgb]{0,0,.69}\blue{$v$}}%
}}}}
\put(6572,-769){\makebox(0,0)[lb]{\smash{{\SetFigFont{9}{10.8}{\familydefault}{\mddefault}{\updefault}{\color[rgb]{0,0,.69}\blue{$v'$}}%
}}}}
\put(6871,-1112){\makebox(0,0)[lb]{\smash{{\SetFigFont{9}{10.8}{\familydefault}{\mddefault}{\updefault}{\color[rgb]{0,0,.69}\blue{$z',\overline{z}',u'$}}%
}}}}
\put(5416,-1879){\makebox(0,0)[lb]{\smash{{\SetFigFont{9}{10.8}{\familydefault}{\mddefault}{\updefault}{\color[rgb]{0,0,0}$N^p$}%
}}}}
\put(3478,-1799){\makebox(0,0)[lb]{\smash{{\SetFigFont{8}{9.6}{\familydefault}{\mddefault}{\updefault}{\color[rgb]{0,0,0}$\Phi_p\circ\tau_p$}%
}}}}
\end{picture}%

%% file: ambiguity-group.pdf_t
\begin{picture}(0,0)%
\includegraphics{ambiguity-group.pdf}%
\end{picture}%
\setlength{\unitlength}{4144sp}%
\begingroup\makeatletter\ifx\SetFigFont\undefined%
\gdef\SetFigFont#1#2#3#4#5{%
  \reset@font\fontsize{#1}{#2pt}%
  \fontfamily{#3}\fontseries{#4}\fontshape{#5}%
  \selectfont}%
\fi\endgroup%
\begin{picture}(5820,2052)(466,-1636)
\put(731,-586){\makebox(0,0)[lb]{\smash{{\SetFigFont{9}{10.8}{\familydefault}{\mddefault}{\updefault}{\color[rgb]{0,0,0}$M_p$}%
}}}}
\put(3016,-1091){\makebox(0,0)[lb]{\smash{{\SetFigFont{9}{10.8}{\familydefault}{\mddefault}{\updefault}{\color[rgb]{0,0,0}normalization~2}%
}}}}
\put(5413,-490){\makebox(0,0)[lb]{\smash{{\SetFigFont{9}{10.8}{\familydefault}{\mddefault}{\updefault}{\color[rgb]{0,0,0}\red{\bf ambiguity}}%
}}}}
\put(3023,-477){\makebox(0,0)[lb]{\smash{{\SetFigFont{9}{10.8}{\familydefault}{\mddefault}{\updefault}{\color[rgb]{0,0,0}normalization~1}%
}}}}
\put(1306,-719){\makebox(0,0)[lb]{\smash{{\SetFigFont{9}{10.8}{\familydefault}{\mddefault}{\updefault}{\color[rgb]{0,0,0}$0$}%
}}}}
\put(5991,269){\makebox(0,0)[lb]{\smash{{\SetFigFont{9}{10.8}{\familydefault}{\mddefault}{\updefault}{\color[rgb]{0,0,0}$v\!=\!z\overline{z}\!+\!{\rm O}(6)$}%
}}}}
\put(6101,-1441){\makebox(0,0)[lb]{\smash{{\SetFigFont{9}{10.8}{\familydefault}{\mddefault}{\updefault}{\color[rgb]{0,0,0}$v'\!=\!z'\overline{z}'\!+\!{\rm O}(6)$}%
}}}}
\end{picture}%

%% file: variation-vector-v.pdf_t
\begin{picture}(0,0)%
\includegraphics{variation-vector-v.pdf}%
\end{picture}%
\setlength{\unitlength}{4144sp}%
\begingroup\makeatletter\ifx\SetFigFont\undefined%
\gdef\SetFigFont#1#2#3#4#5{%
  \reset@font\fontsize{#1}{#2pt}%
  \fontfamily{#3}\fontseries{#4}\fontshape{#5}%
  \selectfont}%
\fi\endgroup%
\begin{picture}(2294,2345)(879,-2773)
\put(2046,-774){\makebox(0,0)[lb]{\smash{{\SetFigFont{9}{10.8}{\familydefault}{\mddefault}{\updefault}{\color[rgb]{0,0,0}\blue{$u$}}%
}}}}
\put(1381,-1129){\makebox(0,0)[lb]{\smash{{\SetFigFont{9}{10.8}{\familydefault}{\mddefault}{\updefault}{\color[rgb]{0,0,0}\green{$\vec{\bf v}$}}%
}}}}
\put(1724,-872){\makebox(0,0)[lb]{\smash{{\SetFigFont{9}{10.8}{\familydefault}{\mddefault}{\updefault}{\color[rgb]{0,0,0}\green{$\vec{\bf v}$}}%
}}}}
\put(2297,-918){\makebox(0,0)[lb]{\smash{{\SetFigFont{9}{10.8}{\familydefault}{\mddefault}{\updefault}{\color[rgb]{0,0,0}\green{$\vec{\bf v}$}}%
}}}}
\put(1896,-1964){\makebox(0,0)[lb]{\smash{{\SetFigFont{9}{10.8}{\familydefault}{\mddefault}{\updefault}{\color[rgb]{0,0,0}\blue{$y$}}%
}}}}
\put(3045,-1648){\makebox(0,0)[lb]{\smash{{\SetFigFont{9}{10.8}{\familydefault}{\mddefault}{\updefault}{\color[rgb]{0,0,0}\blue{$x$}}%
}}}}
\put(2463,-1821){\makebox(0,0)[lb]{\smash{{\SetFigFont{9}{10.8}{\familydefault}{\mddefault}{\updefault}{\color[rgb]{0,0,0}$T_0^cM^p$}%
}}}}
\put(2034,-1755){\makebox(0,0)[lb]{\smash{{\SetFigFont{9}{10.8}{\familydefault}{\mddefault}{\updefault}{\color[rgb]{0,0,0}$0$}%
}}}}
\put(1208,-587){\makebox(0,0)[lb]{\smash{{\SetFigFont{9}{10.8}{\familydefault}{\mddefault}{\updefault}{\color[rgb]{0,0,0}\blue{$T_0M^p$}}%
}}}}
\end{picture}%

%% file: prolong-jet-1-transitive.pdf_t
\begin{picture}(0,0)%
\includegraphics{prolong-jet-1-transitive.pdf}%
\end{picture}%
\setlength{\unitlength}{4144sp}%
\begingroup\makeatletter\ifx\SetFigFont\undefined%
\gdef\SetFigFont#1#2#3#4#5{%
  \reset@font\fontsize{#1}{#2pt}%
  \fontfamily{#3}\fontseries{#4}\fontshape{#5}%
  \selectfont}%
\fi\endgroup%
\begin{picture}(5884,3178)(879,-2976)
\put(6376,-876){\makebox(0,0)[lb]{\smash{{\SetFigFont{9}{10.8}{\familydefault}{\mddefault}{\updefault}{\color[rgb]{0,0,0}\blue{$x_1$}}%
}}}}
\put(5881,-346){\makebox(0,0)[lb]{\smash{{\SetFigFont{9}{10.8}{\familydefault}{\mddefault}{\updefault}{\color[rgb]{0,0,0}\blue{$y_1$}}%
}}}}
\put(5806, 29){\makebox(0,0)[lb]{\smash{{\SetFigFont{9}{10.8}{\familydefault}{\mddefault}{\updefault}{\color[rgb]{0,0,0}\red{$\P_\infty^1$}}%
}}}}
\put(2007, 43){\makebox(0,0)[lb]{\smash{{\SetFigFont{9}{10.8}{\familydefault}{\mddefault}{\updefault}{\color[rgb]{0,0,0}\red{$\P_\infty^1$}}%
}}}}
\put(956,-230){\makebox(0,0)[lb]{\smash{{\SetFigFont{9}{10.8}{\familydefault}{\mddefault}{\updefault}{\color[rgb]{0,0,0}\blue{$J_{1,2}^1$}}%
}}}}
\put(2669,-735){\makebox(0,0)[lb]{\smash{{\SetFigFont{9}{10.8}{\familydefault}{\mddefault}{\updefault}{\color[rgb]{0,0,0}\green{$\vec{\bf v}^{(1)}$}}%
}}}}
\put(1994,-2903){\makebox(0,0)[lb]{\smash{{\SetFigFont{9}{10.8}{\familydefault}{\mddefault}{\updefault}{\color[rgb]{0,0,0}\blue{$0$}}%
}}}}
\put(1342,-2676){\makebox(0,0)[lb]{\smash{{\SetFigFont{9}{10.8}{\familydefault}{\mddefault}{\updefault}{\color[rgb]{0,0,0}\green{$\vec{\bf v}$}}%
}}}}
\put(2481,-2675){\makebox(0,0)[lb]{\smash{{\SetFigFont{9}{10.8}{\familydefault}{\mddefault}{\updefault}{\color[rgb]{0,0,0}\green{$\vec{\bf v}$}}%
}}}}
\put(2970,-2709){\makebox(0,0)[lb]{\smash{{\SetFigFont{9}{10.8}{\familydefault}{\mddefault}{\updefault}{\color[rgb]{0,0,0}\blue{$M^p$}}%
}}}}
\put(1332,-1076){\makebox(0,0)[lb]{\smash{{\SetFigFont{9}{10.8}{\familydefault}{\mddefault}{\updefault}{\color[rgb]{0,0,0}\green{$\vec{\bf v}^{(1)}$}}%
}}}}
\put(2070,-965){\makebox(0,0)[lb]{\smash{{\SetFigFont{9}{10.8}{\familydefault}{\mddefault}{\updefault}{\color[rgb]{0,0,0}\blue{$0$}}%
}}}}
\put(2061,-513){\makebox(0,0)[lb]{\smash{{\SetFigFont{9}{10.8}{\familydefault}{\mddefault}{\updefault}{\color[rgb]{0,0,0}\blue{$\R^2$}}%
}}}}
\put(5250,-1394){\makebox(0,0)[lb]{\smash{{\SetFigFont{9}{10.8}{\familydefault}{\mddefault}{\updefault}{\color[rgb]{0,0,0}\green{${\sf I}_2^{(1)}$}}%
}}}}
\put(6088,-1451){\makebox(0,0)[lb]{\smash{{\SetFigFont{9}{10.8}{\familydefault}{\mddefault}{\updefault}{\color[rgb]{0,0,0}\green{${\sf I}_1^{(1)}$}}%
}}}}
\put(5812,-1055){\makebox(0,0)[lb]{\smash{{\SetFigFont{9}{10.8}{\familydefault}{\mddefault}{\updefault}{\color[rgb]{0,0,0}\blue{$0$}}%
}}}}
\put(5815,-2912){\makebox(0,0)[lb]{\smash{{\SetFigFont{9}{10.8}{\familydefault}{\mddefault}{\updefault}{\color[rgb]{0,0,0}\blue{$0$}}%
}}}}
\end{picture}%

%% file: key-discovery-orbits.pdf_t
\begin{picture}(0,0)%
\includegraphics{key-discovery-orbits.pdf}%
\end{picture}%
\setlength{\unitlength}{4144sp}%
\begingroup\makeatletter\ifx\SetFigFont\undefined%
\gdef\SetFigFont#1#2#3#4#5{%
  \reset@font\fontsize{#1}{#2pt}%
  \fontfamily{#3}\fontseries{#4}\fontshape{#5}%
  \selectfont}%
\fi\endgroup%
\begin{picture}(4881,2844)(654,-2793)
\put(5520,-893){\makebox(0,0)[lb]{\smash{{\SetFigFont{9}{10.8}{\familydefault}{\mddefault}{\updefault}{\color[rgb]{0,0,0}$x_1,\!y_1$}%
}}}}
\put(4818,-473){\makebox(0,0)[lb]{\smash{{\SetFigFont{9}{10.8}{\familydefault}{\mddefault}{\updefault}{\color[rgb]{0,0,0}$x_2,\!y_2$}%
}}}}
\put(1687,-450){\makebox(0,0)[lb]{\smash{{\SetFigFont{7}{8.4}{\familydefault}{\mddefault}{\updefault}{\color[rgb]{0,0,0}$\R^2$}%
}}}}
\put(1466,-105){\makebox(0,0)[lb]{\smash{{\SetFigFont{7}{8.4}{\familydefault}{\mddefault}{\updefault}{\color[rgb]{0,0,0}$\R^2$}%
}}}}
\put(1611,-1775){\makebox(0,0)[lb]{\smash{{\SetFigFont{9}{10.8}{\familydefault}{\mddefault}{\updefault}{\color[rgb]{0,0,0}$\R^2$}%
}}}}
\put(1628,-1985){\makebox(0,0)[lb]{\smash{{\SetFigFont{9}{10.8}{\familydefault}{\mddefault}{\updefault}{\color[rgb]{0,0,0}$0$}%
}}}}
\put(713,-1787){\makebox(0,0)[lb]{\smash{{\SetFigFont{9}{10.8}{\familydefault}{\mddefault}{\updefault}{\color[rgb]{0,0,0}\blue{$J_{1,2}^1$}}%
}}}}
\put(691,-2718){\makebox(0,0)[lb]{\smash{{\SetFigFont{9}{10.8}{\familydefault}{\mddefault}{\updefault}{\color[rgb]{0,0,0}\blue{$M^p$}}%
}}}}
\put(1541,-2699){\makebox(0,0)[lb]{\smash{{\SetFigFont{9}{10.8}{\familydefault}{\mddefault}{\updefault}{\color[rgb]{0,0,0}$0$}%
}}}}
\put(977,-116){\makebox(0,0)[lb]{\smash{{\SetFigFont{9}{10.8}{\familydefault}{\mddefault}{\updefault}{\color[rgb]{0,0,0}\blue{$J_{1,2}^2$}}%
}}}}
\put(5514,-1411){\makebox(0,0)[lb]{\smash{{\SetFigFont{11}{13.2}{\familydefault}{\mddefault}{\updefault}{\color[rgb]{0,0,0}\red{$\Sigma_0^2$}}%
}}}}
\end{picture}%

%% file: translation-normalization.pdf_t
\begin{picture}(0,0)%
\includegraphics{translation-normalization.pdf}%
\end{picture}%
\setlength{\unitlength}{4144sp}%
\begingroup\makeatletter\ifx\SetFigFont\undefined%
\gdef\SetFigFont#1#2#3#4#5{%
  \reset@font\fontsize{#1}{#2pt}%
  \fontfamily{#3}\fontseries{#4}\fontshape{#5}%
  \selectfont}%
\fi\endgroup%
\begin{picture}(6341,1405)(661,-1943)
\put(4994,-688){\makebox(0,0)[lb]{\smash{{\SetFigFont{8}{9.6}{\familydefault}{\mddefault}{\updefault}{\color[rgb]{0,0,0}Normalization}%
}}}}
\put(5311,-846){\makebox(0,0)[lb]{\smash{{\SetFigFont{9}{10.8}{\familydefault}{\mddefault}{\updefault}{\color[rgb]{0,0,0}$\Phi_p$}%
}}}}
\put(2004,-1271){\makebox(0,0)[lb]{\smash{{\SetFigFont{9}{10.8}{\familydefault}{\mddefault}{\updefault}{\color[rgb]{0,0,.69}\blue{$p$}}%
}}}}
\put(826,-1411){\makebox(0,0)[lb]{\smash{{\SetFigFont{9}{10.8}{\familydefault}{\mddefault}{\updefault}{\color[rgb]{0,0,.69}\blue{$M$}}%
}}}}
\put(2335,-922){\makebox(0,0)[lb]{\smash{{\SetFigFont{9}{10.8}{\familydefault}{\mddefault}{\updefault}{\color[rgb]{0,0,0}Translation}%
}}}}
\put(2572,-1045){\makebox(0,0)[lb]{\smash{{\SetFigFont{8}{9.6}{\familydefault}{\mddefault}{\updefault}{\color[rgb]{0,0,0}$\tau_p$}%
}}}}
\put(4258,-1307){\makebox(0,0)[lb]{\smash{{\SetFigFont{9}{10.8}{\familydefault}{\mddefault}{\updefault}{\color[rgb]{0,0,.69}\blue{$0$}}%
}}}}
\put(3069,-1414){\makebox(0,0)[lb]{\smash{{\SetFigFont{9}{10.8}{\familydefault}{\mddefault}{\updefault}{\color[rgb]{0,0,.69}\blue{$M^p$}}%
}}}}
\put(6572,-769){\makebox(0,0)[lb]{\smash{{\SetFigFont{9}{10.8}{\familydefault}{\mddefault}{\updefault}{\color[rgb]{0,0,.69}\blue{$v$}}%
}}}}
\put(4617,-1108){\makebox(0,0)[lb]{\smash{{\SetFigFont{9}{10.8}{\familydefault}{\mddefault}{\updefault}{\color[rgb]{0,0,.69}\blue{$z,\overline{z},u$}}%
}}}}
\put(6508,-1307){\makebox(0,0)[lb]{\smash{{\SetFigFont{9}{10.8}{\familydefault}{\mddefault}{\updefault}{\color[rgb]{0,0,.69}\blue{$0$}}%
}}}}
\put(6871,-1112){\makebox(0,0)[lb]{\smash{{\SetFigFont{9}{10.8}{\familydefault}{\mddefault}{\updefault}{\color[rgb]{0,0,.69}\blue{$z,\overline{z},u$}}%
}}}}
\put(4322,-769){\makebox(0,0)[lb]{\smash{{\SetFigFont{9}{10.8}{\familydefault}{\mddefault}{\updefault}{\color[rgb]{0,0,.69}\blue{$v$}}%
}}}}
\put(5416,-1879){\makebox(0,0)[lb]{\smash{{\SetFigFont{9}{10.8}{\familydefault}{\mddefault}{\updefault}{\color[rgb]{0,0,0}$N^p$}%
}}}}
\put(3399,-1797){\makebox(0,0)[lb]{\smash{{\SetFigFont{8}{9.6}{\familydefault}{\mddefault}{\updefault}{\color[rgb]{0,0,0}$\Phi_p\circ\tau_p=:\varphi$}%
}}}}
\end{picture}%

%% file: f-0-1-2-g-0-1.pdf_t
\begin{picture}(0,0)%
\includegraphics{f-0-1-2-g-0-1.pdf}%
\end{picture}%
\setlength{\unitlength}{4144sp}%
\begingroup\makeatletter\ifx\SetFigFont\undefined%
\gdef\SetFigFont#1#2#3#4#5{%
  \reset@font\fontsize{#1}{#2pt}%
  \fontfamily{#3}\fontseries{#4}\fontshape{#5}%
  \selectfont}%
\fi\endgroup%
\begin{picture}(6330,1026)(353,-803)
\put(2654,-445){\makebox(0,0)[lb]{\smash{{\SetFigFont{9}{10.8}{\familydefault}{\mddefault}{\updefault}{\color[rgb]{0,0,0}$5$}%
}}}}
\put(2189,-444){\makebox(0,0)[lb]{\smash{{\SetFigFont{9}{10.8}{\familydefault}{\mddefault}{\updefault}{\color[rgb]{0,0,0}$4$}%
}}}}
\put(1759,-445){\makebox(0,0)[lb]{\smash{{\SetFigFont{9}{10.8}{\familydefault}{\mddefault}{\updefault}{\color[rgb]{0,0,0}$3$}%
}}}}
\put(1310,-444){\makebox(0,0)[lb]{\smash{{\SetFigFont{9}{10.8}{\familydefault}{\mddefault}{\updefault}{\color[rgb]{0,0,0}$2$}%
}}}}
\put(850,-445){\makebox(0,0)[lb]{\smash{{\SetFigFont{9}{10.8}{\familydefault}{\mddefault}{\updefault}{\color[rgb]{0,0,0}$1$}%
}}}}
\put(404,-444){\makebox(0,0)[lb]{\smash{{\SetFigFont{9}{10.8}{\familydefault}{\mddefault}{\updefault}{\color[rgb]{0,0,0}$0$}%
}}}}
\put(6254,-445){\makebox(0,0)[lb]{\smash{{\SetFigFont{9}{10.8}{\familydefault}{\mddefault}{\updefault}{\color[rgb]{0,0,0}$5$}%
}}}}
\put(5789,-444){\makebox(0,0)[lb]{\smash{{\SetFigFont{9}{10.8}{\familydefault}{\mddefault}{\updefault}{\color[rgb]{0,0,0}$4$}%
}}}}
\put(5359,-445){\makebox(0,0)[lb]{\smash{{\SetFigFont{9}{10.8}{\familydefault}{\mddefault}{\updefault}{\color[rgb]{0,0,0}$3$}%
}}}}
\put(4910,-444){\makebox(0,0)[lb]{\smash{{\SetFigFont{9}{10.8}{\familydefault}{\mddefault}{\updefault}{\color[rgb]{0,0,0}$2$}%
}}}}
\put(4450,-445){\makebox(0,0)[lb]{\smash{{\SetFigFont{9}{10.8}{\familydefault}{\mddefault}{\updefault}{\color[rgb]{0,0,0}$1$}%
}}}}
\put(4004,-444){\makebox(0,0)[lb]{\smash{{\SetFigFont{9}{10.8}{\familydefault}{\mddefault}{\updefault}{\color[rgb]{0,0,0}$0$}%
}}}}
\put(368,-734){\makebox(0,0)[lb]{\smash{{\SetFigFont{9}{10.8}{\familydefault}{\mddefault}{\updefault}{\color[rgb]{0,0,0}$\big\{f_0,\,f_1,\,f_2\big\}$}%
}}}}
\put(3953,-719){\makebox(0,0)[lb]{\smash{{\SetFigFont{9}{10.8}{\familydefault}{\mddefault}{\updefault}{\color[rgb]{0,0,0}$\big\{g_0,\,g_1\big\}$}%
}}}}
\put(4098, 76){\makebox(0,0)[lb]{\smash{{\SetFigFont{9}{10.8}{\familydefault}{\mddefault}{\updefault}{\color[rgb]{0,0,0}$g$}%
}}}}
\put(418, -2){\makebox(0,0)[lb]{\smash{{\SetFigFont{9}{10.8}{\familydefault}{\mddefault}{\updefault}{\color[rgb]{0,0,0}$f$}%
}}}}
\end{picture}%